\pdfoutput=1
\documentclass[reqno]{amsart}
\usepackage{enumerate}
\usepackage[mathscr]{euscript}
\usepackage{xcolor}
\usepackage{array}
\usepackage{amsfonts}
\usepackage{amsmath}
\usepackage{amssymb}
\usepackage{mathtools}
\usepackage{graphicx}
\usepackage{bm}
\usepackage{enumitem}
\usepackage{caption} 
\usepackage{color}
\usepackage{bookmark}
\usepackage[square,numbers,sort&compress]{natbib}
\usepackage{hyperref}
\hypersetup{colorlinks=true,linkcolor=blue,citecolor=red}
\allowdisplaybreaks
\def\Xint#1{\mathchoice
{\XXint\displaystyle\textstyle{#1}}%
{\XXint\textstyle\scriptstyle{#1}}%
{\XXint\scriptstyle\scriptscriptstyle{#1}}%
{\XXint\scriptscriptstyle\scriptscriptstyle{#1}}%
\!\int}
\def\XXint#1#2#3{{\setbox0=\hbox{$#1{#2#3}{\int}$ }
\vcenter{\hbox{$#2#3$ }}\kern-.6\wd0}}

\def\dashint{\Xint-}

\newtheorem{theorem}{Theorem}[section]
\newtheorem{lemma}[theorem]{Lemma}

\newtheorem{remark}[theorem]{Remark}
\theoremstyle{definition}
\newtheorem{definition}[theorem]{Definition}

\numberwithin{equation}{section}

\newcommand{ \mr }{ \mathbb{R} }

\newcommand{\iints}[1]{{\int\hspace{-0.28cm}\int_{#1}}}
\newcommand{\iintss}{{\int\hspace{-0.28cm}\int}}
\newcommand{ \miints }{{\iintss\hspace{-0.56cm} -\hspace{-0.15cm}-}}
\newcommand{\miint}[1]{{\miints_{\hspace{-0.13cm}#1}}}

\begin{document}
\title[A general class of parabolic double phase systems]{Gradient higher integrability for a general class of degenerate parabolic double phase systems}

\author{Bogi Kim}\address{Department of Mathematics, Kyungpook National University, Daegu, 41566, Republic of Korea} \email{rlaqhrl4@knu.ac.kr}
\author{Jehan Oh}\address{Department of Mathematics, Kyungpook National University, Daegu, 41566, Republic of Korea} \email{jehan.oh@knu.ac.kr}
\author{Abhrojyoti Sen}\address{Goethe-Universit\"{a}t Frankfurt, Institut f\"{u}r Mathematik, Robert-Mayer-Str. 10,D-60629 Frankfurt, Germany} \email{sen@math.uni-frankfurt.de}

\subjclass{Primary 35B65; Secondary 35D30, 35K40, 35K55, 35K65}
\date{\today}
\keywords{degenerate parabolic systems, double phase problems, higher integrability, Orlicz growth, interpolation}
\thanks{Bogi Kim is supported by the National Research Foundation of Korea (NRF) grant funded by the Korea government [Grant No. RS-2025-25426375].
Jehan Oh is supported by the National Research Foundation of Korea (NRF) grant funded by the Korea government [Grant Nos. RS-2025-00555316 and RS-2025-25415411].
Abhrojyoti Sen is supported by postdoctoral research grants from the Alexander von Humboldt Foundation, Germany.}

\begin{abstract}
We establish quantitative local higher integrability estimates for weak solutions of degenerate parabolic systems with generalized double phase growth, whose prototype is
\[
 u_t-\operatorname{div}\left(\frac{G'(|Du|)}{|Du|}Du+a(x,t)\frac{H'(|Du|)}{|Du|}Du\right)=0\qquad\text{in }\Omega_T.
\]
Here \(G\) and \(H\) are Young functions satisfying suitable structural assumptions. The nonnegative coefficient \(a\) is allowed to vanish and has a concave modulus of continuity \(\omega\) with respect to the parabolic distance. For \(\Psi(z,\varkappa)=G(\varkappa)+a(z)H(\varkappa)\), we prove that solutions in \(L^\infty(0,T;L^s(\Omega,\mathbb{R}^N))\), \(2\le s<\infty\), satisfy \(\Psi(\cdot,|Du|)\in L^{1+\varepsilon}_{\mathrm{loc}}\) for some \(\varepsilon>0\), provided
\[
 \limsup_{\rho \to 0+}\omega(\rho)\frac{H(\rho^{-\frac{n+s}{s}})}{G(\rho^{-\frac{n+s}{s}})}<\infty.
\]
For bounded solutions, the corresponding condition involves \(H(\rho^{-1})/G(\rho^{-1})\). These conditions extend the known interpolative gap bounds from power growth to general growth functions and continuity moduli. In particular, for \(G(\varkappa)=\varkappa^p\), \(H(\varkappa)=\varkappa^q\), and \(\omega(\rho)\simeq\rho^\alpha\), they recover \(q\le p+s\alpha/(n+s)\) and \(q\le p+\alpha\), respectively. The proof combines an intrinsic stopping-time argument with inverse-function estimates to derive Poincar\'e and reverse H\"{o}lder inequalities adapted to the general growth phases.
\end{abstract}
\maketitle

\section{Introduction and main results}\label{section 1}
We investigate parabolic double phase systems of the form
\begin{equation}    \label{eq : main equation}
    u_t - \operatorname{div} \mathcal{A}(z,Du)=0 \qquad \text{in} \ \, \Omega_T\coloneq \Omega\times (0,T),
\end{equation}
where $\Omega \subset \mr^n$ is a bounded open set with $n \geq 2$ and $\mathcal{A}:\Omega_T \times \mr^{Nn}\rightarrow \mr^{Nn}$ with $N \geq 1$ is a Carath\'{e}odory vector field satisfying the following coercivity and growth bounds: there exist constants $0<\nu\leq L <\infty$ such that
\begin{equation}    \label{cond : double phase bounded condition of integrand}
    \mathcal{A}(z,\xi)\cdot \xi\geq \nu \left[G(|\xi|)+a(z)H(|\xi|)\right] \quad \text{and} \quad |\mathcal{A}(z,\xi)| \leq L \frac{G(|\xi|)+a(z)H(|\xi|)}{|\xi|}
\end{equation}
for all $z \in \Omega_T$ and $\xi \in \mr^{Nn}$. Here $G$ and $H$ are Young functions with $G\prec H$ in the sense of Definition \ref{def : relation of Phi1 and Phi2} below, which means that $H$ grows faster than $G$, and $a:\Omega_T\to[0,\infty)$ is a bounded and continuous modulating coefficient. Throughout the paper we write
\begin{equation}\label{def : definition of Psi}
\Psi(z,\varkappa)\coloneq G(\varkappa)+a(z)H(\varkappa)\qquad\text{for } z\in\Omega_T \text{ and } \varkappa\geq 0.
\end{equation}
The prototype of \eqref{eq : main equation} is obtained for $G(\varkappa)=\varkappa^p$ and $H(\varkappa)=\varkappa^q$ with $2\leq p<q<\infty$, namely the parabolic double phase system
\begin{equation}\label{eq : model parabolic double phase}
u_t-\operatorname{div}\left(|Du|^{p-2}Du+a(z)|Du|^{q-2}Du\right)=0\qquad\text{in }\Omega_T,
\end{equation}
which is the evolutionary counterpart of the elliptic double phase functional
\begin{equation}\label{eq : elliptic double phase functional}
w\mapsto\int_\Omega\left[|Dw|^p+a(x)|Dw|^q\right]dx.
\end{equation}
The functional \eqref{eq : elliptic double phase functional} was introduced by Zhikov \cite{Zhikov1986,Zhikov1995} in the study of homogenization and of the Lavrentiev phenomenon, and it is a model example of the functionals with $(p,q)$-growth considered by Marcellini \cite{Marcellini1989}. The regularity of its minimizers depends on the interplay between the gap $q-p$ and the regularity of the coefficient $a$. If $a\in C^{0,\alpha}(\Omega)$ for some $\alpha\in(0,1]$, then local minimizers have locally H\"{o}lder continuous gradients provided that
$$
\frac{q}{p}\leq 1+\frac{\alpha}{n},
$$
see Colombo and Mingione \cite{Colombo2015} and Baroni, Colombo and Mingione \cite{Baroni2018}, and this gap bound is sharp in view of the counterexamples in \cite{Esposito2004}. For bounded minimizers the gap bound improves to $q\leq p+\alpha$ \cite{Colombo2015a,Baroni2018}, and for minimizers with $u\in L^\gamma_{\mathrm{loc}}(\Omega)$ an interpolating gap bound of the form $q\leq p+\frac{\gamma\alpha}{n+\gamma}$ was obtained by Ok \cite{Ok2020}, see also De Filippis and Mingione \cite{DeFilippisMingione2021AMP} for interpolative gap bounds for general nonautonomous integrals. Generalized double phase functionals, in which the powers are replaced by Young functions $G\prec H$, were studied by Byun and the second author \cite{Byun2020}. There the gap bound is expressed through the modulus of continuity of $a$ and the quotient $H/G$, and it reduces to the gap bounds above in the power case. This Orlicz double phase setting was developed further in \cite{Baasandorj2020,BaasandorjByun2023ASNS}, and it is also covered by the regularity theory in generalized Orlicz spaces of \cite{Harjulehto2019,Haestoe2022a}. Borderline double phase functionals such as $|Dw|^p+a(x)|Dw|^p\log(1+|Dw|)$, which belong to the generalized setting, were treated in \cite{BaroniColomboMingione2016SPMJ}.

The parabolic theory is less complete. Parabolic systems with $(p,q)$-growth were studied by B\"{o}gelein, Duzaar and Marcellini \cite{BogeleinDuzaarMarcellini2013ARMA}. For the parabolic double phase system \eqref{eq : model parabolic double phase}, Kim, Kinnunen and Moring \cite{2023_Gradient_Higher_Integrability_for_Degenerate_Parabolic_Double-Phase_Systems} established the local higher integrability of the gradient under the assumptions
\begin{equation}\label{eq : KKM gap bound}
u\in C(0,T;L^2(\Omega,\mr^N)),\qquad a\in C^{\alpha,\frac{\alpha}{2}}(\Omega_T)\qquad\text{and}\qquad q\leq p+\frac{2\alpha}{n+2},
\end{equation}
where $C^{\alpha,\frac{\alpha}{2}}(\Omega_T)$ denotes the parabolic H\"{o}lder space. This result has been extended to the natural function space of weak solutions by means of the parabolic Lipschitz truncation \cite{KIM2025110738}, to singular systems \cite{KS24}, to multi-phase problems \cite{KimOh2024JDE,S25} and to systems with two modulating coefficients \cite{OhSen2026CVPDE}, and the H\"{o}lder regularity of weak solutions was obtained in \cite{KIM2025113231}. Moreover, the first two authors \cite{kim2025bounded} proved the higher integrability for bounded weak solutions of \eqref{eq : model parabolic double phase} under the gap bound $q\leq p+\alpha$, and for weak solutions with $u\in L^\infty(0,T;L^s(\Omega,\mr^N))$ and $s\geq 2$ under the gap bound
\begin{equation}\label{eq : KO interpolative gap bound}
q\leq p+\frac{s\alpha}{n+s},
\end{equation}
which interpolates between \eqref{eq : KKM gap bound} and the bounded case. The singular counterpart of this interpolation is treated in \cite{kim2026interpolative}.

The higher integrability of gradients goes back to Gehring \cite{Gehring1973} and Giaquinta and Modica \cite{GiaquintaModica1979} in the elliptic case. For degenerate and singular parabolic systems of $p$-Laplacian type the result is due to Kinnunen and Lewis \cite{Kinnunen2000}, whose argument relies on the intrinsic geometry of DiBenedetto \cite{1993_Degenerate_parabolic_equations_DiBenedetto}, and it has been extended to parabolic systems with Orlicz growth by H\"{a}st\"{o} and Ok \cite{Hasto_2021}, see also \cite{OSS24}. The double phase structure requires, in addition, two different intrinsic geometries, one for each phase, as introduced in \cite{2023_Gradient_Higher_Integrability_for_Degenerate_Parabolic_Double-Phase_Systems}.

The aim of this paper is to establish the gradient higher integrability for \eqref{eq : main equation} with general Young functions $G\prec H$, thereby extending the results of \cite{2023_Gradient_Higher_Integrability_for_Degenerate_Parabolic_Double-Phase_Systems,kim2025bounded} from the power case to the Orlicz setting of \cite{Byun2020}. The gap bounds \eqref{eq : KKM gap bound} and \eqref{eq : KO interpolative gap bound} are replaced by conditions on the modulus of continuity of $a$ and on the quotient $H/G$ along the scales $\rho^{-1}$ and $\rho^{-\frac{n+s}{s}}$, see \eqref{cond : main assumption with infty}, \eqref{cond : main assumption with 2} and \eqref{cond : main assumption with s} below.

The recent work of H\"ast\"o and Ok~\cite{HastoOk2026}
establishes gradient higher integrability for parabolic
systems with generalized Orlicz growth under an intrinsic
parabolic $(A1)$ condition. This condition requires the
growth density $\varphi(z,\cdot)$ to be uniformly
comparable at different points of every admissible
intrinsic parabolic cylinder, for gradient levels in a
prescribed energy range, see
\cite[Definition~3.1]{HastoOk2026}.
Their Theorem~1.6 treats both singular and degenerate
growth in the natural energy class, under the additional
restriction
\[
q<p+2\min\left\{\frac{p}{n},1\right\}.
\]
Our approach is different in that, instead of assuming
intrinsic $(A1)$, we exploit additional information on the
solution to obtain higher integrability under
solution-dependent balance conditions. More precisely, for
generalized double phase growth
\[
\Psi(z,\varkappa)=G(\varkappa)+a(z)H(\varkappa),
\]
uniform-in-time $L^s$ control of the solution leads to
a balance between the modulus of continuity $\omega$
of $a$ and the quotient $H/G$ at the scale
$\rho^{-\frac{n+s}{s}}$, while boundedness of the solution
allows the scale $\rho^{-1}$. The principal contribution of our results is to
extend the known solution-dependent interpolative
gap bounds to general $(G,H)$-growth and concave
continuity moduli, complementing the general
framework of~\cite{HastoOk2026}. 
We also note that H\"ast\"o and Ok
\cite[Remark~1.7]{HastoOk2026} identify the range
\[
q\ge p+2\min\left\{\frac{p}{n},1\right\}
\]
as an open problem for their general nonautonomous theory. Although our results do not resolve this problem in the natural energy class, the
solution-dependent balance conditions considered here allow certain double phase systems with growth indices in this range. Indeed, in the
power-growth case with a Lipschitz continuous coefficient, the bounds
\[
q-p\le\frac{s}{n+s}
\qquad\text{and}\qquad q-p\le1
\]
can reach or exceed the above threshold when the spatial dimension is sufficiently large. This illustrates the complementary nature of
the two approaches.

We now describe our assumptions in detail. Let $G, H : [0,\infty) \to [0,\infty)$ be Young functions with $G,\,H\in C^1([0,\infty))$, $G(1)=1=H(1)$ and $G\prec H$, where $\prec$ is defined in Definition \ref{def : relation of Phi1 and Phi2}. Assume that $G,\, H\in \Delta_2\cap \nabla_2$, where $\nabla_2$ and $\Delta_2$ are defined in Definition \ref{def : nabla and Delta}. By Lemma \ref{lem : Phi and Phi'} we may fix constants $\sigma_G=\sigma_G(\Delta_2(G))>1$ and $\sigma_H=\sigma_H(\Delta_2(H))>1$, for instance $\sigma_G=\Delta_2(G)$ and $\sigma_H=\Delta_2(H)$, such that
\begin{equation}\label{cond : pointwise bound of Phi'}
    1\leq \frac{\varkappa G'(\varkappa )}{G(\varkappa )}\leq\sigma_G\quad\text{and}\quad 1\leq\frac{\varkappa H'(\varkappa )}{H(\varkappa )}\leq \sigma_H
\end{equation}
for any $\varkappa\in(0,\infty)$. Since $G,\,H\in C^1((0,\infty))$, writing
$$
\frac{d}{d\varkappa}\log\frac{\Phi(\varkappa)}{\varkappa^{\varsigma}}
=\frac{1}{\varkappa}\left(\frac{\varkappa\Phi'(\varkappa)}{\Phi(\varkappa)}-\varsigma\right),
\qquad \Phi\in\{G,H\},
$$
we see that \eqref{cond : pointwise bound of Phi'} is equivalent to
\begin{equation}\label{cond : monotonicity of Phi over powers}
    (0,\infty)\ni\varkappa\mapsto\frac{\Phi(\varkappa)}{\varkappa}
    \ \text{is non-decreasing}
    \enspace\text{and}\enspace
    (0,\infty)\ni\varkappa\mapsto\frac{\Phi(\varkappa)}{\varkappa^{\sigma_\Phi}}
    \ \text{is non-increasing},
\end{equation}
where $\sigma_\Phi$ denotes $\sigma_G$ or $\sigma_H$ according as $\Phi=G$ or $\Phi=H$. We stress that \eqref{cond : monotonicity of Phi over powers} is a genuine monotonicity and not an almost monotonicity. This is what allows us to pass to the inverse functions in Lemma \ref{lem : properties of Xi}. Together with $\Phi(1)=1$ it also yields
\begin{equation}\label{cond : two sided power bound of Phi}
    \Phi(\varkappa)\leq \varkappa^{\sigma_\Phi}\ \text{ for }\varkappa\geq1
    \qquad\text{and}\qquad
    \Phi(\varkappa)\leq 1\ \text{ for }0<\varkappa\leq 1 .
\end{equation}
We assume that \begin{equation}\label{cond : relation of sigma G and sigma H}\sigma_G\leq \sigma_H.\end{equation}
By \cite[Lemma 4.1]{1961Rutickii}, we have
$$
\frac{\varkappa (G^*(\varkappa ))'}{G^*(\varkappa )}\geq \frac{\sigma_G}{\sigma_G-1}\quad \text{and}\quad\frac{\varkappa (H^*(\varkappa ))'}{H^*(\varkappa )}\geq\frac{\sigma_H}{\sigma_H-1}
$$
and hence we have from \cite[(2.3)]{OSS24} that for $0<c_1<1<c_2$ and $\varkappa,\,\vartheta\in (0,\infty)$
\begin{equation}\label{cond : H* ct and H t with c<1}
\Phi^*(c_1\vartheta)\leq c_1^{\frac{\sigma_\Phi}{\sigma_\Phi-1}} \Phi^*(\vartheta)\quad\text{and}\quad \Phi(c_2\varkappa)\leq c_2^{\sigma_\Phi} \Phi(\varkappa),
\end{equation}
where $\Phi$ denotes either the function $G$ or $H$. The second inequality in \eqref{cond : H* ct and H t with c<1} is also immediate from \eqref{cond : monotonicity of Phi over powers}, and the same argument shows that it remains valid for $0<c_2\leq1$ with the inequality sign reversed. Moreover, by \cite[Corollary 2.4.11]{Harjulehto2019}, the $\nabla_2$ condition for $G$ and $H$ implies that there exist $i_G=i_G(\nabla_2(G))>1$ and $i_H=i_H(\nabla_2(H))>1$ such that $\frac{G(\varkappa)}{\varkappa^{i_G}}$ and $\frac{H(\varkappa)}{\varkappa^{i_H}}$ are almost increasing for $\varkappa\in(0,\infty)$ with constants $\tilde{L}_G=\tilde{L}_G(\nabla_2(G))\geq 1$ and $\tilde{L}_H=\tilde{L}_H(\nabla_2(H))\geq 1$, respectively. Clearly, $i_G\leq \sigma_G$. Without loss of generality, we may assume that
$$
i_H\leq i_G.
$$
Consequently, $\frac{G(\varkappa)}{\varkappa^{i_H}}$ is almost increasing with constant $\tilde{L}\coloneq \max\{\tilde{L}_G,\tilde{L}_H\}$. For $\theta\in(\frac{1}{i_H},1]$, since $\tilde{L}\geq 1$ and $\theta i_H>1$, we see that $\frac{G^\theta(\varkappa)}{\varkappa^{\theta i_H}}$ and $\frac{H^{\theta}(\varkappa)}{\varkappa^{\theta i_H}}$ are almost increasing with constant $\tilde{L}$ and hence 
\begin{equation}\label{cond : ainc_1 of G and H}
    \frac{G^{\theta}(\varkappa)}{\varkappa}\text{ and }\frac{H^{\theta}(\varkappa)}{\varkappa} \text{ are almost increasing with constant }\tilde{L}.
\end{equation}
Now, for $\tilde{G}(\varkappa)\coloneq G(\varkappa^{\frac{1}{2}})$ $(\varkappa\in [0,\infty))$, we suppose that $i_H\geq 2$ and that $\tilde{G}$ is convex. In the power case $G(\varkappa)=\varkappa^p$ both requirements amount to $p\geq 2$, see Remark \ref{rem : special cases}. Since $\tilde{G}$ is linear when $G(\varkappa)=\varkappa^2$, we do not require $\tilde{G}$ to be a Young function. The properties of $\tilde{G}$ used below are collected in Lemmas \ref{lem : quadratic ainc} and \ref{lem : tilde G conjugate}.

Suppose that the modulating coefficient $a : \Omega_T \to [0,\infty)$ is bounded and has a modulus of continuity $\omega:[0,\infty)\rightarrow[0,\infty)$ with respect to the parabolic distance, that is,
$$
|a(z_1)-a(z_2)|\leq \omega\bigl(d_p(z_1,z_2)\bigr)\qquad\text{for all }z_1,\,z_2\in\Omega_T,
$$
where
$$
d_p(z_1,z_2)\coloneq \max\left\{|x_1-x_2|,|t_1-t_2|^{\frac{1}{2}}\right\}\qquad\text{for }z_i=(x_i,t_i),\ i=1,2.
$$
We assume that $\omega$ is concave with $\omega(0)= 0$. Then, we observe that for $c\geq 1$ and $\rho>0$,
\begin{equation}\label{cond : property of concave function}
\omega(c\rho)\leq c\omega(\rho).
\end{equation}
Indeed, for $c\geq 1$ and $\rho>0$, $\omega(\rho)=\omega(\frac{c\rho}{c})=\omega(\frac{1}{c}c\rho+(1-\frac{1}{c})0)\geq \frac{1}{c}\omega(c\rho)$. Now, we shall consider the following conditions:
\begin{equation}    \label{cond : main assumption with infty}
    u\in L^\infty(\Omega_T, \mr^N) \quad \text{and} \quad \limsup_{\rho \to 0+} \omega(\rho) \frac{H(\rho^{-1})}{G(\rho^{-1})} < \infty,
\end{equation}
or
\begin{equation}    \label{cond : main assumption with 2}
u\in L^\infty(0,T;L^2(\Omega, \mr^N)) \quad \text{and} \quad \limsup_{\rho \to 0+} \omega(\rho) \frac{H(\rho^{-\frac{n+2}{2}})}{G(\rho^{-\frac{n+2}{2}})} < \infty.
\end{equation}
Also, we shall consider the following interpolative conditions: for $s \geq 2$,
\begin{equation}    \label{cond : main assumption with s}
u\in L^\infty(0,T;L^s(\Omega, \mr^N)) \quad \text{and} \quad \limsup_{\rho \to 0+} \omega(\rho) \frac{H(\rho^{-\frac{n+s}{s}})}{G(\rho^{-\frac{n+s}{s}})} < \infty.
\end{equation}
Since \eqref{cond : main assumption with 2} is the case $s=2$ of \eqref{cond : main assumption with s}, we will carry out the proofs under either \eqref{cond : main assumption with infty} or \eqref{cond : main assumption with s}. Observe that the second condition in \eqref{cond : main assumption with infty} holds if and only if there exists $\varpi_b>0$ such that
\begin{equation}\label{cond : property of infty condition}
\omega(\rho)\frac{H(\rho^{-1})}{G(\rho^{-1})}\leq\varpi_b\qquad\text{for every }0<\rho\leq 1,
\end{equation}
and the second condition in \eqref{cond : main assumption with s} holds if and only if there exists $\varpi_s>0$ such that
\begin{equation}\label{cond : property of s}
\omega(\rho)\frac{H(\rho^{-\frac{n+s}{s}})}{G(\rho^{-\frac{n+s}{s}})}\leq \varpi_s\qquad\text{for every }0<\rho\leq 1.
\end{equation}

\begin{definition}\label{def : weak solution}
    A function $u:\Omega_T\rightarrow \mr^N$ with
    $$
    u\in C(0,T;L^2(\Omega, \mr^N))\cap L^{\sigma_H}(0,T;W^{1,\sigma_H}(\Omega, \mr^N))
    $$
    is a weak solution to \eqref{eq : main equation} if
    $$
    \iints{\Omega_T} -u\cdot \varphi_t +\mathcal{A}(z,Du)\cdot D\varphi\, dz=0 
    $$
    holds for every $\varphi\in C_0^\infty(\Omega_T, \mr^N)$.
\end{definition}

\begin{remark}\label{rem : function space}
In \cite{2023_Gradient_Higher_Integrability_for_Degenerate_Parabolic_Double-Phase_Systems}, weak solutions to the model case \eqref{eq : model parabolic double phase} are taken in $C(0,T;L^2(\Omega,\mr^N))\cap L^q(0,T;W^{1,q}(\Omega,\mr^N))$, and the exponent $q$ corresponds to $\sigma_H$ in our setting. The space $L^{\sigma_H}(0,T;W^{1,\sigma_H}(\Omega,\mr^N))$ guarantees that $\Psi(\cdot,|Du|)\in L^1(\Omega_T)$ and that the integral in Definition \ref{def : weak solution} is finite, since $|\mathcal{A}(z,Du)|\leq L(1+\|a\|_{L^\infty(\Omega_T)})(1+|Du|^{\sigma_H-1})$ by \eqref{cond : double phase bounded condition of integrand}, \eqref{cond : monotonicity of Phi over powers} and \eqref{cond : two sided power bound of Phi}. In the proofs below it is used only to guarantee that certain integrals are finite, namely in Lemmas \ref{lem : Caccioppoli inequality}, \ref{lem : Slicewise Sobolev inequality} and \ref{lem : Intrinsic interpolation estimate}. As in \cite{KIM2025110738}, we expect that it can be replaced by the natural energy space $C(0,T;L^2(\Omega,\mr^N))\cap L^1(0,T;W^{1,1}(\Omega,\mr^N))$ with $\Psi(\cdot,|Du|)\in L^1(\Omega_T)$ by means of a Lipschitz truncation argument adapted to the Orlicz setting, but we do not pursue this here.
\end{remark}
For $z_0=(x_0,t_0)\in\mr^{n+1}$ and $\varrho,\ell>0$ we write $I_\ell(t_0)\coloneq(t_0-\ell,t_0+\ell)$, $Q_{\varrho,\ell}(z_0)\coloneq B_\varrho(x_0)\times I_\ell(t_0)$ and $Q_\varrho(z_0)\coloneq Q_{\varrho,\varrho^2}(z_0)$, where $B_\varrho(x_0)$ is the open ball in $\mr^n$. We also write
$$
\operatorname{data}_b=\operatorname{data}_b(n,N,\nu,L,\varpi_b,\Delta_2(G),\nabla_2(G),\Delta_2(H),\nabla_2(H),\|u\|_{L^\infty(\Omega_T)}).
$$

\begin{theorem}\label{thm : main theorem with infty}
    Assume that \eqref{cond : main assumption with infty} is satisfied and let $u$ be a weak solution to \eqref{eq : main equation} with \eqref{cond : double phase bounded condition of integrand}. Then there exist constants $\varepsilon_0=\varepsilon_0(\operatorname{data}_b)>0$ and $c=c(\operatorname{data}_b,$ $\operatorname{diam}(\Omega), \|a\|_{L^\infty(\Omega_T)})>1$ such that
    $$
    \begin{aligned}
        &\miint{Q_r(z_0)} \Psi(z,|Du|)^{1+\varepsilon}\, dz\leq c\left(\miint{Q_{2r}(z_0)} \Psi(z,|Du|)\,dz\right)^{1+\frac{\varepsilon \sigma_H}{2}}+c
    \end{aligned}
    $$
    for every $Q_{2r}(z_0)\subset \Omega_T$ and $\varepsilon\in (0,\varepsilon_0)$.
\end{theorem}

We write
$$
\operatorname{data}_s=\operatorname{data}_s(n,N,s,\nu,L,\varpi_s,\Delta_2(G),\nabla_2(G),\Delta_2(H),\nabla_2(H),\|u\|_{L^\infty(0,T;L^s(\Omega))}).
$$
\begin{theorem}\label{thm : main theorem with s}
    Assume that \eqref{cond : main assumption with s} is satisfied and let $u$ be a weak solution to \eqref{eq : main equation} with \eqref{cond : double phase bounded condition of integrand}. Then there exist constants $\varepsilon_0=\varepsilon_0(\operatorname{data}_s)>0$ and $c=c(\operatorname{data}_s,$ $\operatorname{diam}(\Omega), \|a\|_{L^\infty(\Omega_T)})>1$ such that
    $$
    \begin{aligned}
        &\miint{Q_r(z_0)} \Psi(z,|Du|)^{1+\varepsilon}\, dz\leq c\left(\miint{Q_{2r}(z_0)} \Psi(z,|Du|)\,dz\right)^{1+\frac{\varepsilon \sigma_H}{2}}+c
    \end{aligned}
    $$
    for every $Q_{2r}(z_0)\subset \Omega_T$ and $\varepsilon\in (0,\varepsilon_0)$.
\end{theorem}

We remark that in the model case $G(\varkappa)=\varkappa^p$ and $H(\varkappa)=\varkappa^q$ with $a\in C^{\alpha,\frac{\alpha}{2}}(\Omega_T)$, so that $\omega(\rho)=[a]_\alpha\rho^\alpha$ and $\frac{H(\varkappa)}{G(\varkappa)}=\varkappa^{q-p}$, the conditions on $\omega$ in \eqref{cond : main assumption with infty}, \eqref{cond : main assumption with 2} and \eqref{cond : main assumption with s} reduce to the finiteness of
$$
\limsup_{\rho\to 0+}\rho^{\alpha-(q-p)},\qquad \limsup_{\rho\to 0+}\rho^{\alpha-\frac{(n+2)(q-p)}{2}}\qquad\text{and}\qquad \limsup_{\rho\to 0+}\rho^{\alpha-\frac{(n+s)(q-p)}{s}},
$$
that is,
$$
q\leq p+\alpha,\qquad q\leq p+\frac{2\alpha}{n+2}\qquad\text{and}\qquad q\leq p+\frac{s\alpha}{n+s},
$$
respectively. Hence Theorem \ref{thm : main theorem with infty} recovers the result for bounded weak solutions in \cite[Theorem 1.2]{kim2025bounded}, while Theorem \ref{thm : main theorem with s} recovers \cite[Theorem 1.1]{2023_Gradient_Higher_Integrability_for_Degenerate_Parabolic_Double-Phase_Systems} for $s=2$ and \cite[Theorem 1.3]{kim2025bounded} for general $s$. Since $\frac{s\alpha}{n+s}$ increases from $\frac{2\alpha}{n+2}$ to $\alpha$ as $s$ increases from $2$ to $\infty$, the condition \eqref{cond : main assumption with s} interpolates between \eqref{cond : main assumption with 2} and \eqref{cond : main assumption with infty}. In the Orlicz setting the role of the gap $q-p$ is played by the quotient $H/G$ evaluated at the scale $\rho^{-1}$ or $\rho^{-\frac{n+s}{s}}$, in the same way as in the elliptic results of \cite{Byun2020,Baasandorj2020}. We stress that the additional integrability $u\in L^\infty(0,T;L^s(\Omega,\mr^N))$ is uniform in time, which distinguishes the parabolic interpolation from the elliptic one in \cite{Ok2020,DeFilippisMingione2021AMP}. Finally, the convexity of $\tilde{G}(\varkappa)=G(\varkappa^{\frac{1}{2}})$ together with $i_H\geq 2$ plays the role of the restriction $p\geq 2$ in the power case, see Lemma \ref{lem : quadratic ainc} and Remark \ref{rem : special cases}.

\begin{remark}\label{rem : special cases}
The following pairs of Young functions, in which $p\geq 2$, satisfy our assumptions, see also \cite{Byun2020}.
\begin{enumerate}
    \item $G(\varkappa)=\varkappa^p$ and $H(\varkappa)=\varkappa^q$, where $q>p$.
    \item $G(\varkappa)=\varkappa^p$ and $H(\varkappa)=\varkappa^p \log(1+\varkappa)$.
    \item $G(\varkappa)=\varkappa^p$ and $H(\varkappa)=\varkappa^p [\log(1+\varkappa)]^{\gamma}$, where $\gamma>0$.
    \item $G(\varkappa)=\varkappa^p$ and $H(\varkappa)=\varkappa^p \log \log(e+\varkappa)$.
\end{enumerate}
Here $\tilde{G}(\varkappa)=\varkappa^{\frac{p}{2}}$ is convex and one may take $i_H=p$, so that our requirements on $\tilde{G}$ and $i_H$ reduce to $p\geq 2$. In the cases (2)--(4) the normalization $H(1)=1$ is achieved by multiplying $H$ by a constant, which only rescales $a$. These cases are borderline double phase models in the spirit of \cite{BaroniColomboMingione2016SPMJ}, and for them each of the conditions on $\omega$ in \eqref{cond : main assumption with infty}, \eqref{cond : main assumption with 2} and \eqref{cond : main assumption with s} is equivalent to the finiteness of
$$
\limsup_{\rho\to 0+}\omega(\rho)\log\frac{1}{\rho},\qquad \limsup_{\rho\to 0+}\omega(\rho)\left(\log\frac{1}{\rho}\right)^{\gamma}\qquad\text{and}\qquad \limsup_{\rho\to 0+}\omega(\rho)\log\log\frac{1}{\rho},
$$
respectively, as in the elliptic case \cite[Remark 4.13]{Byun2020}.
\end{remark}

A key ingredient in our proof is an estimate of the parabolic
gluing term that retains the actual growth of the two phases.
For this purpose, we use
\[
 \Xi_\Phi(y)\coloneq \frac{y}{\Phi^{-1}(y)},\qquad
 \Xi_\Phi(\Phi(\varkappa))=\frac{\Phi(\varkappa)}{\varkappa}.
\]
Lemma~\ref{lem : properties of Xi} establishes the monotonicity of $\Xi_\Phi$ and the
scaling properties of
$\widehat{\Xi}_\Phi(y)\coloneq \Xi_\Phi\bigl(y^{\frac1\theta}\bigr)$.
Together with the averaging inequality in Lemma~\ref{lem : quasiconcave Jensen}, these properties yield the flux estimate in Lemma~\ref{lem : last term estimate of lemma 2.14}. In the $G$-phase, this estimate uses $\Phi=G$ and $\Phi=AH$, where $A\coloneq \sup_{Q_{10\rho}(z_0)}a$, and the second term is taken to be zero when $A=0$. It is then used in Lemma~\ref{lem : G phase Poincare} to control the gluing term on $G$-intrinsic cylinders and obtain the corresponding Poincar\'e type estimate. Lemma~\ref{lem : H phase Poincare} derives the analogous estimate for the weighted $H$-term from the same argument.
In the $(G,H)$-phase, the comparability established in
Lemma~\ref{lem : Psi z_0 comparability} permits us to work with the frozen density
$\Psi_{z_0}=G+a(z_0)H$. Lemma~\ref{lem : GH phase Poincare} applies the inverse-function argument with $\Phi=\Psi_{z_0}$ to obtain the corresponding
Poincar\'e type estimate. In both phases, the argument preserves
the cancellation
\[
 \frac{\lambda^2}{\Phi(\lambda)}
 \Xi_\Phi(\Phi(\lambda))=\lambda,
 \qquad \Phi\in\{G,\Psi_{z_0}\},
\]
which matches the flux estimate to the intrinsic time scale.
This provides the gluing estimates needed for the reverse H\"{o}lder argument while keeping the full growth functions.

The paper is organized as follows. In Section \ref{section 2} we collect the properties of Young functions, the gap inequalities and the energy estimates for weak solutions. Section \ref{section 3} contains the stopping time argument and the proof that the third case does not occur. In Section \ref{section 4} we prove the reverse H\"{o}lder inequalities in the $G$-phase and in the $(G,H)$-phase, and in Section \ref{section 5} we combine them with a Vitali type covering argument to prove Theorems \ref{thm : main theorem with infty} and \ref{thm : main theorem with s}.

\section{Preliminaries}\label{section 2}
We first set up notation. Recall that $\Psi(z,\varkappa)=G(\varkappa)+a(z)H(\varkappa)$ for $z\in\Omega_T$ and $\varkappa\geq 0$, see \eqref{def : definition of Psi}. For a fixed point $z_0 \in \Omega_T$, we denote
\begin{equation}\label{def : definition of H with a fixed center z_0}
    \Psi_{z_0}(\varkappa) \coloneq \Psi(z_0,\varkappa)=G(\varkappa)+a(z_0)H(\varkappa)
    \qquad \text{for } \varkappa\geq 0.
\end{equation}
We set intrinsic cylinders
\begin{equation}\label{eq : definition of p-intrinsic cylinder}
Q_\rho^\lambda(z_0)\coloneq B_\rho(x_0)\times \left(t_0-\frac{\lambda^2}{G(\lambda)}\rho^2,t_0+\frac{\lambda^2}{G(\lambda)}\rho^2 \right)\eqcolon B_\rho(x_0)\times I_{\rho}^\lambda(t_0)
\end{equation}
and
\begin{equation}\label{eq : definition of p,q-intrinsic cylinder}
\widetilde{Q}_\rho^\lambda(z_0)\coloneq B_\rho(x_0)\times\left(t_0-\frac{\lambda^2}{\Psi_{z_0}(\lambda)}\rho^2,t_0+\frac{\lambda^2}{\Psi_{z_0}(\lambda)}\rho^2\right)\eqcolon B_\rho(x_0)\times J_\rho^\lambda(t_0).
\end{equation}
For convenience, we write 
\begin{equation*}
    \operatorname{data}=\begin{cases}
        \operatorname{data}_b &\text{if \eqref{cond : main assumption with infty} holds,}\\
        \operatorname{data}_s &\text{if \eqref{cond : main assumption with s} holds}.
    \end{cases}
\end{equation*}
Next, we denote the super-level set by
\begin{equation}\label{def : definition of Upsilon}
    \Upsilon(\Lambda)\coloneq \{z\in\Omega_T: \Psi(z,|Du(z)|)>\Lambda\}.
\end{equation}

For any measurable set $E\subset \Omega_T$ with $0<|E|<\infty$ and integrable function $f\in L^1(\Omega_T)$, we denote the integral average of $f$ over $E$ by
$$
f_E=\frac{1}{|E|}\iints{E} f\, dz=\miint{E} f \, dz,
$$
where $|E|$ denotes the $(n+1)$-dimensional Lebesgue measure of $E$.

\subsection{Young functions and Orlicz spaces}
A function $\Phi: (0,\infty)\rightarrow \mr$ is almost increasing if there exists a constant $\tilde{L}\geq 1$ such that $\Phi(\varkappa_1)\leq \tilde{L}\Phi(\varkappa_2)$ for all $0<\varkappa_1<\varkappa_2$. A Young function $\Phi : [0,\infty) \rightarrow [0,\infty)$ is an increasing convex function satisfying
\begin{equation*}
\Phi(0)=0, \ \lim_{\varkappa \to \infty} \Phi(\varkappa) = \infty, \ \lim_{\varkappa \to 0+} \frac{\Phi(\varkappa)}{\varkappa} = 0, \, \lim_{\varkappa \to \infty} \frac{\Phi(\varkappa)}{\varkappa} = \infty.
\end{equation*}

\begin{definition}\label{def : nabla and Delta}
Let $\Phi$ be a Young function.
\begin{enumerate}
\item $\Phi$ is said to satisfy the $\Delta_2$-condition, denoted by $\Phi \in \Delta_2$, if there exists a positive number $\Delta_2(\Phi)$ such that $\Phi(2\varkappa) \leq \Delta_2(\Phi) \, \Phi(\varkappa)$ for all $\varkappa \geq 0$.
\item $\Phi$ is said to satisfy the $\nabla_2$-condition, denoted by $\Phi \in \nabla_2$, if there exists a positive number $\nabla_2(\Phi)>1$ such that $\Phi(\nabla_2(\Phi) \, \varkappa) \geq 2\nabla_2(\Phi) \, \Phi(\varkappa)$ for all $\varkappa \geq 0$.
\item We write $\Phi \in \Delta_2 \cap \nabla_2$ if $\Phi \in \Delta_2$ and $\Phi \in \nabla_2$.
\end{enumerate}
\end{definition}

We note that if $\Phi \in \Delta_2$, then $\Delta_2(\Phi)>2$.
Indeed, by the convexity of $\Phi$, we get
\begin{equation}
\label{2-02}
\Phi(2\varkappa) \leq \Delta_2(\Phi) \, \Phi(\varkappa) \leq \frac{\Delta_2(\Phi)}{2} \, \Phi(2\varkappa) \ \ \textrm{for all} \ \varkappa \geq 0,
\end{equation}
and hence $\Delta_2(\Phi) \geq 2$.
If $\Delta_2(\Phi)=2$, then it follows from (\ref{2-02}) that $\Phi(2\varkappa)=2\Phi(\varkappa)$ for all $\varkappa \geq 0$, and so $\Phi(\varkappa) \equiv \Phi(1) \, \varkappa$ is not a Young function.
Thus $\Delta_2(\Phi)>2$. Also, we obtain 
\begin{equation}\label{eq : subadditivity of Phi}
\Phi(\varkappa_1+\varkappa_2)\leq \frac{\Delta_2(\Phi)}{2}(\Phi(\varkappa_1)+\Phi(\varkappa_2))\quad \text{for }\varkappa_1,\varkappa_2>0.
\end{equation}
Indeed, the definition of $\Delta_2$ implies that
$$
\Phi(\varkappa_1+\varkappa_2)\leq \Phi\left(2\cdot\frac{\varkappa_1+\varkappa_2}{2}\right)\leq \Delta_2(\Phi)\Phi\left(\frac{\varkappa_1+\varkappa_2}{2}\right)\leq \frac{\Delta_2(\Phi)}{2}(\Phi(\varkappa_1)+\Phi(\varkappa_2))
$$
for $\varkappa_1,\varkappa_2>0$.

For a given Young function $\Phi$, we define the complementary Young function $\Phi^*$ of $\Phi$ by
\begin{equation*}
\Phi^*(\vartheta) = \sup \{\varkappa \vartheta-\Phi(\varkappa) : \varkappa \geq 0 \}.
\end{equation*}
We remark that $\Phi^*$ satisfies all the conditions to be a Young function and that $(\Phi^*)^*=\Phi$.
Moreover, $\Phi \in \nabla_2$ if and only if $\Phi^* \in \Delta_2$ with $2\nabla_2(\Phi)=\Delta_2(\Phi^*)$. The following properties of Young functions in $\Delta_2\cap\nabla_2$ are taken from \cite[Lemma 2.2]{Byun2020}.

\begin{lemma}
\label{lem202}
Let $\Phi$ be a Young function with $\Phi \in \Delta_2 \cap \nabla_2$.
\begin{enumerate}
\item For any $\mu\geq 1$ and $\varkappa \geq 0$, we have
\begin{equation}\label{eq : Phi(Lambda t) vs Phi(t)}
\Phi(\mu \varkappa) \leq \Delta_2(\Phi) \mu^{\log_2 \Delta_2(\Phi)} \, \Phi(\varkappa),
\end{equation}
and for any $0<\mu\leq 1$ and $\varkappa\geq 0$, we have
\begin{equation}
\Phi(\mu \varkappa) \leq 2\nabla_2(\Phi) \mu^{1+\log_{\nabla_2(\Phi)}2} \, \Phi(\varkappa).
\end{equation}
\item (Young's inequality) For any $\varepsilon \in (0,1]$, there exists a positive constant $c$ depending only on $\Delta_2(\Phi)$, $\nabla_2(\Phi)$ and $\varepsilon$ such that
\begin{equation}
\label{Youngs-ineq}
\varkappa_1\varkappa_2 \leq \varepsilon \Phi(\varkappa_1) + c \Phi^*(\varkappa_2), \quad \forall \varkappa_1,\varkappa_2 \geq 0.
\end{equation}
\item If $\Phi \in C^1([0,\infty))$, then for any $\varkappa \geq 0$, we have
\begin{equation}\label{eq : comparison with Phi and Phi'}
c_1^{-1} \Phi(\varkappa) \leq \varkappa\Phi'(\varkappa) \leq c_1 \Phi(\varkappa)
\end{equation}
and
\begin{equation}\label{eq : Phi*(Phi')}
\Phi^* \left( \Phi'(\varkappa) \right) \leq c_2 \Phi(\varkappa)
\end{equation}
for some constants $c_1,c_2\geq 1$ depending only on $\Delta_2(\Phi)$ and $\nabla_2(\Phi)$.
\item (A modified form of Young's inequality) If $\Phi \in C^1([0,\infty))$, then for any $\varepsilon \in (0,1]$, there exists a positive constant $c$ depending only on $\Delta_2(\Phi)$, $\nabla_2(\Phi)$ and $\varepsilon$ such that
\begin{equation}
\label{modified-Young-ineq}
\varkappa_1 \Phi'(\varkappa_2) \leq \varepsilon \Phi(\varkappa_1) + c \Phi(\varkappa_2), \quad \forall \varkappa_1,\varkappa_2 \geq 0.
\end{equation}
\end{enumerate}
\end{lemma}

\begin{lemma}\label{lem : Phi and Phi'}
    Let $\Phi\in C^1([0,\infty))$ be a Young function satisfying the $\Delta_2$-condition. Then
    $$
    \Phi(\varkappa)\leq \varkappa\Phi'(\varkappa)\leq \Delta_2(\Phi)\Phi(\varkappa).
    $$
\end{lemma}
\begin{proof}
    This follows from \cite[Lemma 2.6.6]{2011Diening} together with the definition of the $\Delta_2$-condition.
\end{proof}

\begin{definition}\label{def : relation of Phi1 and Phi2}
Let $\Phi_1$, $\Phi_2$ be Young functions.
We shall write
\begin{equation*}
\Phi_1 \prec \Phi_2
\end{equation*}
if $\Phi_2 \circ \Phi_1^{-1}$ is a Young function.
\end{definition}
The next two lemmas are \cite[Lemmas 2.5 and 2.7]{Byun2020}.

\begin{lemma}
\label{lem205}
Let $\Phi_1$, $\Phi_2$ be Young functions with $\Phi_1 \prec \Phi_2$.
Then
\begin{equation}
\label{ele}
\Phi_1(\varkappa) \leq \frac{1}{\left( \Phi_2 \circ \Phi_1^{-1} \right)(1)} \, \Phi_2(\varkappa), \quad \forall \varkappa \geq \Phi_1^{-1}(1).
\end{equation}
\end{lemma}

\begin{lemma}
\label{lem208}
Let $\Phi_1$, $\Phi_2$ be Young functions with $\Phi_1 \prec \Phi_2$.
Then the function
\begin{equation*}
\varkappa \mapsto \left( \frac{\Phi_2}{\Phi_1} \right)(\varkappa) = \frac{\Phi_2(\varkappa)}{\Phi_1(\varkappa)}
\end{equation*}
is non-decreasing.
\end{lemma}

\begin{lemma}\label{lem : Gap inequality with s}
    If \eqref{cond : main assumption with s} holds, then for any $0<\rho\leq 1$
    $$
    \omega(\rho)H(\varkappa)\leq \varpi_s\left[1+\left(\frac{\varkappa}{M}\right)^{\frac{s}{n+s}}\right]G(\varkappa)\qquad \text{for all }\varkappa\geq 0,
    $$
    where $M\coloneq \rho^{-\frac{n+s}{s}}$.
\end{lemma}
\begin{proof}
    If $\omega\equiv0$, we are done. Now, we assume that $\omega(\tau_0)>0$ for some $\tau_0> 0$. Since $\omega$ is a nonnegative concave function with $\omega(0)=0$, we have $\omega(\tau)>0$ for all $\tau>0$.

    Note that the function $\tau\mapsto \frac{\omega(\tau)}{\tau}$ is non-increasing on $(0,\infty)$. Indeed, for $\tau_1<\tau_2$ in $(0,\infty)$, \eqref{cond : property of concave function} implies
    $$
    \omega(\tau_2)=\omega\left(\frac{\tau_2}{\tau_1}\tau_1\right)\leq \frac{\tau_2}{\tau_1}\omega(\tau_1).
    $$

    If $\varkappa\leq M$, then, since $G\prec H$, \eqref{cond : property of s} and Lemma \ref{lem208} give 
    $$
    \omega(\rho)\frac{H(\varkappa)}{G(\varkappa)}\leq \omega(\rho)\frac{H(M)}{G(M)}\leq\varpi_s.
    $$

    If $\varkappa>M$, we set $\rho_\varkappa\coloneq \varkappa^{-\frac{s}{n+s}}$, so that $\varkappa=\rho_\varkappa^{-\frac{n+s}{s}}$ and $0<\rho_\varkappa<\rho\leq 1$. By \eqref{cond : property of s} with $\rho=\rho_\varkappa$ and the monotonicity of $\tau\mapsto \frac{\omega(\tau)}{\tau}$, we have
    $$
    \begin{aligned}
    \omega(\rho)H(\varkappa)=\omega(\rho_\varkappa)G(\varkappa)\frac{\omega(\rho)}{\omega(\rho_\varkappa)}\cdot\frac{H(\varkappa)}{G(\varkappa)}&\leq \omega(\rho_\varkappa)G(\varkappa)\frac{\rho}{\rho_\varkappa}\cdot\frac{\varpi_s}{\omega(\rho_\varkappa)}\\
    &=\varpi_s\left(\frac{\varkappa}{M}\right)^{\frac{s}{n+s}}G(\varkappa).
    \end{aligned}
    $$
    Hence, combining both cases leads to the desired conclusion.
\end{proof}
\begin{lemma}\label{lem : Gap inequality with infty}
    If \eqref{cond : main assumption with infty} holds, then for any $0<\rho\leq 1$
    $$
    \omega(\rho)H(\varkappa)\leq \varpi_b\left[1+\rho\varkappa\right]G(\varkappa)\quad \text{for all }\varkappa\geq 0.
    $$
\end{lemma}
\begin{proof}
    As in the above proof, we assume that $\omega(\tau)>0$ for all $\tau>0$.

    If $\varkappa\leq \rho^{-1}$, then, since $G\prec H$, \eqref{cond : property of infty condition} and Lemma \ref{lem208} give
    $$
    \omega(\rho)\frac{H(\varkappa)}{G(\varkappa)}\leq \omega(\rho)\frac{H(\rho^{-1})}{G(\rho^{-1})}\leq \varpi_b.
    $$
    Hence $\omega(\rho)H(\varkappa)\leq \varpi_b G(\varkappa)$.

    If $\varkappa > \rho^{-1}$, we set $\rho_\varkappa\coloneq \varkappa^{-1}$, so that $\varkappa=\rho_{\varkappa}^{-1}$ and $0<\rho_\varkappa<\rho\leq 1$. Since $\tau\mapsto \frac{\omega(\tau)}{\tau}$ is non-increasing, by \eqref{cond : property of infty condition} with $\rho=\rho_\varkappa$, we get
    $$
    \omega(\rho)H(\varkappa)=G(\varkappa)\frac{\omega(\rho)}{\omega(\rho_\varkappa)}\omega(\rho_\varkappa)\frac{H(\varkappa)}{G(\varkappa)}\leq \varpi_b\rho\varkappa G(\varkappa).
    $$
    Thus, we have the desired conclusion.
\end{proof}
\begin{lemma}\label{lem : Young function of Psi_z_0}
    The function $\varkappa\mapsto\Psi_{z_0}(\varkappa^\frac{1}{2})$, with $\Psi_{z_0}$ as in \eqref{def : definition of H with a fixed center z_0}, is convex on $[0,\infty)$.
\end{lemma}
\begin{proof}
    Since $\tilde{G}(\varkappa)=G(\varkappa^\frac{1}{2})$ is convex by assumption, it suffices to show that $\varkappa\mapsto H(\varkappa^\frac{1}{2})$ is convex. Since $G\prec H$, the function $H\circ G^{-1}$ is a Young function, and in particular it is convex and non-decreasing. For $\varkappa\geq 0$ we have
    $$
    H(\varkappa^\frac{1}{2})=(H\circ G^{-1})(\tilde{G}(\varkappa)),
    $$
    and the composition of a non-decreasing convex function with a convex function is convex.
\end{proof}
The following lemma is the Orlicz counterpart of the restriction $p\geq 2$ in the power case.
\begin{lemma}\label{lem : quadratic ainc}
    Both $\displaystyle (0,\infty)\ni \varkappa\mapsto \frac{G(\varkappa)}{\varkappa^2}$ and
    $\displaystyle (0,\infty)\ni \varkappa\mapsto \frac{H(\varkappa)}{\varkappa^2}$ are
    non-decreasing. Consequently, for $\Phi\in\{G,H\}$, $0<c\leq 1$ and $\varkappa,\,y>0$,
    \begin{equation}\label{eq : quadratic scaling}
        \Phi(c\varkappa)\leq c^2\Phi(\varkappa)
        \qquad\text{and}\qquad
        \Phi^{-1}(cy)\geq c^{\frac12}\Phi^{-1}(y).
    \end{equation}
\end{lemma}
\begin{proof}
    By assumption $\tilde{G}(\varkappa)=G(\varkappa^{\frac12})$ is convex with $\tilde{G}(0)=0$, and therefore
    $\varkappa\mapsto\tilde{G}(\varkappa)/\varkappa$ is non-decreasing. Replacing $\varkappa$
    by $\varkappa^2$ gives the assertion for $G$. Since $G\prec H$, Lemma \ref{lem208} shows
    that $H/G$ is non-decreasing, and hence
    $$
    \frac{H(\varkappa)}{\varkappa^2}=\frac{H(\varkappa)}{G(\varkappa)}\cdot\frac{G(\varkappa)}{\varkappa^2}
    $$
    is a product of two nonnegative non-decreasing functions. The first inequality in
    \eqref{eq : quadratic scaling} is immediate from the monotonicity of
    $\Phi(\varkappa)/\varkappa^2$, and the second one follows by applying the first one with
    $\varkappa=c^{\frac12}\Phi^{-1}(y)$, which gives
    $\Phi\bigl(c^{\frac12}\Phi^{-1}(y)\bigr)\leq cy$.
\end{proof}
The next lemma replaces Young's inequality \eqref{Youngs-ineq} for the function $\tilde{G}$, which need not belong to $\nabla_2$.
\begin{lemma}\label{lem : tilde G conjugate}
    Let $\tilde{G}^*(\vartheta)\coloneq \sup_{\varkappa\geq 0}\,[\varkappa\vartheta-\tilde{G}(\varkappa)]$ for $\vartheta\geq 0$ be the conjugate function of $\tilde{G}$. Then
    \begin{equation}\label{eq : tilde G property}
        \tilde{G}^*\left(\frac{G(\lambda)}{\lambda^2}\right)\leq G(\lambda)\qquad\text{for every }\lambda>0,
    \end{equation}
    and for every $\varepsilon\in(0,1]$ there exists a constant $c=c(\varepsilon,\Delta_2(G))>0$ such that
    \begin{equation}\label{eq : Young inequality for tilde G}
        \varkappa_1\varkappa_2\leq \varepsilon\tilde{G}^*(\varkappa_1)+c\,\tilde{G}(\varkappa_2)\qquad\text{for all }\varkappa_1,\,\varkappa_2\geq 0.
    \end{equation}
\end{lemma}
\begin{proof}
    As in the proof of Lemma \ref{lem : quadratic ainc}, the function $\varkappa\mapsto\tilde{G}(\varkappa)/\varkappa$ is non-decreasing on $(0,\infty)$, and $\tilde{G}(\lambda^2)/\lambda^2=G(\lambda)/\lambda^2$. Hence $\varkappa G(\lambda)/\lambda^2-\tilde{G}(\varkappa)\leq 0$ for $\varkappa\geq\lambda^2$, and therefore
    $$
    \tilde{G}^*\left(\frac{G(\lambda)}{\lambda^2}\right)=\sup_{0\leq\varkappa\leq\lambda^2}\left[\varkappa\frac{G(\lambda)}{\lambda^2}-\tilde{G}(\varkappa)\right]\leq \sup_{0\leq\varkappa\leq\lambda^2}\varkappa\frac{G(\lambda)}{\lambda^2}=G(\lambda),
    $$
    which is \eqref{eq : tilde G property}. Next, being a supremum of affine functions, $\tilde{G}^*$ is convex, and $\tilde{G}^*(0)=0$ since $\tilde{G}\geq 0=\tilde{G}(0)$. Thus $\tilde{G}^*(\varepsilon\varkappa_1)\leq\varepsilon\tilde{G}^*(\varkappa_1)$ for $\varepsilon\in(0,1]$, and the definition of $\tilde{G}^*$ gives
    $$
    \varkappa_1\varkappa_2=(\varepsilon\varkappa_1)\frac{\varkappa_2}{\varepsilon}\leq \tilde{G}^*(\varepsilon\varkappa_1)+\tilde{G}\left(\frac{\varkappa_2}{\varepsilon}\right)\leq \varepsilon\tilde{G}^*(\varkappa_1)+\tilde{G}\left(\frac{\varkappa_2}{\varepsilon}\right).
    $$
    Finally, since $G$ is increasing and satisfies the $\Delta_2$-condition, we have $\tilde{G}(2\varkappa)=G(2^{\frac12}\varkappa^{\frac12})\leq G(2\varkappa^{\frac12})\leq\Delta_2(G)\tilde{G}(\varkappa)$, and iterating this inequality yields $\tilde{G}(\varkappa_2/\varepsilon)\leq c(\varepsilon,\Delta_2(G))\tilde{G}(\varkappa_2)$. This proves \eqref{eq : Young inequality for tilde G}.
\end{proof}

\subsection{Energy lemmas}

To prove the main theorems, we refer to three energy lemmas from \cite{2023_Gradient_Higher_Integrability_for_Degenerate_Parabolic_Double-Phase_Systems}. The following lemma provides a Caccioppoli type inequality.
\begin{lemma} \label{lem : Caccioppoli inequality}
    Let $u$ be a weak solution to \eqref{eq : main equation}. Then there exists a positive constant $c=c(n,\nu,L,\Delta_2(G),\nabla_2(G),\Delta_2(H),\nabla_2(H))>1$ such that
    $$
    \begin{aligned}
        &\sup_{t\in (t_0-\tau,t_0+\tau)}\dashint_{B_r(x_0)} \frac{|u-u_{Q_{r,\tau}(z_0)}|^2}{\tau}\, dx+\miint{Q_{r,\tau}(z_0)} \Psi(z,|Du|)\, dz\\
        &\quad \leq c \miint{Q_{R, \ell}\left(z_0\right)}\left[G\left(\frac{\left|u-u_{Q_{R, \ell}\left(z_0\right)}\right|}{R-r}\right)+a(z) H\left(\frac{\left|u-u_{Q_{R, \ell}\left(z_0\right)}\right|}{R-r}\right)\right] dz \\
        &\qquad+c \miint{Q_{R, \ell}\left(z_0\right)} \frac{\left|u-u_{Q_{R, \ell}\left(z_0\right)}\right|^2}{\ell-\tau} \, dz
    \end{aligned}
    $$
    for every $Q_{R, \ell}\left(z_0\right)=B_R(x_0) \times(t_0-\ell, t_0+\ell) \subset \Omega_T$, with $R, \ell>0,\, r \in[R / 2, R)$ and $\tau \in [\ell / 4, \ell )$.
\end{lemma}
\begin{proof}
    Let $\eta\in C_0^\infty (B_R(x_0))$ be a cut-off function with
    $$
    0\leq \eta \leq 1,\quad \eta \equiv 1 \text{ in }B_r(x_0)\quad \text{and}\quad \|D\eta\|_{L^\infty}\leq \frac{2}{R-r}.
    $$
    For $\tau\in[\frac{\ell}{4},\ell)$, let $h_0>0$ be sufficiently small so that there exists a cut-off function $\zeta\in C_0^\infty(I_{\ell-h_0}(t_0))$ with
    $$
    0\leq \zeta\leq 1, \quad \zeta\equiv 1\text{ in }I_\tau(t_0)\quad\text{and}\quad \|\partial_t\zeta\|_{L^\infty}\leq \frac{3}{\ell-\tau}.
    $$
    Let $t_* \in I_\tau(t_0)$ and $\delta\in (0,h_0)$. We define $\zeta_\delta$ as 
    $$
    \zeta_\delta (t)=\begin{cases}
        1,\qquad &t\in(-\infty,t_*-\delta),\\
        1-\frac{t-t_*+\delta}{\delta}, &t\in [t_*-\delta,t_*],\\
        0,&t\in (t_*,\infty).
    \end{cases}
    $$
    For $h\in (0,h_0)$, the Steklov averages satisfy
    \begin{equation}    \label{eq : main equation w.r.t Steklov average}
    \partial_t [u-u_{Q_{R,\ell}(z_0)}]_h-\operatorname{div} [\mathcal{A}(\cdot,Du)]_h = 0
    \end{equation}
    in $B_R(x_0)\times I_{\ell-h}(t_0)$. Observe that
    $$
    \begin{aligned}
    [u-u_{Q_{R,\ell}(z_0)}]_h \eta^{\sigma_H}\zeta^2\zeta_\delta \in &\,W_0^{1,2}(I_{\ell-h};L^2(B_R(x_0),\mr^N))\\
    &\qquad \cap\, L^{\sigma_H}(I_{\ell-h};W_0^{1,\sigma_H}(B_R(x_0),\mr^N)).
    \end{aligned}
    $$
    Using $\varphi=[u-u_{Q_{R,\ell}(z_0)}]_h \eta^{\sigma_H}\zeta^2\zeta_\delta$ as a test function in \eqref{eq : main equation w.r.t Steklov average}, we obtain
    \begin{equation}\label{eq : I+II=0 in Caccioppoli inequality}
    \mathrm{I}+\mathrm{II}=\miint{Q_{R,\ell}(z_0)} \partial_t [u-u_{Q_{R,\ell}(z_0)}]_h\cdot \varphi \, dz + \miint{Q_{R,\ell}(z_0)} [\mathcal{A}(\cdot,Du)]_h\cdot D\varphi\, dz =0.
    \end{equation}
    We claim that $\mathrm{II}$ is finite under the assumption $|Du|\in L^{\sigma_H}(\Omega_T)$. It follows from \eqref{cond : double phase bounded condition of integrand} and \eqref{eq : comparison with Phi and Phi'} that
    \begin{equation}    \label{cond : growth condition of A with G' and H'}
    |[\mathcal{A}(z,Du)]_h|\leq \tilde{c}[G'(|Du|)+a(z)H'(|Du|)]_h
    \end{equation}
    for some $\tilde{c}=\tilde{c}(L,\Delta_2(G),\nabla_2(G),\Delta_2(H),\nabla_2(H))>1$.
    Then, we have
    $$
        \mathrm{II}\leq c\miint{Q_{R,\ell}(z_0)} (G'(|Du|))_h|D\varphi|\, dz+c\miint{Q_{R,\ell}(z_0)} (aH'(|Du|))_h|D\varphi|\, dz.
    $$
    Then we get from \eqref{modified-Young-ineq} that
    $$
    \begin{aligned}
        \mathrm{II} &\leq c\miint{Q_{R,\ell}(z_0)} \left[(G(|Du|))_h+G(|D\varphi|)\right] dz\\
        &\quad + c\miint{Q_{R,\ell}(z_0)} \left[(aH(|Du|))_h+a_hH(|D\varphi|)\right] dz.
    \end{aligned}
    $$
    Since $G(\varkappa)+H(\varkappa)\leq 2(1+\varkappa^{\sigma_H})$ for $\varkappa\geq 0$ by \eqref{cond : two sided power bound of Phi} and \eqref{cond : relation of sigma G and sigma H}, and since $a$ is bounded and $|Du|, |D\varphi|\in L^{\sigma_H}(\Omega_T)$, the properties of the Steklov average show that $\mathrm{II}$ is finite.

    \textbf{Estimate of $\mathrm{I}$:} By integration by parts, we obtain
    $$
    \begin{aligned}
    \mathrm{I}&=\miint{Q_{R,\ell}(z_0)}\frac{1}{2} (\partial_t|[u-u_{Q_{R,\ell}(z_0)}]_h|^2) \eta^{\sigma_H}\zeta^2\zeta_\delta\, dz\\
    &= - \miint{Q_{R,\ell}(z_0)} |[u-u_{Q_{R,\ell}(z_0)}]_h|^2 \eta^{\sigma_H}\zeta \zeta_\delta \partial_t\zeta\, dz\\
    &\qquad - \miint{Q_{R,\ell}(z_0)}  \frac{1}{2} |[u-u_{Q_{R,\ell}(z_0)}]_h|^2 \eta^{\sigma_H} \zeta^2 \partial_t \zeta_\delta\, dz.
    \end{aligned}
    $$
    By the definition of $\zeta$, we get
    $$
    - \miint{Q_{R,\ell}(z_0)} |[u-u_{Q_{R,\ell}(z_0)}]_h|^2 \eta^{\sigma_H}\zeta \zeta_\delta \partial_t\zeta\, dz \geq -3\miint{Q_{R,\ell}(z_0)} \frac{|[u-u_{Q_{R,\ell}(z_0)}]_h|^2}{\ell-\tau}\, dz.
    $$
    In addition, the definition of $\zeta_\delta$ implies
    $$
    \begin{aligned}
    &- \miint{Q_{R,\ell}(z_0)}  \frac{1}{2} |[u-u_{Q_{R,\ell}(z_0)}]_h|^2 \eta^{\sigma_H} \zeta^2 \partial_t \zeta_\delta\, dz\\
    &\qquad = \frac{1}{|Q_{R,\ell}|}\dashint_{t_*-\delta}^{t_*}\int_{B_R(x_0)} \frac{1}{2}|[u-u_{Q_{R,\ell}(z_0)}]_h|^2 \eta^{\sigma_H} \zeta^2 \, dxdt\\
    &\qquad \geq \frac{1}{|Q_{R,\ell}|} \dashint_{t_*-\delta}^{t_*}\int_{B_r(x_0)} \frac{1}{2}|[u-u_{Q_{R,\ell}(z_0)}]_h|^2 \, dxdt.
    \end{aligned}
    $$
    Thus, we get
    \begin{equation}\label{eq : estimate I in Caccioppoli inequality}
    \begin{aligned}
    \lim_{h\rightarrow 0^+}\lim_{\delta\rightarrow 0^+} \mathrm{I}&\geq -3\miint{Q_{R,\ell}(z_0)} \frac{|u-u_{Q_{R,\ell}(z_0)}|^2}{\ell-\tau}\, dz\\
    &\qquad + \frac{1}{2|Q_{R,\ell}|}\int_{B_r(x_0)} |u(x,t_*)-u_{Q_{R,\ell}(z_0)}|^2\, dx.
    \end{aligned}
    \end{equation}

    \textbf{Estimate of $\mathrm{II}$:} We observe that
    $$
    \begin{aligned}
    \mathrm{II}&=\miint{Q_{R,\ell}(z_0)} [\mathcal{A}(\cdot,Du)]_h\cdot [Du]_h \eta^{\sigma_H}\zeta^2\zeta_\delta\, dz\\
    &\qquad +\sigma_H\miint{Q_{R,\ell}(z_0)} [\mathcal{A}(\cdot,Du)]_h \cdot [u-u_{Q_{R,\ell}(z_0)}]_h \eta^{\sigma_H-1} D\eta \zeta^2\zeta_\delta\, dz.
    \end{aligned}
    $$
    Using \eqref{cond : double phase bounded condition of integrand} yields
    $$
    \begin{aligned}
    &\lim_{h\rightarrow 0^+}\lim_{\delta\rightarrow 0^+} \miint{Q_{R,\ell}(z_0)} [\mathcal{A}(\cdot,Du)]_h\cdot [Du]_h \eta^{\sigma_H}\zeta^2\zeta_\delta\, dz\\
    &\geq \frac{\nu}{|Q_{R,\ell}|}\int_{I_{\ell}(t_0)\cap (-\infty,t_*)}\int_{B_R(x_0)} (G(|Du|)+a(z)H(|Du|))\eta^{\sigma_H}\zeta^2\,dxdt.
    \end{aligned}
    $$
    Furthermore, it follows from \eqref{cond : growth condition of A with G' and H'} and the definition of $\eta$ that
    $$
    \begin{aligned}
        &\lim_{h\rightarrow0^+}\lim_{\delta\rightarrow0^+} \miint{Q_{R,\ell}(z_0)} [\mathcal{A}(\cdot,Du)]_h\cdot [u-u_{Q_{R,\ell}(z_0)}]_h \eta^{\sigma_H-1}D\eta \zeta^2\zeta_\delta\,dz\\
        &\quad \geq -\frac{2\tilde{c}}{|Q_{R,\ell}|}\int_{I_\ell(t_0)\cap(-\infty,t_*)}\int_{B_R(x_0)} G'(|Du|)\eta^{\sigma_H-1}\zeta^2\frac{|u-u_{Q_{R,\ell}(z_0)}|}{R-r}\,dxdt\\
        &\qquad -\frac{2\tilde{c}}{|Q_{R,\ell}|}\int_{I_\ell(t_0)\cap(-\infty,t_*)}\int_{B_R(x_0)} a(z)H'(|Du|)\eta^{\sigma_H-1}\zeta^2\frac{|u-u_{Q_{R,\ell}(z_0)}|}{R-r}\,dxdt
    \end{aligned}
    $$
    for some $\tilde{c}=\tilde{c}(L,\Delta_2(G),\nabla_2(G),\Delta_2(H),\nabla_2(H))>1$. By \eqref{cond : relation of sigma G and sigma H} and \eqref{eq : Phi*(Phi')}, we note that 
    $$
    \frac{\sigma_G(\sigma_H-1)}{\sigma_G-1}\geq \sigma_H
    $$
    and 
    $$
    G^*(G'(\varkappa))\leq c_G G(\varkappa)\quad\text{and}\quad H^*(H'(\varkappa))\leq c_H H(\varkappa)
    $$
    for some $c_G=c_G(\Delta_2(G),\nabla_2(G))>1$ and $c_H=c_H(\Delta_2(H),\nabla_2(H))>1$. It follows from these facts, \eqref{cond : H* ct and H t with c<1} and \eqref{Youngs-ineq} that
    $$
    \begin{aligned}
    &G'(|Du|)\eta^{\sigma_H-1}\frac{|u-u_{Q_{R,\ell}(z_0)}|}{R-r}\\
    &\qquad\leq \frac{\nu}{8\sigma_H\tilde{c}c_G} G^*(G'(|Du|)\eta^{\sigma_H-1})+c G\left(\frac{|u-u_{Q_{R,\ell}(z_0)}|}{R-r}\right)\\
    &\qquad\leq \frac{\nu}{8\sigma_H\tilde{c}c_G}\eta^{\frac{\sigma_G(\sigma_H-1)}{\sigma_G-1}} G^*(G'(|Du|))+c G\left(\frac{|u-u_{Q_{R,\ell}(z_0)}|}{R-r}\right)\\
    &\qquad\leq\frac{\nu}{8\sigma_H\tilde{c}}\eta^{\sigma_H} G(|Du|)+c G\left(\frac{|u-u_{Q_{R,\ell}(z_0)}|}{R-r}\right)
    \end{aligned}
    $$
    and
    $$
    \begin{aligned}
    &H'(|Du|)\eta^{\sigma_H-1}\frac{|u-u_{Q_{R,\ell}(z_0)}|}{R-r}\\
    &\qquad \leq \frac{\nu}{8\sigma_H\tilde{c}c_H} H^*(H'(|Du|)\eta^{\sigma_H-1})+cH\left(\frac{|u-u_{Q_{R,\ell}(z_0)}|}{R-r}\right)\\
    &\qquad \leq \frac{\nu}{8\sigma_H\tilde{c}c_H} \eta^{\sigma_H} H^*(H'(|Du|))+cH\left(\frac{|u-u_{Q_{R,\ell}(z_0)}|}{R-r}\right)\\
    &\qquad \leq \frac{\nu}{8\sigma_H\tilde{c}} \eta^{\sigma_H} H(|Du|)+cH\left(\frac{|u-u_{Q_{R,\ell}(z_0)}|}{R-r}\right)
    \end{aligned}
    $$
    for some $c=c(\nu,L,\Delta_2(G),\nabla_2(G),\Delta_2(H),\nabla_2(H))>1$. Hence, we get
    $$
    \begin{aligned}
        &\lim_{h\rightarrow0^+}\lim_{\delta\rightarrow0^+} \miint{Q_{R,\ell}(z_0)} [\mathcal{A}(\cdot,Du)]_h \cdot [u-u_{Q_{R,\ell}(z_0)}]_h \eta^{\sigma_H-1} D\eta \zeta^2\zeta_\delta\, dz\\
        &\quad \geq -\frac{\nu}{4\sigma_H|Q_{R,\ell}|}\int_{I_\ell(t_0)\cap(-\infty,t_*)}\int_{B_R(x_0)} [G(|Du|)+a(z)H(|Du|)]\eta^{\sigma_H}\zeta^2\,dxdt\\
        &\qquad -c\miint{Q_{R,\ell}(z_0)} \left[G\left(\frac{|u-u_{Q_{R,\ell}(z_0)}|}{R-r}\right)+a(z)H\left(\frac{|u-u_{Q_{R,\ell}(z_0)}|}{R-r}\right)\right]\,dz.
        \end{aligned}
    $$
    Therefore, recalling the factor $\sigma_H$ in front of the second integral in the decomposition of $\mathrm{II}$ and using $\nu-\sigma_H\cdot\frac{\nu}{4\sigma_H}\geq \frac{\nu}{2}$, we see that
    \begin{equation}\label{eq : estimate II in Caccioppoli inequality}
    \begin{aligned}
    \lim_{h\rightarrow 0^+}\lim_{\delta\rightarrow 0^+} \mathrm{II}&\geq \frac{\nu}{2|Q_{R,\ell}|}\int_{I_\ell(t_0)\cap(-\infty,t_*)}\int_{B_R(x_0)} [G(|Du|)+a(z)H(|Du|)]\eta^{\sigma_H}\zeta^2\,dxdt\\
    &\quad-c\miint{Q_{R,\ell}(z_0)} \left[G\left(\frac{|u-u_{Q_{R,\ell}(z_0)}|}{R-r}\right)+a(z)H\left(\frac{|u-u_{Q_{R,\ell}(z_0)}|}{R-r}\right)\right]\,dz.
    \end{aligned}
    \end{equation}

    Combining \eqref{eq : I+II=0 in Caccioppoli inequality}, \eqref{eq : estimate I in Caccioppoli inequality} and \eqref{eq : estimate II in Caccioppoli inequality} implies that
    $$
    \begin{aligned}
        &\frac{1}{|Q_{R,\ell}|}\int_{B_r(x_0)} |u(x,t_*)-u_{Q_{R,\ell}(z_0)}|^2\, dx\\
        &\qquad\quad+ \frac{1}{|Q_{R,\ell}|}\int_{I_\ell(t_0)\cap(-\infty,t_*)}\int_{B_R(x_0)} [G(|Du|)+a(z)H(|Du|)]\eta^{\sigma_H}\zeta^2\,dxdt\\
        &\quad\leq c\miint{Q_{R,\ell}(z_0)} \left[G\left(\frac{|u-u_{Q_{R,\ell}(z_0)}|}{R-r}\right)+a(z)H\left(\frac{|u-u_{Q_{R,\ell}(z_0)}|}{R-r}\right)\right]\,dz\\
        &\qquad\quad+ c\miint{Q_{R,\ell}(z_0)} \frac{|u-u_{Q_{R,\ell}(z_0)}|^2}{\ell-\tau}\, dz.
    \end{aligned}
    $$
    Since $t_*\in I_\tau(t_0)$ is arbitrary, $|B_R|\simeq_n |B_r|$ and $|I_\ell|\simeq |I_\tau|$, we have
    $$
    \begin{aligned}
        &\sup_{t\in(t_0-\tau,t_0+\tau)}\dashint_{B_r(x_0)} \frac{|u-u_{Q_{R,\ell}(z_0)}|^2}{\tau}\, dx\\
        &\qquad\quad+ \miint{Q_{r,\tau}(z_0)} [G(|Du|)+a(z)H(|Du|)]\,dxdt\\
        &\quad\leq c\miint{Q_{R,\ell}(z_0)} \left[G\left(\frac{|u-u_{Q_{R,\ell}(z_0)}|}{R-r}\right)+a(z)H\left(\frac{|u-u_{Q_{R,\ell}(z_0)}|}{R-r}\right)\right]\,dz\\
        &\qquad\quad+ c\miint{Q_{R,\ell}(z_0)} \frac{|u-u_{Q_{R,\ell}(z_0)}|^2}{\ell-\tau}\, dz
    \end{aligned}
    $$
    for some $c=c(n,\nu,L,\Delta_2(G),\nabla_2(G),\Delta_2(H),\nabla_2(H))>1$. Finally, we replace the mean $u_{Q_{R,\ell}(z_0)}$ in the first term on the left-hand side by $u_{Q_{r,\tau}(z_0)}$. For every $t\in(t_0-\tau,t_0+\tau)$ we have
    $$
    \dashint_{B_r(x_0)}|u-u_{Q_{r,\tau}(z_0)}|^2\, dx\leq 2\dashint_{B_r(x_0)}|u-u_{Q_{R,\ell}(z_0)}|^2\, dx+2\left|u_{Q_{r,\tau}(z_0)}-u_{Q_{R,\ell}(z_0)}\right|^2
    $$
    and, by Jensen's inequality,
    $$
    \begin{aligned}
    \left|u_{Q_{r,\tau}(z_0)}-u_{Q_{R,\ell}(z_0)}\right|^2&\leq \miint{Q_{r,\tau}(z_0)}|u-u_{Q_{R,\ell}(z_0)}|^2\, dz\\
    &\leq \sup_{t\in(t_0-\tau,t_0+\tau)}\dashint_{B_r(x_0)}|u-u_{Q_{R,\ell}(z_0)}|^2\, dx.
    \end{aligned}
    $$
    Hence the first term on the left-hand side of the assertion is bounded by four times the corresponding term in the last display, and the proof is complete.
\end{proof}

The next lemma is a gluing lemma. For this, for a function $f\in L^1(Q_{R,\ell}(z_0),\mr^N)$ we denote its spatial integral average over $B_R(x_0)$ by
$$
f_{B_R(x_0)}(t)\coloneq\dashint_{B_R(x_0)} f(x,t)\, dx,
$$
and we write $(u\eta)_{B_R(x_0)}(t)$ when $f=u\eta$.
\begin{lemma} \label{lem : gluing lemma}
    Let $u$ be a weak solution to \eqref{eq : main equation}, and let $\eta \in C_0^{\infty}(B_R(x_0))$ be a function such that
    \begin{equation}    \label{cond : cut-off function eta}
    \eta \geq 0, \quad \dashint_{B_R(x_0)} \eta \, dx=1 \quad \text { and } \quad \|\eta\|_{L^{\infty}}+R\|D \eta\|_{L^{\infty}} \leq c(n).
    \end{equation}
    Then there exists a positive constant $c=c(n,L,\Delta_2(G),\nabla_2(G),\Delta_2(H),\nabla_2(H))$ such that
\begin{multline*}
    \sup_{t_1, t_2 \in (t_0-\ell, t_0+\ell)} |(u \eta)_{B_R(x_0)}(t_2)-(u \eta)_{B_R(x_0)}(t_1)| \\
    \leq c \frac{\ell}{R} \miint{Q_{R, \ell}(z_0)}\left[G'(|D u|)+a(z)H'(|D u|)
    \right] dz
\end{multline*}
    for every $Q_{R, \ell}(z_0) \subset \Omega_T$ with $R, \ell>0$.
\end{lemma}
\begin{proof}
    Let $t_1,\,t_2\in (t_0-\ell,t_0+\ell)$ with $t_1<t_2$. For a sufficiently small $\delta\in(0,1)$, we define $\zeta_\delta\in W_0^{1,\infty}(t_0-\ell,t_0+\ell)$ as follows:
    $$
    \zeta_\delta(t)=\begin{cases}
        0,&t_0-\ell\leq t\leq t_1-\delta,\\
        \frac{t-t_1+\delta}{\delta},&t_1-\delta<t<t_1,\\
        1,&t_1\leq t\leq t_2,\\
        \frac{t_2+\delta-t}{\delta},&t_2<t<t_2+\delta,\\
        0,&t_2+\delta\leq t \leq t_0+\ell.
    \end{cases}
    $$
    Testing \eqref{eq : main equation} with $\eta\zeta_\delta\in W_0^{1,\infty}(Q_{R,\ell}(z_0))$ and using \eqref{cond : double phase bounded condition of integrand} and \eqref{eq : comparison with Phi and Phi'}, we obtain
    $$
    \begin{aligned}
        &\left|\dashint_{t_1-\delta}^{t_1}\dashint_{B_R(x_0)} u\eta\, dxdt-\dashint_{t_2}^{t_2+\delta}\dashint_{B_R(x_0)} u\eta\, dxdt\right|\\
        \qquad & \leq c\int_{t_1-\delta}^{t_2+\delta}\dashint_{B_R(x_0)} [G'(|Du|)+a(z)H'(|Du|)]|D\eta||\zeta_\delta|\, dxdt
    \end{aligned}
    $$
    for some $c=c(L,\Delta_2(G),\nabla_2(G),\Delta_2(H),\nabla_2(H))>1$. Letting $\delta\rightarrow0^+$, we obtain
    $$
    \begin{aligned}
        |(u\eta)_{B_R(x_0)}(t_1)-(u\eta)_{B_R(x_0)}(t_2)|\leq c\frac{\ell}{R}\miint{Q_{R,\ell}(z_0)} [G'(|Du|)+a(z)H'(|Du|)]\, dz,
    \end{aligned}
    $$
    where $c>1$ depends on $n,L,\Delta_2(G),\nabla_2(G),\Delta_2(H)$ and $\nabla_2(H)$.
\end{proof}
The following lemma, known as the parabolic Poincar\'{e} inequality, is derived by applying Lemma \ref{lem : gluing lemma} as in the proof of \cite[Lemma 2.5]{2023_Gradient_Higher_Integrability_for_Degenerate_Parabolic_Double-Phase_Systems}.
\begin{lemma} \label{lem : semi-Parabolic Poincare inequality}
    Let $u$ be a weak solution to \eqref{eq : main equation}. Then there exists a positive constant $c=c(n,N,L,\Delta_2(G),\nabla_2(G),\Delta_2(H),\nabla_2(H))>1$ such that
\begin{align*}
        &\miint{Q_{R, \ell}(z_0)} \Phi^\theta\left(\frac{|u-u_{Q_{R, \ell}(z_0)}|}{R}\right) \,dz \\
        &\qquad\leq c \miint{Q_{R, \ell}(z_0)} \Phi^\theta\left(|D u|\right) \,dz\\
        &\qquad\quad+c\Phi^\theta\left(\frac{\ell}{R^2} \miint{Q_{R, \ell}\left(z_0\right)} \left[G'(|D u|)+a(z)H'(|D u|)\right] d z\right)
\end{align*}
    for every $Q_{R, \ell}(z_0)=B_R(x_0) \times(t_0-\ell, t_0+\ell) \subset \Omega_T$ with $R, \ell>0$ and $\theta \in(\frac{1}{i_H}, 1]$, where $\Phi$ denotes either the function $G$ or $H$.
\end{lemma}
\begin{proof}
    By \eqref{eq : subadditivity of Phi}, we have
    $$
    \begin{aligned}
    &\miint{Q_{R, \ell}(z_0)} \Phi^\theta\left(\frac{|u-u_{Q_{R, \ell}(z_0)}|}{R}\right) \,dz\leq c\miint{Q_{R, \ell}(z_0)} \Phi^\theta\left(\frac{|u-u_{B_{R}(x_0)}(t)|}{R}\right) \,dz\\
    &\qquad\qquad+c\dashint_{I_{\ell}(t_0)} \Phi^\theta\left(\frac{|u_{Q_{R, \ell}(z_0)}-u_{B_R(x_0)}(t)|}{R}\right) \,dt
    \end{aligned}
    $$
    for some constant $c=c(\Delta_2(\Phi))>1$. By \eqref{cond : ainc_1 of G and H}, the Sobolev--Poincar\'{e} inequality \cite[Corollary 7.4.3]{Harjulehto2019} is available, and it provides a constant $\beta=\beta(n,N,\nabla_2(\Phi))\in(0,1]$ such that
    $$
    \miint{Q_{R, \ell}(z_0)} \Phi^\theta\left(\frac{\beta|u-u_{B_{R}(x_0)}(t)|}{2R}\right) \,dz\leq \miint{Q_{R, \ell}(z_0)} \Phi^\theta\left(|Du|\right) \,dz.
    $$
    By \eqref{eq : Phi(Lambda t) vs Phi(t)}, we see that
    $$
    \begin{aligned}
        \miint{Q_{R, \ell}(z_0)} \Phi^\theta\left(\frac{|u-u_{B_{R}(x_0)}(t)|}{R}\right) \,dz&=\miint{Q_{R, \ell}(z_0)} \Phi^\theta\left(\frac{2}{\beta}\cdot\frac{\beta|u-u_{B_{R}(x_0)}(t)|}{2R}\right) \,dz\\
        &\leq c\miint{Q_{R, \ell}(z_0)} \Phi^\theta\left(\frac{\beta|u-u_{B_{R}(x_0)}(t)|}{2R}\right) \,dz\\
        &\leq c\miint{Q_{R, \ell}(z_0)} \Phi^\theta\left(|Du|\right) \,dz
    \end{aligned}
    $$
    for some $c=c(n,N,\nabla_2(\Phi),\Delta_2(\Phi))>1$.

    Now, by Jensen's inequality \cite[Lemma 4.3.1]{Harjulehto2019} and the above argument, there exist $\tilde{\beta}=\tilde{\beta}(N,\nabla_2(\Phi))\in(0,1)$ and $c=c(N,\tilde{L},\nabla_2(\Phi),\Delta_2(\Phi))>1$ such that
    $$
    \begin{aligned}
    &\dashint_{I_{\ell}(t_0)} \Phi^\theta\left(\frac{|u_{Q_{R, \ell}(z_0)}-u_{B_R(x_0)}(t)|}{R}\right) \,dt\\
    &\qquad\qquad\leq c\dashint_{I_{\ell}(t_0)} \Phi^\theta\left(\frac{\tilde{\beta}|u_{Q_{R, \ell}(z_0)}-u_{B_R(x_0)}(t)|}{R}\right) \,dt\\
    &\qquad\qquad\leq c\dashint_{I_{\ell}(t_0)}\dashint_{I_{\ell}(t_0)} \Phi^\theta\left(\frac{|u_{B_R (x_0)}(t_1)-u_{B_R(x_0)}(t_2)|}{R}\right) \,dt_1dt_2.
    \end{aligned}
    $$
    For $\eta\in C_0^\infty(B_R(x_0))$ satisfying \eqref{cond : cut-off function eta}, we have from \eqref{eq : subadditivity of Phi} that
    $$
    \begin{aligned}
        &\dashint_{I_{\ell}(t_0)}\dashint_{I_{\ell}(t_0)} \Phi^\theta\left(\frac{|u_{B_{R}(x_0)}(t_1)-u_{B_R(x_0)}(t_2)|}{R}\right) \,dt_1dt_2\\
        &\qquad\qquad \leq c \dashint_{I_\ell(t_0)} \Phi^\theta\left(\frac{|(u\eta)_{B_R(x_0)}(t)-u_{B_R(x_0)}(t)|}{R}\right)\,dt\\
        &\qquad\qquad\quad + c \sup_{t_1,t_2\in I_\ell(t_0)} \Phi^\theta\left(\frac{|(u\eta)_{B_R(x_0)}(t_1)-(u\eta)_{B_R(x_0)}(t_2)|}{R}\right).
    \end{aligned}
    $$
    By Lemma \ref{lem : gluing lemma}, the second term is bounded by
    \begin{multline*}
        \sup_{t_1,t_2\in I_\ell(t_0)} \Phi^\theta\left(\frac{|(u\eta)_{B_R(x_0)}(t_1)-(u\eta)_{B_R(x_0)}(t_2)|}{R}\right)\\
        \leq c\Phi^\theta\left(\frac{\ell}{R^2} \miint{Q_{R, \ell}\left(z_0\right)} \left[G'(|D u|)+a(z)H'(|D u|)\right] d z\right).
    \end{multline*}
    Moreover, Jensen's inequality as above gives
    \begin{multline*}
    \dashint_{I_\ell(t_0)} \Phi^\theta\left(\frac{|(u\eta)_{B_R(x_0)}(t)-u_{B_R(x_0)}(t)|}{R}\right)\,dt\\
    \leq c\miint{Q_{R,\ell}(z_0)}\Phi^\theta\left(\frac{|u(x,t)-u_{B_R(x_0)}(t)|}{R}\right)\,dz.
    \end{multline*}
    By applying the Sobolev--Poincar\'{e} inequality as above, we have
    $$
    \dashint_{I_\ell(t_0)} \Phi^\theta\left(\frac{|(u\eta)_{B_R(x_0)}(t)-u_{B_R(x_0)}(t)|}{R}\right)\,dt\leq c \miint{Q_{R,\ell}(z_0)} \Phi^\theta(|Du|)\, dz.
    $$
    Combining these inequalities, we obtain the conclusion.
\end{proof}
\section{Stopping time argument}\label{section 3}
In this section we use a stopping time argument to choose an intrinsic cylinder centered at almost every point of a superlevel set of $\Psi(\cdot,|Du|)$. We then distinguish three cases according to the size of $a(\cdot)$ on these cylinders and show that, under \eqref{cond : main assumption with infty} or \eqref{cond : main assumption with s}, the third case cannot occur, so that only the $G$-phase and the $(G,H)$-phase remain.

We put
\begin{equation}\label{def : lambda_0 and Lambda_0}
    \lambda_0^2\coloneq \miint{Q_{2r}(z_0)} \Psi(z,|Du|)\, dz+1\quad\text{and}\quad \Lambda_0\coloneq G(\lambda_0)+\sup_{z\in Q_{2r}(z_0)} a(z) H(\lambda_0),
\end{equation}
where $Q_{2r}(z_0)=B_{2r}(x_0)\times (t_0-(2r)^2,t_0+(2r)^2)$. Throughout Sections \ref{section 3}--\ref{section 5} we assume that
$$
2r\leq 1.
$$
All radii at which \eqref{cond : property of infty condition}, \eqref{cond : property of s} and Lemmas \ref{lem : Gap inequality with s} and \ref{lem : Gap inequality with infty} are applied in Sections \ref{section 3} and \ref{section 4} are bounded by $4\rho_w$, where $\rho_w<(r_2-r_1)/(2\kappa)\leq r/20$ is the radius constructed below, so that they are at most $r/5\leq 1$. The case $2r>1$ is reduced to the case $2r\leq 1$ at the end of Section \ref{section 5}. Let
\begin{equation}\label{def : definition of K and kappa}
    K\coloneq \begin{cases}
        1+80c_b\varpi_b\quad &\text{if \eqref{cond : main assumption with infty} holds},\\
        1+80c_s\varpi_s\quad &\text{if \eqref{cond : main assumption with s} holds},
    \end{cases}\quad\text{and}\quad \kappa\coloneq (10\tilde{L}K)^{\frac{\sigma_H+2}{i_H}}\geq 10.
\end{equation}
Here, $c_b$ is defined in Lemma \ref{lem : no occurance with infty} and $c_s$ is defined in Lemma \ref{lem : no occurance with s}. For $\Upsilon(\Lambda)$ as in \eqref{def : definition of Upsilon} and $\varrho\in[r,2r]$, we write
$$
\Upsilon(\Lambda,\varrho)\coloneq \Upsilon(\Lambda)\cap Q_\varrho(z_0).
$$
Now, we apply a stopping time argument. Let $r\leq r_1<r_2\leq 2r$ and 
$$
\Lambda>\left(\frac{4\kappa r}{r_2-r_1}\right)^\frac{\sigma_H(n+2)}{2}\Lambda_0.
$$
For any $w\in\Upsilon(\Lambda,r_1)$, we choose $\lambda_w>0$ such that
\begin{equation}\label{eq : definition of lambda_w}
\Lambda=G(\lambda_w)+a(w)H(\lambda_w)=\Psi_w(\lambda_w). 
\end{equation}
We claim that
$$
\lambda_w>\left(\frac{4\kappa r}{r_2-r_1}\right)^\frac{n+2}{2}\lambda_0.
$$
Indeed, if not, we get
$$
\begin{aligned}
    \Lambda=\Psi_w(\lambda_w)&=G(\lambda_w)+a(w)H(\lambda_w)\\
    &\leq G\left(\left(\frac{4\kappa r}{r_2-r_1}\right)^\frac{n+2}{2}\lambda_0\right)+a(w)H\left(\left(\frac{4\kappa r}{r_2-r_1}\right)^\frac{n+2}{2}\lambda_0\right).
\end{aligned}
$$
Since 
$$
\left(\frac{4\kappa r}{r_2-r_1}\right)^\frac{n+2}{2}\geq 1,
$$
we have
$$
\begin{aligned}
    \Lambda&\leq \left(\frac{4\kappa r}{r_2-r_1}\right)^{\frac{\sigma_H(n+2)}{2}}\left(G(\lambda_0)+a(w)H(\lambda_0)\right)< \Lambda,
\end{aligned}
$$
which is a contradiction. Moreover, since $\lambda_w\geq\lambda_0\geq 1$, Lemma \ref{lem : quadratic ainc} gives $\frac{\lambda_w^2}{G(\lambda_w)}\leq 1$, and therefore
$$
Q_\sigma^{\lambda_w}(w)\subset Q_{r_2}(z_0)\subset Q_{2r}(z_0)\qquad\text{for every }\sigma\in(0,r_2-r_1),
$$
because $|x-x_0|<\sigma+r_1<r_2$ and $|t-t_0|<\sigma^2+r_1^2\leq r_2^2$ for $(x,t)\in Q_\sigma^{\lambda_w}(w)$. Therefore, for any $\sigma\in[(r_2-r_1)/(2\kappa),r_2-r_1)$, we have
\begin{align}
    \miint{Q_\sigma^{\lambda_w}(w)} \Psi(z,|Du|)\, dz&\leq \frac{G(\lambda_w)}{\lambda_w^2}\left(\frac{2r}{\sigma}\right)^{n+2}\miint{Q_{2r}(z_0)} \Psi(z,|Du|)\, dz \nonumber\\\label{eq : estimate of Psi in Q_s^lambda as G}
    &\leq\frac{G(\lambda_w)}{\lambda_w^2}\left(\frac{4\kappa r}{r_2-r_1}\right)^{n+2}\lambda_0^2< G(\lambda_w).
\end{align}
By \eqref{eq : definition of lambda_w}, $w\in \Upsilon(\Lambda,r_1)$ implies $w\in \Upsilon(G(\lambda_w),r_1)$. From now on we assume in addition that $w$ is a Lebesgue point of $\Psi(\cdot,|Du|)$ with respect to the parabolic cylinders $Q_\varrho(w)$, which is the case for almost every $w\in\Upsilon(\Lambda,r_1)$. Since $\frac{\lambda_w^2}{\Psi_w(\lambda_w)}\leq\frac{\lambda_w^2}{G(\lambda_w)}\leq 1$, the means of $\Psi(\cdot,|Du|)$ over $Q_\sigma^{\lambda_w}(w)$ and over $\widetilde{Q}_\sigma^{\lambda_w}(w)$ then converge to $\Psi(w,|Du(w)|)$ as $\sigma\to0$. Hence there exists $\rho_w\in (0,(r_2-r_1)/(2\kappa))$ such that
\begin{equation}\label{cond : p-phase condition equal}
    \miint{Q_{\rho_w}^{\lambda_w}(w)} \Psi(z,|Du|)\, dz = G(\lambda_w)
\end{equation}
and
\begin{equation}\label{cond : p-phase condition leq}
    \miint{Q_{\sigma}^{\lambda_w}(w)} \Psi(z,|Du|)\, dz < G(\lambda_w)
\end{equation}
for any $\sigma\in (\rho_w,r_2-r_1)$. Since the first inequality in \eqref{eq : estimate of Psi in Q_s^lambda as G} holds for every $\sigma\in(0,r_2-r_1)$, it follows from this inequality with $\sigma=\rho_w$ and \eqref{cond : p-phase condition equal} that
\begin{equation}
    \lambda_w\leq \left(\frac{2r}{\rho_w}\right)^{\frac{n+2}{2}}\lambda_0.
\end{equation}

For $K>1$ defined in \eqref{def : definition of K and kappa}, we consider the following three cases:
\begin{enumerate}   [label= ($\arabic*$)]
    \item\label{case : p-intrinsic}  $\displaystyle KG(\lambda_w)\geq \sup_{Q_{10\rho_w}(w)}a(\cdot)H(\lambda_w)$,
    \item\label{case : p,q-intrinsic} $\displaystyle KG(\lambda_w)\leq \sup_{Q_{10\rho_w}(w)}a(\cdot)H(\lambda_w)$\qquad and\qquad $\displaystyle\sup_{Q_{10\rho_w}(w)}a(\cdot)\geq 40\omega(\rho_w)$,
    \item\label{case : no occurance} $\displaystyle KG(\lambda_w)\leq \sup_{Q_{10\rho_w}(w)}a(\cdot)H(\lambda_w)$\qquad and\qquad $\displaystyle \sup_{Q_{10\rho_w}(w)}a(\cdot)\leq 40\omega(\rho_w)$.
\end{enumerate}
\textbf{Case} \ref{case : p-intrinsic}: Using \eqref{cond : p-phase condition equal} and \eqref{cond : p-phase condition leq} and replacing $(w,\rho_w,\lambda_w)$ with $(z_0,\rho,\lambda)$, we have
\begin{equation}\label{cond : properties in p intrinsic case}
\left\{
\begin{aligned}
    &KG(\lambda)\geq \sup_{Q_{10\rho}(z_0)}a(\cdot)H(\lambda),\\
    &\miint{Q_{\rho}^{\lambda}(z_0)} \Psi(z,|Du|)\, dz = G(\lambda),\\
    &\miint{Q_{\sigma}^{\lambda}(z_0)} \Psi(z,|Du|)\, dz < G(\lambda)\quad \text{for any }\sigma\in(\rho,2\kappa\rho].
\end{aligned}
\right.
\end{equation}
\textbf{Case} \ref{case : p,q-intrinsic}: Since $d_p(z_1,z_2)\leq 20\rho_w$ for all $z_1,\,z_2\in Q_{10\rho_w}(w)$, it follows from \ref{case : p,q-intrinsic}$_2$ and \eqref{cond : property of concave function} that 
$$
40\omega(\rho_w)\leq \sup_{Q_{10\rho_w}(w)} a(\cdot)\leq \inf_{Q_{10\rho_w}(w)}a(\cdot)+20\omega(\rho_w).
$$
Then we have
$$
\sup_{Q_{10\rho_w}(w)}a(\cdot)\leq \inf_{Q_{10\rho_w}(w)}a(\cdot)+20\omega(\rho_w) \leq 2\inf_{Q_{10\rho_w}(w)}a(\cdot).
$$
Thus, we obtain
\begin{equation}\label{eq : comparison of a in p,q phase}
\frac{a(w)}{2}\leq a(\tilde{w})\leq 2a(w) \quad \text{for every }\tilde{w}\in Q_{10\rho_w}(w).
\end{equation}
Combining \ref{case : p,q-intrinsic}$_1$ and \eqref{eq : comparison of a in p,q phase}, we have $KG(\lambda_w)\leq 2a(w)H(\lambda_w)$, and hence $a(w)>0$. Since $\widetilde{Q}_\sigma^{\lambda_w}(w)\subsetneq Q_\sigma^{\lambda_w}(w)$, by \eqref{cond : p-phase condition equal} and \eqref{cond : p-phase condition leq}, we have
$$
\miint{\widetilde{Q}_\sigma^{\lambda_w}(w)} \Psi(z,|Du|)\, dz\leq \frac{\Psi_w(\lambda_w)}{G(\lambda_w)}\miint{Q_\sigma^{\lambda_w}(w)} \Psi(z,|Du|)\,dz <\Psi_w(\lambda_w)
$$
for any $\sigma\in (\rho_w,r_2-r_1)$. Since $w\in \Upsilon(\Lambda,r_1)$ and $\Lambda=\Psi_w(\lambda_w)$, we can choose $\varsigma_w\in (0,\rho_w]$ such that
\begin{equation}\label{cond : p,q-phase condition equal}
\miint{\widetilde{Q}_{\varsigma_w}^{\lambda_w}(w)} \Psi(z,|Du|)\, dz = \Psi_{w}(\lambda_w)
\end{equation}
and
\begin{equation}\label{cond : p,q-phase condition leq}
    \miint{\widetilde{Q}_{\sigma}^{\lambda_w}(w)} \Psi(z,|Du|)\, dz < \Psi_{w}(\lambda_w)
\end{equation}
for every $\sigma\in(\varsigma_w,r_2-r_1)$. Therefore, by replacing $(w,\varsigma_w,\lambda_w)$ in \eqref{eq : comparison of a in p,q phase}, \eqref{cond : p,q-phase condition equal} and \eqref{cond : p,q-phase condition leq} with $(z_0,\rho,\lambda)$, we have
\begin{equation}\label{cond : properties in p,q intrinsic case}
\left\{
\begin{aligned}
    &\frac{a(z_0)}{2}\leq a(z)\leq 2a(z_0)\quad \text{for any }z\in \widetilde{Q}^\lambda_{4\rho}(z_0),\\
    &\miint{\widetilde{Q}_{\rho}^{\lambda}(z_0)} \Psi(z,|Du|)\, dz = \Psi_{z_0}(\lambda),\\
    &\miint{\widetilde{Q}_{\sigma}^{\lambda}(z_0)} \Psi(z,|Du|)\, dz < \Psi_{z_0}(\lambda)\quad \text{for any }\sigma\in(\rho,2\kappa\rho].
\end{aligned}
\right.
\end{equation}
\textbf{Case} \ref{case : no occurance}: First, we show that this case cannot occur for bounded weak solutions.
\begin{lemma}\label{lem : no occurance with infty}
    Let $u$ be a weak solution to \eqref{eq : main equation}. If
    \begin{equation}\label{eq : sup a < omega(rho)}
        \sup_{Q_{10\rho_w}(w)} a(\cdot)\leq 40\omega(\rho_w)
    \end{equation}
    and \eqref{cond : main assumption with infty} hold, then there exists a constant $c_b=c_b(n,\nu,L,\varpi_b,\Delta_2(G),\nabla_2(G),$ $\Delta_2(H), \nabla_2(H),\|u\|_{L^\infty(\Omega_T)})>1$ such that
    $$
    \left(\frac{H}{G}\right)(\lambda_w)\leq c_b\left(\frac{H}{G}\right)(\rho_w^{-1}).
    $$
\end{lemma}
\begin{proof}
    It follows from Lemma \ref{lem : Caccioppoli inequality} and \eqref{cond : p-phase condition equal} that
    \begin{align*}
        G(\lambda_w)&=\miint{Q_{\rho_w}^{\lambda_w}(w)} \Psi(z,|Du|)\, dz\\
        &\leq c\miint{Q_{2\rho_w}^{\lambda_w}(w)}\left[G\left(\frac{\Big|u-u_{Q_{2\rho_w}^{\lambda_w}(w)}\Big|}{\rho_w}\right)+a(z)H\left(\frac{\Big|u-u_{Q_{2\rho_w}^{\lambda_w}(w)}\Big|}{\rho_w}\right)\right]\,dz\\
        &\qquad +c\frac{G(\lambda_w)}{\lambda_w^2}\miint{Q_{2\rho_w}^{\lambda_w}(w)} \frac{\Big|u-u_{Q_{2\rho_w}^{\lambda_w}(w)}\Big|^2}{\rho_w^2}\, dz=\mathrm{I}+\mathrm{II}
    \end{align*}
    for some $c=c(n,\nu,L,\Delta_2(G),\nabla_2(G),\Delta_2(H),\nabla_2(H))>1$. First, we estimate $\mathrm{I}$. Notice that
    $$
    \frac{\Big|u-u_{Q_{2\rho_w}^{\lambda_w}(w)}\Big|}{\rho_w}\leq 2\|u\|_{L^\infty(\Omega_T;\mr^N)}\rho_w^{-1}.
    $$
    Then we get from \eqref{eq : sup a < omega(rho)} that
    $$
    \begin{aligned}
        \mathrm{I}&\leq cG(\rho_w^{-1})\left[1+\omega(\rho_w)\frac{H(\rho_w^{-1})}{G(\rho_w^{-1})}\right]\leq cG(\rho_w^{-1}),
    \end{aligned}
    $$
    where $c>1$ depends on $n,\nu,L,\varpi_b,\Delta_2(G),\nabla_2(G),\Delta_2(H),\nabla_2(H)$ and $\|u\|_{L^\infty(\Omega_T)}$.

    Now, we estimate $\mathrm{II}$. Observe that
    $$
    \begin{aligned}
        \mathrm{II}&\leq c\frac{G(\lambda_w)}{\lambda_w^2}\miint{Q_{2\rho_w}^{\lambda_w}(w)}\frac{|u|^2}{\rho_w^2}\, dz\leq c\frac{G(\lambda_w)}{\lambda_w^2}\rho_w^{-2}
    \end{aligned}
    $$
    for some $c=c(n,\nu,L,\Delta_2(G),\nabla_2(G),\Delta_2(H),\nabla_2(H),\|u\|_{L^\infty(\Omega_T)})>1$. Applying Lemma \ref{lem : tilde G conjugate} with $\varepsilon=\frac{1}{2}$, $\varkappa_1=\frac{G(\lambda_w)}{\lambda_w^2}$ and $\varkappa_2=c\rho_w^{-2}$, and then using \eqref{eq : tilde G property} and the $\Delta_2$-condition of $G$, we have
    $$
    \mathrm{II}\leq \frac{1}{2} \tilde{G}^*\left(\frac{G(\lambda_w)}{\lambda_w^2}\right)+c \tilde{G}(c\rho_w^{-2})\leq \frac{1}{2}G(\lambda_w)+c G(\rho_w^{-1}).
    $$
    Combining the above inequalities, we obtain
    $$
    G(\lambda_w) \leq cG(\rho_w^{-1})
    $$
    for some $c=c(n,\nu,L,\varpi_b,\Delta_2(G),\nabla_2(G),\Delta_2(H),\nabla_2(H),\|u\|_{L^\infty(\Omega_T)})>1$. By \cite[Lemma 2.3.3]{Harjulehto2019} and the definition of $\nabla_2(G)$, we conclude that
    $$
    \lambda_w\leq c \rho_w^{-1}
    $$
    for some $c=c(n,\nu,L,\varpi_b,\Delta_2(G),\nabla_2(G),\Delta_2(H),\nabla_2(H),\|u\|_{L^\infty(\Omega_T)})>1$. Since $G\prec H$, by Lemma \ref{lem208} and the definitions of $\Delta_2(H)$ and $\nabla_2(G)$, we get
    $$
    \left(\frac{H}{G}\right)(\lambda_w)\leq c_b \left(\frac{H}{G}\right)(\rho_w^{-1}).
    $$
\end{proof}
Now, we show that Case \ref{case : no occurance} never occurs whenever \eqref{cond : main assumption with infty} is satisfied. If \ref{case : no occurance} holds, then
$$
KG(\lambda_w)\leq \sup_{Q_{10\rho_w}(w)}a(\cdot)H(\lambda_w)\leq 40\omega(\rho_w)H(\lambda_w).
$$
When \eqref{cond : main assumption with infty} holds, it follows from Lemma \ref{lem : no occurance with infty}, \eqref{cond : property of infty condition} and \eqref{def : definition of K and kappa} that
$$
\begin{aligned}
KG(\lambda_w)\leq 40\omega(\rho_w)H(\lambda_w)&=40\omega(\rho_w)\left(\frac{H}{G}\right)(\lambda_w)G(\lambda_w)\\
&\leq 40c_b\omega(\rho_w)\left(\frac{H}{G}\right)(\rho_w^{-1})G(\lambda_w)\\
&\leq 40c_b\varpi_bG(\lambda_w)\leq \frac{K}{2}G(\lambda_w),
\end{aligned}
$$
which is a contradiction.

Next, we prove that Case \ref{case : no occurance} cannot occur for weak solutions that satisfy \eqref{cond : main assumption with s}. For this, we fix 
\begin{equation}\label{eq : choice of theta}
\theta\in\left(\max\left\{\frac{1}{i_H},1-\frac{s^2}{n(n+s)\sigma_G}\right\},1\right).
\end{equation}
Then the choice of $\theta$ guarantees, for $\mu=0$ and $\mu=\frac{s}{n+s}$,
\begin{equation}\label{eq : choice of sigma_mu}
\sigma_\mu\coloneq \mu n + (1-\theta)n\sigma_G\leq s.
\end{equation}
Indeed, \eqref{eq : choice of theta} yields that
$$
\sigma_0\leq \sigma_{\frac{s}{n+s}}=\frac{ns}{n+s}+(1-\theta)n\sigma_G\leq \frac{ns}{n+s}+\frac{s^2}{n+s}=s.
$$
\begin{lemma}\label{lem : Slicewise Sobolev inequality}
    Let $B=B_{2\rho}(x_0)\subset \Omega$. Let $v\in W^{1,\sigma_G}(B)$ be nonnegative with $|Dv|\leq g$ a.e. for some $g\in L^{\sigma_G}(B)$, and $\theta$ be as in \eqref{eq : choice of theta}. Then
    $$
    \left(\dashint_B G^{\frac{\theta n}{n-1}}(v)\, dx\right)^\frac{n-1}{n}\leq c \left(\dashint_B G^\theta(2\rho g)\, dx + \dashint_B G^\theta(v)\, dx\right)
    $$
    for some $c=c(n,\Delta_2(G),\nabla_2(G),\nabla_2(H))>1$.
\end{lemma}
\begin{proof}
    Set $f\coloneq G^\theta(v)$. Since $G\in C^1((0,\infty))$ and $\varkappa G'(\varkappa)\leq \sigma_G G(\varkappa)$, and since \eqref{cond : ainc_1 of G and H} gives $\frac{G^\theta(\varkappa)}{\varkappa}\leq \tilde{L}G^\theta(1)$ for $\varkappa\leq 1$, the function $G^\theta$ is locally Lipschitz on $[0,\infty)$ with $(G^\theta)'(\varkappa)\leq \sigma_G\frac{G^\theta(\varkappa)}{\varkappa}$ for $\varkappa>0$. For $k\in\mathbb{N}$ the function $G^\theta$ is Lipschitz on $[0,k]$, and hence the chain rule gives $f_k\coloneq G^\theta(\min\{v,k\})\in W^{1,1}(B)$ with
    $$
    2\rho|Df_k|=2\rho\theta G^{\theta-1}(v)G'(v)|Dv|\chi_{\{v<k\}}\leq 2\rho\sigma_G\frac{G^\theta(v)}{v}g.
    $$
    If $2\rho g\leq v$, then $2\rho\frac{G^\theta(v)}{v}g\leq G^\theta(v)$. Otherwise, \eqref{cond : ainc_1 of G and H} gives $\frac{G^\theta(v)}{v}\leq \tilde{L}\frac{G^\theta(2\rho g)}{2\rho g}$, and hence $2\rho\frac{G^{\theta}(v)}{v}g\leq \tilde{L}G^\theta(2\rho g)$. Therefore, we have
    $$
    2\rho |Df_k|\leq c(\sigma_G,\tilde{L}) \left[G^\theta(2\rho g)+G^\theta(v)\right].
    $$
    
    By \eqref{cond : two sided power bound of Phi}, $G^\theta(\varkappa)\leq 1+\varkappa^{\sigma_G}$ for $\varkappa\geq0$, so the right-hand side and $f$ are integrable over $B$ because $g$ and $v$ belong to $L^{\sigma_G}(B)$. Letting $k\to\infty$ by dominated convergence, we obtain $f\in W^{1,1}(B)$ together with the same bound for $2\rho|Df|$. By the Sobolev inequality in $W^{1,1}(B)$,
    $$
    \begin{aligned}
    \left(\dashint_B |f|^{\frac{n}{n-1}}\, dx\right)^\frac{n-1}{n}&\leq c(n)\left(2\rho\dashint_B |Df|\, dx+\dashint_B |f|\, dx\right)\\
    &\leq c(n,\sigma_G,\tilde{L})\left(\dashint_B G^\theta(2\rho g)\, dx+\dashint_B G^\theta(v)\, dx\right).
    \end{aligned}
    $$
\end{proof}
Now, let $u$ be a weak solution to \eqref{eq : main equation} with $u\in L^\infty(0,T;L^s(\Omega,\mr^N))$. In the remainder of this section we abbreviate
$$
Q\coloneq Q_{2\rho_w}^{\lambda_w}(w),\qquad B\coloneq B_{2\rho_w}(x_w)\qquad\text{and}\qquad I\coloneq I_{2\rho_w}^{\lambda_w}(t_w),
$$
and we set
$$
u_B(t)\coloneq \dashint_B u(x,t)\,dx\qquad\text{and}\qquad v_b(z)\coloneq \frac{|u(z)-u_B(t)|}{2\rho_w}
$$
with $w=(x_w,t_w)$ and $z=(x,t)\in Q$.
\begin{lemma}\label{lem : Slicewise L^s bounds}
    There exists a constant $c=c(n,s,\|u\|_{L^\infty(0,T;L^s(\Omega,\mr^N))})>1$ such that for a.e. $t\in I$
    $$
    \dashint_{B} v_b(x,t)^s\, dx\leq c\rho_w^{-(n+s)}
    $$
    and
    $$
    \frac{|u_{B}(t)|+|u_{Q}|}{2\rho_w}\leq c \rho_w^{-\frac{n+s}{s}}.
    $$
\end{lemma}
\begin{proof}
    For a.e. $t$, we obtain
    $$
    \begin{aligned}
    |u_B(t)|\leq \left(\dashint_B |u|^s\, dx\right)^\frac{1}{s}&\leq |B|^{-\frac{1}{s}}\|u\|_{L^\infty(0,T;L^s(\Omega,\mr^N))}\\
    &\leq c(n,s,\|u\|_{L^\infty(0,T;L^s(\Omega,\mr^N))})\rho_w^{-\frac{n}{s}}
    \end{aligned}
    $$
    and so,
    $$
    |u_Q|\leq c\rho_w^{-\frac{n}{s}}.
    $$
    Next, 
    $$
    \dashint_B v_b^s\, dx\leq 2^{s-1}\rho_w^{-s}\dashint_B(|u|^s+|u_B(t)|^s)\, dx\leq c\rho_w^{-(n+s)}.
    $$
\end{proof}
\begin{lemma}\label{lem : Intrinsic interpolation estimate}
    Set
    $$
    \mathscr{X}\coloneq \miint{Q} G\left(v_b\right)\, dz
    $$
    and 
    $$
    \mathscr{Y}\coloneq \miint{Q} v_b^{\frac{s}{n+s}} G\left(v_b\right)\, dz.
    $$
    Then for every $\varepsilon\in(0,1]$,
    $$
    \mathscr{X}+\rho_w\mathscr{Y}\leq \varepsilon G(\lambda_w)+c_\varepsilon G\left(\rho_w^{-\frac{n+s}{s}}\right),
    $$
    for some $c_\varepsilon=c_\varepsilon(\varepsilon,n,s,\Delta_2(G),\nabla_2(G),\nabla_2(H),\|u\|_{L^\infty(0,T;L^s(\Omega,\mr^N))})>1$.
\end{lemma}
\begin{proof}
    Note that $\mathscr{X}<\infty$, since $u\in L^{\sigma_H}(0,T;W^{1,\sigma_H}(\Omega;\mr^N))$ and, by \eqref{cond : two sided power bound of Phi} and \eqref{cond : relation of sigma G and sigma H}, $G(\varkappa)\leq 1+\varkappa^{\sigma_H}$ for any $\varkappa>0$.

    Fix $\mu\in\{0,\frac{s}{n+s}\}$ and let $\sigma_\mu\leq s$ be as in \eqref{eq : choice of sigma_mu}. For a.e. $t\in I$, we have
    $$
    \begin{aligned}
        \dashint_{B} v_b^{\mu} G\left(v_b\right)\, dx & = \dashint_{B} G^{\theta}(v_b)v_b^{\mu} G^{1-\theta}(v_b)\, dx\\
        &\leq \left(\dashint_B G^{\frac{\theta n}{n-1}}(v_b)\, dx\right)^{\frac{n-1}{n}}\left(\dashint_B v_b^{\mu n}G^{(1-\theta)n}(v_b)\, dx\right)^\frac{1}{n}.
    \end{aligned}
    $$

    Next, consider $\mathscr{Z}(\varkappa)\coloneq \varkappa^{\mu n} G^{(1-\theta)n}(\varkappa)$ for $\varkappa>0$. Clearly, the function $\mathscr{Z}$ is non-decreasing and, by \eqref{cond : H* ct and H t with c<1}, satisfies $\mathscr{Z}(c\varkappa)\leq c^{\sigma_\mu}\mathscr{Z}(\varkappa)$ for $c\geq 1$ with $\sigma_\mu=\mu n + (1-\theta)n\sigma_G\leq s$. It follows that $\mathscr{Z}(v_b)\leq \mathscr{Z} \left(\rho_w^{-\frac{n+s}{s}}\right)\left[1+v_b^{\sigma_\mu}\rho_w^{\frac{\sigma_\mu(n+s)}{s}}\right]$. Indeed, if $v_b\leq \rho_w^{-\frac{n+s}{s}}$, since $\mathscr{Z}$ is non-decreasing, we are done, and otherwise, $$\mathscr{Z}(v_b)=\mathscr{Z}\left(v_b\rho_w^{\frac{n+s}{s}}\rho_w^{-\frac{n+s}{s}}\right)\leq v_b^{\sigma_\mu}\rho_w^{\frac{\sigma_\mu(n+s)}{s}}\mathscr{Z}(\rho_w^{-\frac{n+s}{s}}).$$
    Hence, H\"{o}lder's inequality and Lemma \ref{lem : Slicewise L^s bounds} give, for a.e. $t$,
    $$
    \dashint_B \mathscr{Z}(v_b)\, dx\leq \mathscr{Z}\left(\rho_w^{-\frac{n+s}{s}}\right)\left[1+\rho_w^{\frac{\sigma_\mu(n+s)}{s}}\left(\dashint_B v_b^s\, dx\right)^\frac{\sigma_\mu}{s}\right]\leq c \mathscr{Z}\left(\rho_w^{-\frac{n+s}{s}}\right)
    $$
    for some $c=c(n,s,\|u\|_{L^\infty(0,T;L^s(\Omega,\mr^N))})$. Thus,
    $$
    \left(\dashint_B v_b^{\mu n}G^{(1-\theta)n}(v_b)\, dx\right)^{\frac{1}{n}}\leq c \rho_w^{-\frac{\mu(n+s)}{s}}G^{1-\theta}\left(\rho_w^{-\frac{n+s}{s}}\right).
    $$

    Moreover, for a.e. $t\in I$, the function $v_b(\cdot,t)$ belongs to $W^{1,\sigma_G}(B)$ with $|Dv_b|\leq |Du|/(2\rho_w)\in L^{\sigma_G}(B)$, since $\sigma_G\leq\sigma_H$. Then Lemma \ref{lem : Slicewise Sobolev inequality} yields
    $$
    \left(\dashint_B G^{\frac{\theta n}{n-1}}(v_b)\, dx\right)^{\frac{n-1}{n}}\leq c\dashint_B G^\theta(|Du|)\, dx + c\dashint_B G^\theta(v_b)\, dx.
    $$
    Combining the above inequalities, we have
    $$
    \begin{aligned}
        \dashint_{B} v_b^{\mu} G\left(v_b\right)\, dx&\leq c\rho_w^{-\frac{\mu(n+s)}{s}}G^{1-\theta}\left(\rho_w^{-\frac{n+s}{s}}\right)\left[\dashint_B G^\theta(|Du|)\, dx + \dashint_B G^\theta(v_b)\, dx\right]
    \end{aligned}
    $$
    and so, taking the average over time and applying H\"{o}lder's inequality with exponents $\frac{1}{\theta}$ and $\frac{1}{1-\theta}$ and \eqref{cond : p-phase condition leq}, we obtain
    \begin{equation}\label{eq : estimate v^mu G to G and X}
    \begin{aligned}
    \miint{Q} v_b^\mu G(v_b)\, dz&\leq c\rho_w^{-\frac{\mu(n+s)}{s}}G^{1-\theta}\left(\rho_w^{-\frac{n+s}{s}}\right)\left[\left(\miint{Q} G(|Du|)\, dz\right)^\theta + \mathscr{X}^\theta\right]\\
    &\leq c\rho_w^{-\frac{\mu(n+s)}{s}}G^{1-\theta}\left(\rho_w^{-\frac{n+s}{s}}\right)\left[G^\theta(\lambda_w) + \mathscr{X}^\theta\right]
    \end{aligned}
    \end{equation}
    for some $c=c(n,s,\sigma_G,\tilde{L},\|u\|_{L^\infty(0,T;L^s(\Omega,\mr^N))})$.

    Finally, by Young's inequality with exponents $\frac{1}{\theta}$ and $\frac{1}{1-\theta}$, for any $\tilde{\varepsilon}\in (0,1]$,
    $$
    cG^{1-\theta}\left(\rho_w^{-\frac{n+s}{s}}\right)G^\theta(\lambda_w)\leq \tilde{\varepsilon}G(\lambda_w)+c(\tilde{\varepsilon})G\left(\rho_w^{-\frac{n+s}{s}}\right)
    $$
    and
    $$
    cG^{1-\theta}\left(\rho_w^{-\frac{n+s}{s}}\right)\mathscr{X}^\theta\leq \frac{1}{2}\mathscr{X}+cG\left(\rho_w^{-\frac{n+s}{s}}\right).
    $$
    Since $\mathscr{X}<\infty$, \eqref{eq : estimate v^mu G to G and X} with $\mu=0$ gives
    $$
    \mathscr{X}\leq 2\tilde{\varepsilon}G(\lambda_w)+c(\tilde{\varepsilon})G\left(\rho_w^{-\frac{n+s}{s}}\right).
    $$
    Next, \eqref{eq : estimate v^mu G to G and X} with $\mu=\frac{s}{n+s}$ implies that
    $$
    \rho_w\mathscr{Y}\leq c G^{1-\theta}\left(\rho_w^{-\frac{n+s}{s}}\right)\left[G^\theta(\lambda_w)+\mathscr{X}^\theta\right].
    $$
    Since $\mathscr{X}^\theta\leq (2\tilde{\varepsilon}G(\lambda_w))^{\theta}+c(\tilde{\varepsilon})G^\theta\left(\rho_w^{-\frac{n+s}{s}}\right)$, by Young's inequality, we obtain
    $$
    \rho_w\mathscr{Y}\leq c\tilde{\varepsilon} G\left(\lambda_w\right)+c(\tilde{\varepsilon})G\left(\rho_w^{-\frac{n+s}{s}}\right).
    $$
    Choosing $\tilde{\varepsilon}$ appropriately in terms of $\varepsilon$ finishes the proof.
\end{proof}
We also write $$v(z)\coloneq \frac{|u(z)-u_Q|}{2\rho_w}.$$
\begin{lemma}\label{lem : no occurance with s}
    If
    $$
    \sup_{Q_{10\rho_w}(w)} a(\cdot)\leq 40 \omega(\rho_w)
    $$
    and \eqref{cond : main assumption with s} hold, then there exists a constant $c_s=c_s(n,s,\nu,L,\varpi_s,\Delta_2(G),\nabla_2(G),$ $\Delta_2(H),\nabla_2(H),\|u\|_{L^\infty(0,T;L^s(\Omega,\mr^N))})>1$ such that
    $$
    \left(\frac{H}{G}\right)(\lambda_w)\leq c_s \left(\frac{H}{G}\right)\left(\rho_w^{-\frac{n+s}{s}}\right).
    $$
\end{lemma}
\begin{proof}
    By Lemma \ref{lem : Caccioppoli inequality} and \eqref{cond : p-phase condition equal}, we get
    $$
    \begin{aligned}
    G(\lambda_w)\leq c \miint{Q}[G(v)+a(z)H(v)]\, dz+c \frac{G(\lambda_w)}{\lambda_w^2}\miint{Q} v^2\, dz\eqcolon \mathrm{I}_G+\mathrm{I}_H+\mathrm{I}\mathrm{I}
    \end{aligned}
    $$
    for some $c=c(n,\nu,L,\Delta_2(G),\nabla_2(G),\Delta_2(H),\nabla_2(H))>1$. Lemma \ref{lem : Slicewise L^s bounds} implies $v\leq v_b+c\rho_w^{-\frac{n+s}{s}}$. Then, by \eqref{eq : subadditivity of Phi}, for $\Phi\in\{G,H\}$,
    \begin{equation}\label{eq : v leq v_b}
    \Phi(v)\leq c\Phi(v_b)+c\Phi\left(\rho_w^{-\frac{n+s}{s}}\right)\qquad \text{in }Q.
    \end{equation}
    By \eqref{eq : v leq v_b} with $\Phi=G$,
    $$
    \mathrm{I}_G\leq c \mathscr{X}+cG\left(\rho_w^{-\frac{n+s}{s}}\right).
    $$
    
    By assumption, we have $a(z)\leq 40\omega(\rho_w)$ for $z\in Q\subset Q_{2\rho_w}(w)$. Hence \eqref{eq : v leq v_b} with $\Phi=H$, Lemma \ref{lem : Gap inequality with s} and \eqref{cond : property of s} imply that
    $$
    \begin{aligned}
        \mathrm{I}_H &\leq c\omega(\rho_w)\miint{Q} H(v_b)\, dz + c\omega(\rho_w)H\left(\rho_w^{-\frac{n+s}{s}}\right)\\
        &\leq c\varpi_s \miint{Q} \left[1+v_b^{\frac{s}{n+s}}\rho_w\right]G(v_b)\, dz + c \varpi_sG\left(\rho_w^{-\frac{n+s}{s}}\right)\\
        &=c\varpi_s[\mathscr{X}+\rho_w\mathscr{Y}]+c\varpi_sG\left(\rho_w^{-\frac{n+s}{s}}\right)
    \end{aligned}
    $$
    for some $c=c(n,s,\nu,L,\Delta_2(G),\nabla_2(G),\Delta_2(H),\nabla_2(H),\|u\|_{L^\infty(0,T;L^s(\Omega,\mr^N))})>1$.

    Since $s\geq 2$, H\"{o}lder's inequality and Lemma \ref{lem : Slicewise L^s bounds} imply that
    $$
    \miint{Q} v^2\, dz\leq \left(\miint{Q} |v|^s\, dz\right)^{\frac{2}{s}}\leq c\rho_w^{-\frac{2(n+s)}{s}}
    $$
    for some $c=c(n,s,\|u\|_{L^\infty(0,T;L^s(\Omega,\mr^N))})$. By Lemma \ref{lem : tilde G conjugate} with $\varepsilon=\frac{1}{4}$, \eqref{eq : tilde G property} and the $\Delta_2$-condition of $G$, we obtain
    $$
    \mathrm{II}\leq c\frac{G(\lambda_w)}{\lambda_w^2}\rho_w^{-\frac{2(n+s)}{s}}\leq \frac{1}{4}G(\lambda_w)+cG(\rho_w^{-\frac{n+s}{s}}).
    $$
    Thus, combining the three estimates and using Lemma \ref{lem : Intrinsic interpolation estimate} give
    $$
    G(\lambda_w)\leq c G(\rho_w^{-\frac{n+s}{s}}).
    $$
    From here on, the desired result follows by arguing similarly to the final part of the proof of Lemma \ref{lem : no occurance with infty}.
\end{proof}
Now, we prove that Case \ref{case : no occurance} never occurs when \eqref{cond : main assumption with s} holds. If \ref{case : no occurance} and \eqref{cond : main assumption with s} are satisfied, we get from Lemma \ref{lem : no occurance with s}, \eqref{cond : property of s} and \eqref{def : definition of K and kappa} that
$$
\begin{aligned}
    KG(\lambda_w)\leq 40\omega(\rho_w)H(\lambda_w)&=40\omega(\rho_w)\left(\frac{H}{G}\right)(\lambda_w)G(\lambda_w)\\
    &\leq 40c_s\omega(\rho_w)\left(\frac{H}{G}\right)\left(\rho_w^{-\frac{n+s}{s}}\right)G(\lambda_w)\\
    &\leq 40c_s\varpi_sG(\lambda_w)\leq \frac{K}{2}G(\lambda_w),
\end{aligned}
$$
which is a contradiction. Consequently, under either \eqref{cond : main assumption with infty} or \eqref{cond : main assumption with s}, only Cases \ref{case : p-intrinsic} and \ref{case : p,q-intrinsic} occur, and these are treated separately in the next section.

\section{Reverse H\"{o}lder inequalities}\label{section 4}
Throughout this section we fix a point $w\in\Upsilon(\Lambda,r_1)$ from the stopping time argument of Section \ref{section 3} and, as in \eqref{cond : properties in p intrinsic case} and \eqref{cond : properties in p,q intrinsic case}, we rename $(w,\rho_w,\lambda_w)$ or $(w,\varsigma_w,\lambda_w)$ as $(z_0,\rho,\lambda)$ with $z_0=(x_0,t_0)$. Thus, in this section $z_0$ denotes a point of $\Upsilon(\Lambda,r_1)$ and not the center of the cylinder $Q_{2r}(z_0)$ fixed in Section \ref{section 3}, which plays no role until Section \ref{section 5}. In this notation \eqref{eq : definition of lambda_w} reads $\Lambda=\Psi_{z_0}(\lambda)$, and in particular $G(\lambda)\leq\Lambda$. We aim to prove the reverse H\"{o}lder inequalities separately in the intrinsic cylinders \eqref{eq : definition of p-intrinsic cylinder} and \eqref{eq : definition of p,q-intrinsic cylinder}. First, we recall the Gagliardo--Nirenberg inequality in Orlicz form from \cite[Lemma 2.13]{Hasto_2021}.
\begin{lemma}\label{lem : Gagliardo-Nirenberg inequality}
    Suppose that $\Phi:[0,\infty)\rightarrow[0,\infty)$ satisfies $\Phi(0)=0$, $\displaystyle \lim_{\varkappa\rightarrow 0^+}\Phi(\varkappa)=0$, $\displaystyle \lim_{\varkappa\rightarrow \infty} \Phi(\varkappa)=\infty$,
    $$
    (0,\infty)\ni\varkappa\mapsto \frac{\Phi(\varkappa)}{\varkappa} \text{ is almost increasing}
    $$
    and    $$(0,\infty)\ni\varkappa\mapsto \frac{\Phi(\varkappa)}{\varkappa^{q_1}} \text{ is almost decreasing}
    $$
    with constant $\tilde{L}$ for some $q_1\geq 1$. Then, for any ball $B_r\subset\mr^n$, $f\in W^{1,1}(B_r,\mr^N)$, $q_2\in [1,n)$, $\gamma,\,q_3>0$ and $\theta \in [0,1]$, we have
    $$
    \begin{aligned}
    \left(\dashint_{B_r}\Phi\left(\frac{|f|}{r}\right)^\gamma\, dx\right)^\frac{1}{\gamma}&\leq c\left(\dashint_{B_r}\left[\Phi(|Df|)^{q_2}+\Phi\left(\frac{|f|}{r}\right)^{q_2}\right] dx\right)^{\frac{\theta}{q_2}}\\
    &\qquad \times\Phi\left(\left(\dashint_{B_r}\left(\frac{|f|}{r}\right)^{q_3}\, dx\right)^{\frac{1}{q_3}}\right)^{1-\theta}
    \end{aligned}
    $$
    for some $c=c(n,\tilde{L},q_1,q_3)>0$, provided that
    $$
    \frac{1}{\gamma}\geq \theta\left(\frac{1}{q_2}-\frac{1}{n}\right)+(1-\theta)\frac{q_1}{q_3}.
    $$
\end{lemma}
\begin{remark}
    Since $\varkappa^{\sigma_H-2}=\frac{H(\varkappa)}{\varkappa^2}\cdot\frac{\varkappa^{\sigma_H}}{H(\varkappa)}$ is non-decreasing by Lemma \ref{lem : quadratic ainc} and \eqref{cond : monotonicity of Phi over powers}, we have $\sigma_H\geq 2$. Hence, by Lemma \ref{lem : quadratic ainc}, \eqref{cond : monotonicity of Phi over powers} and \eqref{cond : relation of sigma G and sigma H}, the functions $G^{\frac{1}{2}}$ and $H^{\frac{1}{2}}$ satisfy the hypotheses of Lemma \ref{lem : Gagliardo-Nirenberg inequality} with $q_1=\frac{\sigma_H}{2}\geq 1$, and so does $\Psi_{z_0}^{\frac{1}{2}}$ by Lemma \ref{lem : Psi z_0 comparability} below.
\end{remark}
We also recall a standard iteration lemma, see \cite[Lemma 6.1]{2003_Giusti_Direct_methods_in_the_calculus_of_variations}.
\begin{lemma}\label{lem : iteration lemma}
    Let $0<\rho<\tau<\infty$ and let $g:[\rho,\tau]\rightarrow[0,\infty)$ be a bounded function. Suppose that
    $$
    g(\rho_1)\leq\vartheta g(\rho_2)+\frac{A}{(\rho_2-\rho_1)^{\gamma}}+B
    $$
    is satisfied for all $0<\rho\leq \rho_1<\rho_2\leq \tau$, where $\vartheta\in(0,1)$, $A,\, B\geq 0$ and $\gamma>0$. Then there exists a positive constant $c$ depending on $\vartheta$ and $\gamma$ such that
    $$
    g(\rho)\leq c\left(\frac{A}{(\tau-\rho)^\gamma}+B\right).
    $$
\end{lemma}
\subsection{The \texorpdfstring{$G$}{G}-phase case} We assume \ref{case : p-intrinsic} and estimate the last term in Lemma \ref{lem : semi-Parabolic Poincare inequality}.

In the model case $G(\varkappa)=\varkappa^p$, $H(\varkappa)=\varkappa^q$ the corresponding
estimate is obtained by splitting $a=a^{\frac1q}a^{\frac{q-1}{q}}$ and using
$a^{\frac{q-1}{q}}\varkappa^{q-1}=(a\varkappa^q)^{1-\frac1q}$, see
\cite{2023_Gradient_Higher_Integrability_for_Degenerate_Parabolic_Double-Phase_Systems}.
The last identity is available only for power functions: in the Orlicz setting one merely has
$\frac{H(\varkappa)}{\varkappa}\leq H^{1-\frac{1}{\sigma_H}}(\varkappa)$ for $\varkappa\geq1$,
and the resulting loss $\bigl(\varkappa^{\sigma_H}/H(\varkappa)\bigr)^{1/\sigma_H}$ destroys the
cancellation of the powers of $\lambda$ in Lemma \ref{lem : G phase Poincare} below unless
$i_H=\sigma_H$. We therefore avoid splitting $a$ into powers and use instead the sharp
pointwise bounds \eqref{eq : Xi pointwise}.

We first record an elementary Jensen type inequality for quasi-concave functions.

\begin{lemma}\label{lem : quasiconcave Jensen}
    Let $\Theta:[0,\infty)\rightarrow[0,\infty)$ be non-decreasing with $\Theta(0)=0$ and such
    that $y\mapsto \Theta(y)/y$ is non-increasing on $(0,\infty)$. Then, for every measurable
    set $E$ with $0<|E|<\infty$ and every nonnegative $f\in L^1(E)$,
    $$
    \miint{E} \Theta(f)\, dz\leq 2\,\Theta\left(\miint{E} f\, dz\right).
    $$
\end{lemma}
\begin{proof}
    Set 
    $$
    m\coloneq \miint{E}f\, dz
    $$ 
    and assume $m>0$, the case $m=0$ being trivial. On
    $\{f\leq m\}$ we have $\Theta(f)\leq \Theta(m)$, while on $\{f>m\}$ the monotonicity of
    $\Theta(y)/y$ gives $\Theta(f)=\frac{\Theta(f)}{f}f\leq \frac{\Theta(m)}{m}f$. Hence
    \[
    \miint{E}\Theta(f)\, dz\leq \Theta(m)+\frac{\Theta(m)}{m}\miint{E} f\, dz=2\Theta(m). \qedhere
    \]
\end{proof}

Throughout this subsection, $\theta$ denotes an arbitrary number with
\begin{equation}\label{eq : choice of theta in G phase}
    \theta\in\left(\max\left\{\frac{1}{i_H},\,1-\frac{1}{\sigma_H}\right\},1\right],
\end{equation}
which is a non-empty range since $i_H>1$ and $\sigma_H<\infty$. Note that
\eqref{eq : choice of theta in G phase} is admissible in Lemma
\ref{lem : semi-Parabolic Poincare inequality}. We write
\begin{equation}\label{def : A in G phase}
    A\coloneq \sup_{Q_{10\rho}(z_0)}a(\cdot).
\end{equation}
The construction that follows will be used with $\Phi=G$ and $\Phi=A\,H$ in the present
subsection, and with $\Phi=\Psi_{z_0}$ in Subsection \ref{subsec : GH phase}. We call a Young
function $\Phi$ \emph{admissible} if
\begin{equation}\label{cond : admissible Phi}
    (0,\infty)\ni\varkappa\mapsto\frac{\Phi(\varkappa)}{\varkappa^{2}}
    \ \text{ is non-decreasing}
    \;\text{and}\;
    (0,\infty)\ni\varkappa\mapsto\frac{\Phi(\varkappa)}{\varkappa^{\sigma_H}}
    \ \text{ is non-increasing},
\end{equation}
and for an admissible $\Phi$ we set
\begin{equation}\label{def : Xi}
    \Xi_\Phi(y)\coloneq \frac{y}{\Phi^{-1}(y)}\quad (y>0),
    \qquad \Xi_\Phi(0)\coloneq 0,
    \qquad
    \widehat{\Xi}_\Phi(y)\coloneq \Xi_\Phi\bigl(y^{\frac1\theta}\bigr)\quad (y\geq0).
\end{equation}
Both $G$ and $H$ are admissible: the first condition in \eqref{cond : admissible Phi} is Lemma
\ref{lem : quadratic ainc}, while the second one is
\eqref{cond : monotonicity of Phi over powers} for $\Phi=H$ and follows for $\Phi=G$ from
$\frac{G(\varkappa)}{\varkappa^{\sigma_H}}=\frac{G(\varkappa)}{\varkappa^{\sigma_G}}
\varkappa^{\sigma_G-\sigma_H}$ together with
\eqref{cond : relation of sigma G and sigma H}. Clearly $c\,\Phi$ is admissible whenever
$\Phi$ is admissible and $c>0$. We abbreviate
\begin{equation}\label{def : Xi A}
    \Xi_A\coloneq \Xi_{A H}\quad (A>0),\qquad \Xi_A\coloneq 0\quad (A=0),\qquad \widehat{\Xi}_A(y)\coloneq \Xi_A\bigl(y^{\frac1\theta}\bigr),
\end{equation}
so that $\Xi_A(y)=y/H^{-1}(y/A)$ when $A>0$.

\begin{lemma}\label{lem : properties of Xi}
    Let $\Phi$ be an admissible Young function and let $\theta$ be as in
    \eqref{eq : choice of theta in G phase}. Then $\Xi_\Phi$ is non-decreasing on
    $[0,\infty)$,
    \begin{equation}\label{eq : Xi pointwise}
        \frac{\Phi(\varkappa)}{\varkappa}=\Xi_\Phi\bigl(\Phi(\varkappa)\bigr)
        \qquad\text{for every }\varkappa>0,
    \end{equation}
    \begin{equation}\label{eq : Xi quasiconcave}
        (0,\infty)\ni y\mapsto \frac{\widehat{\Xi}_\Phi(y)}{y}\ \text{ is non-increasing},
    \end{equation}
    and
    \begin{equation}\label{eq : Xi scaling}
        \widehat{\Xi}_\Phi(\epsilon y)\leq \epsilon^{\frac{1}{2\theta}}\widehat{\Xi}_\Phi(y)
        \qquad\text{for all } 0<\epsilon\leq 1 \text{ and } y>0 .
    \end{equation}
    If in addition $A>0$ and $0\leq b\leq A$, then, with $\Xi_A$ as in \eqref{def : Xi A},
    \begin{equation}\label{eq : Xi A pointwise}
        \frac{b\,H(\varkappa)}{\varkappa}\leq \Xi_{A}\bigl(b\,H(\varkappa)\bigr)
        \qquad\text{for every }\varkappa>0 .
    \end{equation}
\end{lemma}
\begin{proof}
    Writing $y=\Phi(m)$ we get $\Xi_\Phi(y)=\Phi(m)/m$, which is non-decreasing in $m$ by the
    convexity of $\Phi$ together with $\Phi(0)=0$. Since $m=\Phi^{-1}(y)$ is non-decreasing in
    $y$, the function $\Xi_\Phi$ is non-decreasing. The identity
    \eqref{eq : Xi pointwise} is immediate from $\Phi^{-1}(\Phi(\varkappa))=\varkappa$.

    We turn to \eqref{eq : Xi quasiconcave}. Substituting $\tilde{y}=y^{\frac1\theta}$ we obtain
    $$
    \frac{\widehat{\Xi}_\Phi(y)}{y}=\frac{\tilde{y}}{\Phi^{-1}(\tilde{y})\,\tilde{y}^{\theta}}
    =\frac{\tilde{y}^{1-\theta}}{\Phi^{-1}(\tilde{y})},
    $$
    so that \eqref{eq : Xi quasiconcave} amounts to $\tilde{y}\mapsto \Phi^{-1}(\tilde{y})/\tilde{y}^{1-\theta}$ being
    non-decreasing. The second condition in \eqref{cond : admissible Phi} says precisely that
    $\tilde{y}\mapsto\Phi^{-1}(\tilde{y})/\tilde{y}^{\frac{1}{\sigma_H}}$ is non-decreasing, whence
    $$
    \frac{\Phi^{-1}(\tilde{y})}{\tilde{y}^{1-\theta}}=\frac{\Phi^{-1}(\tilde{y})}{\tilde{y}^{\frac{1}{\sigma_H}}}\cdot
    \tilde{y}^{\frac{1}{\sigma_H}-(1-\theta)}
    $$
    is non-decreasing because $\theta\geq 1-\frac{1}{\sigma_H}$ by
    \eqref{eq : choice of theta in G phase}.

    Next, the first condition in \eqref{cond : admissible Phi} gives
    $\Phi(c\varkappa)\leq c^2\Phi(\varkappa)$ for $0<c\leq1$, hence
    $\Phi^{-1}(cy)\geq c^{\frac12}\Phi^{-1}(y)$, and therefore, for $0<\epsilon\leq1$,
    $$
    \widehat{\Xi}_\Phi(\epsilon y)=\frac{\epsilon^{\frac1\theta}y^{\frac1\theta}}
    {\Phi^{-1}\bigl(\epsilon^{\frac1\theta}y^{\frac1\theta}\bigr)}
    \leq \frac{\epsilon^{\frac1\theta}y^{\frac1\theta}}
    {\epsilon^{\frac{1}{2\theta}}\Phi^{-1}\bigl(y^{\frac1\theta}\bigr)}
    =\epsilon^{\frac{1}{2\theta}}\widehat{\Xi}_\Phi(y),
    $$
    which is \eqref{eq : Xi scaling}.

    Finally, let $A>0$ and $0\leq b\leq A$. We may assume $b>0$. Since
    $b\,H(\varkappa)/A\leq H(\varkappa)$ we have
    $H^{-1}\bigl(b\,H(\varkappa)/A\bigr)\leq \varkappa$, and hence
    \[
    \Xi_A\bigl(b\,H(\varkappa)\bigr)
    =\frac{b\,H(\varkappa)}{H^{-1}\!\left(\frac{b\,H(\varkappa)}{A}\right)}
    \geq \frac{b\,H(\varkappa)}{\varkappa}. \qedhere
    \]
\end{proof}

\begin{remark}\label{rem : Xi quasiconcave consequences}
    Since $\widehat{\Xi}_\Phi$ is non-decreasing, vanishes at the origin and satisfies
    \eqref{eq : Xi quasiconcave}, Lemma \ref{lem : quasiconcave Jensen} applies to it. Moreover
    \eqref{eq : Xi quasiconcave} gives $\widehat{\Xi}_\Phi(\mu y)\leq
    \mu\widehat{\Xi}_\Phi(y)$ for every $\mu\geq1$ and $y>0$.
\end{remark}

We are now ready to estimate the last term in Lemma
\ref{lem : semi-Parabolic Poincare inequality}. Observe that, by Lemma
\ref{lem : quadratic ainc} and $\lambda\geq \lambda_0\geq 1$, we have
$\frac{\lambda^2}{G(\lambda)}\leq 1$, so that
\begin{equation}\label{eq : inclusion of Q sigma lambda}
    Q_\sigma^\lambda(z_0)\subset Q_{10\rho}(z_0)\qquad\text{for every }\sigma\in[2\rho,4\rho].
\end{equation}

\begin{lemma}\label{lem : last term estimate of lemma 2.14}
    Let $u$ be a weak solution to \eqref{eq : main equation} and assume that
    $Q_{4\rho}^\lambda(z_0)\subset \Omega_T$ satisfies \eqref{cond : properties in p intrinsic case}. Then,
    for $\sigma\in[2\rho,4\rho]$ and $\theta$ as in
    \eqref{eq : choice of theta in G phase}, there exists a constant
    $c=c(\Delta_2(G),\Delta_2(H))>1$ such that
    $$
    \begin{aligned}
        &\miint{Q_\sigma^\lambda(z_0)}\left[G'(|Du|)+a(z)H'(|Du|)\right] dz\\
        &\qquad\leq c\,\widehat{\Xi}_G\!\left(\miint{Q_\sigma^\lambda(z_0)} G^\theta(|Du|)\, dz\right)
        +c\,\widehat{\Xi}_A\!\left(\miint{Q_\sigma^\lambda(z_0)} \bigl(a(z)H(|Du|)\bigr)^\theta dz\right).
    \end{aligned}
    $$
\end{lemma}
\begin{proof}
    By Lemma \ref{lem : Phi and Phi'} we have
    $\Phi'(\varkappa)\leq \Delta_2(\Phi)\frac{\Phi(\varkappa)}{\varkappa}$ for
    $\Phi\in\{G,H\}$ and $\varkappa>0$. We combine this with \eqref{eq : Xi pointwise} applied
    to $\Phi=G$ and with \eqref{eq : Xi A pointwise} applied to $b=a(z)$, which is legitimate
    because $a(z)\leq A$ on $Q_\sigma^\lambda(z_0)$ by
    \eqref{eq : inclusion of Q sigma lambda} and \eqref{def : A in G phase}. Recalling
    \eqref{def : Xi}, this yields the pointwise bound
    $$
    G'(|Du|)+a(z)H'(|Du|)
    \leq c\,\widehat{\Xi}_G\bigl(G^\theta(|Du|)\bigr)
    +c\,\widehat{\Xi}_A\Bigl(\bigl(a(z)H(|Du|)\bigr)^\theta\Bigr)
    $$
    with $c=\max\{\Delta_2(G),\Delta_2(H)\}$. By Remark
    \ref{rem : Xi quasiconcave consequences}, Lemma \ref{lem : quasiconcave Jensen} applies to
    each of the two terms and yields the assertion.
\end{proof}

\begin{lemma}\label{lem : G phase Poincare}
    Let $u$ be a weak solution to \eqref{eq : main equation} and assume that
    $Q_{4\rho}^\lambda(z_0)\subset \Omega_T$ satisfies \eqref{cond : properties in p intrinsic case}. Then,
    for $\sigma\in[2\rho,4\rho]$ and $\theta$ as in
    \eqref{eq : choice of theta in G phase}, there exists a constant
    $c=c(\operatorname{data})>1$ such that
    $$
    \miint{Q_\sigma^\lambda(z_0)}G^{\theta}\left(\frac{\bigl|u-u_{Q_\sigma^\lambda(z_0)}\bigr|}{\sigma}\right) dz
    \leq c\miint{Q_\sigma^\lambda(z_0)} \Psi^\theta(z,|Du|) \, dz.
    $$
\end{lemma}
\begin{proof}
    Throughout the proof we abbreviate $Q\coloneq Q_\sigma^\lambda(z_0)$ and put
    $$
    Y\coloneq \miint{Q} \Psi^\theta(z,|Du|)\, dz .
    $$
    We may assume $Y>0$, since otherwise $Du=0$ almost everywhere in $Q$ and both sides vanish.

    We first show that $Y\leq G^\theta(\lambda)$. Since
    $\kappa\geq 10$ we have $2\kappa\rho\geq 20\rho>\sigma$, so
    \eqref{cond : properties in p intrinsic case}$_3$ applies and gives
    $$\miint{Q}\Psi(z,|Du|)\, dz\leq G(\lambda).$$ 
    As $\theta\leq 1$, Jensen's inequality yields
    \begin{equation}\label{eq : Y leq G lambda}
        Y\leq \left(\miint{Q}\Psi(z,|Du|)\, dz\right)^\theta\leq G^\theta(\lambda).
    \end{equation}
    We write $\epsilon\coloneq Y/G^\theta(\lambda)\in(0,1]$.

    Next, we estimate the two $\widehat{\Xi}$-terms. Since $G^\theta(|Du|)\leq
    \Psi^\theta(z,|Du|)$ and $\bigl(a(z)H(|Du|)\bigr)^\theta\leq \Psi^\theta(z,|Du|)$
    pointwise, the monotonicity of $\widehat{\Xi}_G$ and $\widehat{\Xi}_A$ together with
    Lemma \ref{lem : last term estimate of lemma 2.14} gives
    $$
    \miint{Q}\left[G'(|Du|)+a(z)H'(|Du|)\right] dz
    \leq c\,\widehat{\Xi}_G(Y)+c\,\widehat{\Xi}_A(Y).
    $$
    By \eqref{eq : Xi scaling} and \eqref{eq : Y leq G lambda},
    $$
    \widehat{\Xi}_G(Y)=\widehat{\Xi}_G\bigl(\epsilon G^\theta(\lambda)\bigr)
    \leq \epsilon^{\frac{1}{2\theta}}\widehat{\Xi}_G\bigl(G^\theta(\lambda)\bigr)
    =\epsilon^{\frac{1}{2\theta}}\,\Xi_G\bigl(G(\lambda)\bigr)
    =\epsilon^{\frac{1}{2\theta}}\,\frac{G(\lambda)}{\lambda}.
    $$
    For $\widehat{\Xi}_A$ we may assume $A>0$, since otherwise $\widehat{\Xi}_A\equiv0$. This
    is the only step of the proof that uses \eqref{cond : properties in p intrinsic case}$_1$,
    which gives $A\,H(\lambda)\leq KG(\lambda)$,
    that is, $G(\lambda)/A\geq H(\lambda)/K$, and therefore the second inequality in
    \eqref{eq : quadratic scaling} applied with $c=1/K\leq1$ gives
    $$
    H^{-1}\!\left(\frac{G(\lambda)}{A}\right)\geq H^{-1}\!\left(\frac{H(\lambda)}{K}\right)
    \geq K^{-\frac12}H^{-1}\bigl(H(\lambda)\bigr)=K^{-\frac12}\lambda .
    $$
    Consequently, again by \eqref{eq : Xi scaling} and \eqref{eq : Y leq G lambda},
    $$
    \widehat{\Xi}_A(Y)\leq \epsilon^{\frac{1}{2\theta}}\,\Xi_A\bigl(G(\lambda)\bigr)
    =\epsilon^{\frac{1}{2\theta}}\,\frac{G(\lambda)}{H^{-1}\!\left(\frac{G(\lambda)}{A}\right)}
    \leq K^{\frac12}\epsilon^{\frac{1}{2\theta}}\,\frac{G(\lambda)}{\lambda}.
    $$
    Altogether we arrive at
    \begin{equation}\label{eq : G phase gluing term}
        \frac{\lambda^2}{G(\lambda)}\miint{Q}\left[G'(|Du|)+a(z)H'(|Du|)\right] dz
        \leq c\bigl(1+K^{\frac12}\bigr)\,\epsilon^{\frac{1}{2\theta}}\,\lambda .
    \end{equation}

    Finally, we apply Lemma
    \ref{lem : semi-Parabolic Poincare inequality} with $\Phi=G$, $R=\sigma$ and
    $\ell=\frac{\lambda^2}{G(\lambda)}\sigma^2$, so that $Q_{R,\ell}(z_0)=Q$ and
    $\frac{\ell}{R^2}=\frac{\lambda^2}{G(\lambda)}$. This gives
    $$
    \begin{aligned}
    \miint{Q}G^{\theta}\left(\frac{\bigl|u-u_{Q}\bigr|}{\sigma}\right) dz
    &\leq c\miint{Q} G^\theta(|Du|)\, dz\\
    &\qquad+cG^\theta\left(\frac{\lambda^2}{G(\lambda)}\miint{Q}\left[G'(|Du|)+a(z)H'(|Du|)\right] dz\right).
    \end{aligned}
    $$
    The first term on the right-hand side is bounded by $cY$. For the second one we use
    \eqref{eq : G phase gluing term}, then \eqref{eq : Phi(Lambda t) vs Phi(t)} to absorb the
    constant $c(1+K^{\frac12})\geq1$, and finally the first inequality in
    \eqref{eq : quadratic scaling} with $\epsilon^{\frac{1}{2\theta}}\leq 1$:
    $$
    G^\theta\left(c\bigl(1+K^{\frac12}\bigr)\epsilon^{\frac{1}{2\theta}}\lambda\right)
    \leq cG^\theta\left(\epsilon^{\frac{1}{2\theta}}\lambda\right)
    \leq c\left(\epsilon^{\frac{1}{2\theta}}\right)^{2\theta}G^\theta(\lambda)
    =c\,\epsilon\,G^\theta(\lambda)=cY .
    $$
    Combining the last three displays completes the proof.
\end{proof}
Similarly, we have the following lemma.
\begin{lemma}\label{lem : H phase Poincare}
    Let $u$ be a weak solution to \eqref{eq : main equation} and assume that
    $Q_{4\rho}^\lambda(z_0)\subset \Omega_T$ satisfies \eqref{cond : properties in p intrinsic case}. Then,
    for $\sigma\in[2\rho,4\rho]$ and $\theta$ as in
    \eqref{eq : choice of theta in G phase}, there exists a constant
    $c=c(\operatorname{data})>1$ such that
    $$
    \inf_{Q_\sigma^\lambda (z_0)}a(\cdot)^\theta\miint{Q_\sigma^\lambda(z_0)} H^{\theta}\left(\frac{\bigl|u-u_{Q_\sigma^\lambda(z_0)}\bigr|}{\sigma}\right) dz
    \leq c\miint{Q_\sigma^\lambda(z_0)} \Psi^\theta(z,|Du|) \, dz.
    $$
\end{lemma}
\begin{proof}
    As before, we denote $Q\coloneq Q_\sigma^\lambda(z_0)$ and assume
    $$
    Y\coloneq \miint{Q} \Psi^\theta(z,|Du|)\, dz>0.
    $$
    Recall that
    $$\miint{Q}\Psi(z,|Du|)\, dz\leq G(\lambda).$$ 
    Since $\theta\leq 1$, Jensen's inequality gives
    \begin{equation}\label{eq : Y leq G lambda in H phase}
        Y\leq \left(\miint{Q}\Psi(z,|Du|)\, dz\right)^\theta\leq G^\theta(\lambda).
    \end{equation}
    We write $\epsilon\coloneq Y/G^\theta(\lambda)\in(0,1]$. We apply Lemma
    \ref{lem : semi-Parabolic Poincare inequality} with $\Phi=H$, $R=\sigma$ and
    $\ell=\frac{\lambda^2}{G(\lambda)}\sigma^2$, so that $Q_{R,\ell}(z_0)=Q$ and
    $\frac{\ell}{R^2}=\frac{\lambda^2}{G(\lambda)}$. This gives
    $$
    \begin{aligned}
    &\inf_{Q}a(\cdot)^\theta\miint{Q} H^{\theta}\left(\frac{\bigl|u-u_{Q}\bigr|}{\sigma}\right) dz\\
    &\qquad\leq c\inf_{Q}a(\cdot)^\theta\miint{Q} H^\theta(|Du|)\, dz\\
    &\qquad\qquad+c\inf_{Q}a(\cdot)^\theta H^\theta\left(\frac{\lambda^2}{G(\lambda)}\miint{Q}\left[G'(|Du|)+a(z)H'(|Du|)\right] dz\right).
    \end{aligned}
    $$
    By the proof of Lemma \ref{lem : G phase Poincare}, \eqref{eq : Y leq G lambda in H phase} and \eqref{cond : properties in p intrinsic case}$_1$, we get
    $$
    \begin{aligned}
    &\inf_{Q}a(\cdot)^\theta\miint{Q} H^{\theta}\left(\frac{\bigl|u-u_{Q}\bigr|}{\sigma}\right) dz\\
    &\qquad\leq cY + c\epsilon \inf_{Q}a(\cdot)^\theta  H^\theta(\lambda)\leq cY,
    \end{aligned}
    $$
    and hence we obtain the conclusion.
\end{proof}
\subsubsection{\texorpdfstring{$G$}{G}-phase with \eqref{cond : main assumption with s}}
Now, we assume that \eqref{cond : main assumption with s} holds.
\begin{lemma}\label{lem : estimate S in p case with s}
    Let $u$ be a weak solution to \eqref{eq : main equation}. If $Q_{2\kappa\rho}^\lambda (z_0)\subset \Omega_T$ satisfies \eqref{cond : properties in p intrinsic case}, there exists a constant $c=c(\operatorname{data})$ such that
    $$
    S(u,Q_{2\rho}^\lambda (z_0))=\sup_{I_{2\rho}^\lambda (t_0)} \dashint_{B_{2\rho}(x_0)} \frac{|u-u_{Q_{2\rho}^\lambda (z_0)}|^2}{(2\rho)^2}\,dx\leq c\lambda^2.
    $$
\end{lemma}
\begin{proof}
    We apply Lemma \ref{lem : Caccioppoli inequality} with $R=4\rho$, $r=2\rho$, $\ell=\frac{\lambda^2}{G(\lambda)}(4\rho)^2$ and $\tau=\frac{\lambda^2}{G(\lambda)}(2\rho)^2$, so that $Q_{R,\ell}(z_0)=Q_{4\rho}^\lambda(z_0)$ and $Q_{r,\tau}(z_0)=Q_{2\rho}^\lambda(z_0)$. Since $R-r=2\rho$ and $\ell-\tau=\frac{3}{4}\ell$, and since $\frac{G(\varkappa)}{\varkappa^{\sigma_G}}$ and $\frac{H(\varkappa)}{\varkappa^{\sigma_H}}$ are non-increasing with $\sigma_G\leq \sigma_H$, we obtain
    $$
    \begin{aligned}
        &\frac{G(\lambda)}{\lambda^2}S(u,Q_{2\rho}^\lambda (z_0))\\
        &\quad\leq c\miint{Q_{4\rho}^\lambda (z_0)}\left[G\left(\frac{|u-u_{Q_{4\rho}^\lambda (z_0)}|}{4\rho}\right)+a(z)H\left(\frac{|u-u_{Q_{4\rho}^\lambda (z_0)}|}{4\rho}\right)\right] \,dz\\
        &\quad\quad + \frac{cG(\lambda)}{\lambda^2}\miint{Q_{4\rho}^\lambda (z_0)} \frac{|u-u_{Q_{4\rho}^\lambda (z_0)}|^2}{(4\rho)^2}\, dz
    \end{aligned}
    $$
    for some $c=c(n,\nu,L,\Delta_2(G),\nabla_2(G),\Delta_2(H),\nabla_2(H))>1$.
    By Lemma \ref{lem : G phase Poincare} with $\theta =1$ and \eqref{cond : properties in p intrinsic case}$_3$, we obtain
    \begin{equation}\label{eq : estimate of G of S in p phase}
    \miint{Q_{4\rho}^\lambda (z_0)} G\left(\frac{|u-u_{Q_{4\rho}^\lambda (z_0)}|}{4\rho}\right)\, dz\leq cG(\lambda),
    \end{equation}
    where $c=c(\operatorname{data}).$ Next, we note that
    \begin{equation}\label{eq : estimate on a H of S in p phase}
    \begin{aligned}
        &\miint{Q_{4\rho}^\lambda (z_0)} a(z)H\left(\frac{|u-u_{Q_{4\rho}^\lambda (z_0)}|}{4\rho}\right)\, dz\\
        &\quad\leq \inf_{Q_{4\rho}^\lambda (z_0)}a(\cdot)\miint{Q_{4\rho}^\lambda (z_0)} H\left(\frac{|u-u_{Q_{4\rho}^\lambda (z_0)}|}{4\rho}\right)\, dz\\
        &\quad\quad + 2\omega(4\rho)\miint{Q_{4\rho}^\lambda (z_0)}H\left(\frac{|u-u_{Q_{4\rho}^\lambda (z_0)}|}{4\rho}\right)\, dz.
    \end{aligned}
    \end{equation}
    Here we used that $a(z)\leq \inf_{Q_{4\rho}^\lambda(z_0)}a(\cdot)+2\omega(4\rho)$ for $z\in Q_{4\rho}^\lambda(z_0)$, which follows from \eqref{cond : property of concave function}, because $\frac{\lambda^2}{G(\lambda)}\leq 1$ implies $d_p(z_1,z_2)\leq 8\rho$ for all $z_1,\,z_2\in Q_{4\rho}^\lambda(z_0)$.
    By Lemma \ref{lem : H phase Poincare} with $\theta =1$ and \eqref{cond : properties in p intrinsic case}$_3$, we get
    $$
    \inf_{Q_{4\rho}^\lambda (z_0)}a(\cdot)\miint{Q_{4\rho}^\lambda (z_0)} H\left(\frac{|u-u_{Q_{4\rho}^\lambda (z_0)}|}{4\rho}\right)\, dz\leq c G(\lambda).
    $$
    Now, we estimate the second term on the right-hand side of \eqref{eq : estimate on a H of S in p phase}. By Lemma \ref{lem : Gap inequality with s}, we have
    \begin{equation}\label{eq : (4.17)}
    \begin{aligned}
        &\omega(4\rho)\miint{Q_{4\rho}^\lambda (z_0)} H\left(\frac{|u-u_{Q_{4\rho}^\lambda (z_0)}|}{4\rho}\right)\, dz\\
        &\qquad\leq c\miint{Q_{4\rho}^\lambda (z_0)} G\left(\frac{|u-u_{Q_{4\rho}^\lambda (z_0)}|}{4\rho}\right)\, dz\\
        &\qquad\qquad + c \miint{Q_{4\rho}^\lambda (z_0)} 4\rho \left(\frac{|u-u_{Q_{4\rho}^\lambda (z_0)}|}{4\rho}\right)^\frac{s}{n+s}G\left(\frac{|u-u_{Q_{4\rho}^\lambda (z_0)}|}{4\rho}\right)\, dz.
    \end{aligned}
    \end{equation}
    The first term has already been estimated in \eqref{eq : estimate of G of S in p phase}, so we focus on the second term. Note that
    $$
    \left(\frac{|u-u_{Q_{4\rho}^\lambda (z_0)}|}{4\rho}\right)^\frac{s}{n+s}\leq \left(\frac{|u-u_{B_{4\rho} (x_0)}(t)|}{4\rho}\right)^\frac{s}{n+s}+\left(\frac{|u_{B_{4\rho}(x_0)}(t)-u_{Q_{4\rho}^\lambda (z_0)}|}{4\rho}\right)^\frac{s}{n+s}.
    $$
    Since Lemma \ref{lem : Slicewise L^s bounds}, applied with $(w,\rho_w,\lambda_w)$ replaced by $(z_0,2\rho,\lambda)$, gives
    $$
    4\rho\left(\frac{|u_{B_{4\rho}(x_0)}(t)-u_{Q_{4\rho}^\lambda (z_0)}|}{4\rho}\right)^\frac{s}{n+s}\leq c,
    $$
    we get
    \begin{equation}\label{eq : (4.18)}
    \begin{aligned}
        &\miint{Q_{4\rho}^\lambda (z_0)} 4\rho \left(\frac{|u-u_{Q_{4\rho}^\lambda (z_0)}|}{4\rho}\right)^\frac{s}{n+s}G\left(\frac{|u-u_{Q_{4\rho}^\lambda (z_0)}|}{4\rho}\right)\, dz\\
        &\qquad \leq c \miint{Q_{4\rho}^\lambda (z_0)} G\left(\frac{|u-u_{Q_{4\rho}^\lambda (z_0)}|}{4\rho}\right)\, dz\\
        &\qquad\qquad +c\miint{Q_{4\rho}^\lambda (z_0)} 4\rho \left(\frac{|u-u_{B_{4\rho} (x_0)}(t)|}{4\rho}\right)^\frac{s}{n+s}G\left(\frac{|u-u_{Q_{4\rho}^\lambda (z_0)}|}{4\rho}\right)\, dz.
    \end{aligned}
    \end{equation}
    Note that
    $$
    G\left(\frac{|u-u_{Q_{4\rho}^\lambda (z_0)}|}{4\rho}\right)\leq c \left[G\left(\frac{|u-u_{B_{4\rho} (x_0)}(t)|}{4\rho}\right)+G\left(\frac{|u_{B_{4\rho}(x_0)}(t)-u_{Q_{4\rho}^\lambda (z_0)}|}{4\rho}\right)\right].
    $$
    By H\"{o}lder's inequality with exponents $n+s$ and $\frac{n+s}{n+s-1}$, we obtain
    $$
    \begin{aligned}
        &\miint{Q_{4\rho}^\lambda (z_0)} 4\rho \left(\frac{|u-u_{B_{4\rho} (x_0)}(t)|}{4\rho}\right)^\frac{s}{n+s}G\left(\frac{|u-u_{B_{4\rho} (x_0)}(t)|}{4\rho}\right)\, dz\\
        &\qquad \leq 4\rho\dashint_{I_{4\rho}^\lambda (t_0)}\left(\dashint_{B_{4\rho} (x_0)} \left(\frac{|u-u_{B_{4\rho} (x_0)}(t)|}{4\rho}\right)^s\, dx\right)^\frac{1}{n+s}\\
        &\qquad\qquad\times\left(\dashint_{B_{4\rho} (x_0)}G^\frac{n+s}{n+s-1}\left(\frac{|u-u_{B_{4\rho} (x_0)}(t)|}{4\rho}\right)\, dx\right)^\frac{n+s-1}{n+s}\,dt.
    \end{aligned}
    $$
    As in the proof of Lemma \ref{lem : Slicewise L^s bounds}, for a.e. $t\in I_{4\rho}^\lambda (t_0)$ we have
    $$
    4\rho\left(\dashint_{B_{4\rho} (x_0)} \left(\frac{|u-u_{B_{4\rho} (x_0)}(t)|}{4\rho}\right)^s\, dx\right)^\frac{1}{n+s}\leq c.
    $$
    Since $\frac{n+s}{n+s-1}<\frac{n}{n-1}$, by \cite[Corollary 7.4.3]{Harjulehto2019}, we have
    $$
    \begin{aligned}
    \dashint_{I_{4\rho}^\lambda (t_0)}\left(\dashint_{B_{4\rho} (x_0)}G^\frac{n+s}{n+s-1}\left(\frac{|u-u_{B_{4\rho} (x_0)}(t)|}{4\rho}\right)\, dx\right)^\frac{n+s-1}{n+s}\,dt&\leq c\miint{Q_{4\rho}^\lambda (z_0)} G(|Du|)\, dz\\
    &\leq cG(\lambda).
    \end{aligned}
    $$
    Hence, combining these inequalities implies that
    \begin{equation}\label{eq : (4.19)}
    \begin{aligned}
        &\miint{Q_{4\rho}^\lambda (z_0)} 4\rho \left(\frac{|u-u_{B_{4\rho} (x_0)}(t)|}{4\rho}\right)^\frac{s}{n+s}G\left(\frac{|u-u_{B_{4\rho} (x_0)}(t)|}{4\rho}\right)\, dz\leq c G(\lambda).
    \end{aligned}
    \end{equation}
    Also, since $\frac{G(\varkappa)}{\varkappa}$ is non-decreasing and
    $$
    \frac{|u_{B_{4\rho}(x_0)}(t)-u_{Q_{4\rho}^\lambda(z_0)}|}{4\rho}\leq \dashint_{B_{4\rho}(x_0)}\frac{|u-u_{Q_{4\rho}^\lambda(z_0)}|}{4\rho}\, dx,
    $$
    by H\"{o}lder's inequality and Jensen's inequality \cite[Lemma 4.3.1]{Harjulehto2019}, 
    \begin{equation}\label{eq : (4.20)}
    \begin{aligned}
        &\miint{Q_{4\rho}^\lambda (z_0)} 4\rho \left(\frac{|u-u_{B_{4\rho} (x_0)}(t)|}{4\rho}\right)^\frac{s}{n+s}G\left(\frac{|u_{B_{4\rho} (x_0)}(t)-u_{Q_{4\rho}^\lambda(z_0)}|}{4\rho}\right)\, dz\\
        &\quad\leq 4\rho \dashint_{I_{4\rho}^\lambda (t_0)} G\left(\frac{|u_{B_{4\rho} (x_0)}(t)-u_{Q_{4\rho}^\lambda(z_0)}|}{4\rho}\right)\dashint_{B_{4\rho}(x_0)}\left(\frac{|u-u_{B_{4\rho} (x_0)}(t)|}{4\rho}\right)^\frac{s}{n+s}\, dx \,dt\\
        &\quad \leq 4\rho \operatorname*{ess\,sup}_{t\in I_{4\rho}^\lambda (t_0)} \left(\dashint_{B_{4\rho}(x_0)}\left(\frac{|u-u_{B_{4\rho} (x_0)}(t)|}{4\rho}\right)^s\, dx\right)^\frac{1}{n+s}\\
        &\qquad \quad \times \miint{Q_{4\rho}^\lambda (z_0)} G\left(\frac{|u-u_{Q_{4\rho}^\lambda(z_0)}|}{4\rho}\right)\, dz \leq cG(\lambda).
    \end{aligned}
    \end{equation}
    Combining \eqref{eq : (4.17)}, \eqref{eq : (4.18)}, \eqref{eq : (4.19)} and \eqref{eq : (4.20)}, we have
    $$
    \omega(4\rho)\miint{Q_{4\rho}^\lambda (z_0)} H\left(\frac{|u-u_{Q_{4\rho}^\lambda (z_0)}|}{4\rho}\right)\, dz\leq c G(\lambda).
    $$
    Lastly, since $\tilde{G}$ is convex, Jensen's inequality \cite[Lemma 4.3.1]{Harjulehto2019} and \eqref{eq : estimate of G of S in p phase} give
    $$
    \begin{aligned}
        G\left(\left(\miint{Q_{4\rho}^\lambda (z_0)} \frac{|u-u_{Q_{4\rho}^\lambda (z_0)}|^2}{(4\rho)^2}\, dz\right)^{\frac{1}{2}}\right)
        &=\tilde{G}\left(\miint{Q_{4\rho}^\lambda (z_0)} \frac{|u-u_{Q_{4\rho}^\lambda (z_0)}|^2}{(4\rho)^2}\, dz\right)\\
        &\leq \miint{Q_{4\rho}^\lambda (z_0)} G\left(\frac{|u-u_{Q_{4\rho}^\lambda (z_0)}|}{4\rho}\right) dz\leq cG(\lambda).
    \end{aligned}
    $$
    Since $\frac{G(\varkappa)}{\varkappa}$ is non-decreasing, we have $cG(\lambda)\leq G(c\lambda)$, and hence
    $$
    \miint{Q_{4\rho}^\lambda (z_0)} \frac{|u-u_{Q_{4\rho}^\lambda (z_0)}|^2}{(4\rho)^2}\, dz\leq c\lambda^2.
    $$
    Combining the above inequalities, we obtain $S(u,Q_{2\rho}^\lambda (z_0))\leq c\lambda^2$.
\end{proof}
\begin{lemma}\label{lem : estimate of Psi(|u-(u)_Q|) in G intrinsic case with s}
    Let $u$ be a weak solution to \eqref{eq : main equation}. If $Q_{2\kappa\rho}^\lambda(z_0)\subset \Omega_T$ satisfies \eqref{cond : properties in p intrinsic case}, then there exists a constant $c=c(\operatorname{data})$ such that for any $\theta\in(\theta_0,1)$, where $\theta_0\coloneq\frac{n(\sigma_H(n+s)+2)}{(n+s)(n\sigma_H+4)}\in(0,1)$,
    $$
    \begin{aligned}
        &\miint{Q_{2\rho}^\lambda(z_0)} \left[G\left(\frac{|u-u_{Q_{2\rho}^\lambda(z_0)}|}{2\rho}\right)+a(z)H\left(\frac{|u-u_{Q_{2\rho}^\lambda(z_0)}|}{2\rho}\right)\right]\, dz\\
        &\quad\leq cG^{1-\theta}(\lambda) \miint{Q_{2\rho}^\lambda(z_0)}\left[G^\theta\left(\frac{|u-u_{Q_{2\rho}^\lambda(z_0)}|}{2\rho}\right)+G^\theta(|Du|)\right]\, dz\\
        &\qquad +cG^{1-\theta}(\lambda)\inf_{Q_{2\rho}^\lambda(z_0)}a(\cdot)^\theta\miint{Q_{2\rho}^\lambda(z_0)} \left[ H^\theta\left(\frac{|u-u_{Q_{2\rho}^\lambda(z_0)}|}{2\rho}\right)+H^\theta(|Du|)\right]\, dz.
    \end{aligned}
    $$
\end{lemma}
\begin{proof}
    We note that
    $$
    \begin{aligned}
        &\miint{Q_{2\rho}^\lambda(z_0)} \left[G\left(\frac{|u-u_{Q_{2\rho}^\lambda(z_0)}|}{2\rho}\right)+a(z)H\left(\frac{|u-u_{Q_{2\rho}^\lambda(z_0)}|}{2\rho}\right)\right]\, dz\\
        &\quad\leq \miint{Q_{2\rho}^\lambda(z_0)} G\left(\frac{|u-u_{Q_{2\rho}^\lambda(z_0)}|}{2\rho}\right)\, dz\\
        &\qquad + \inf_{Q_{2\rho}^\lambda(z_0)}a(\cdot)\miint{Q_{2\rho}^\lambda(z_0)} H\left(\frac{|u-u_{Q_{2\rho}^\lambda(z_0)}|}{2\rho}\right)\, dz\\
        &\qquad +c\omega(\rho) \miint{Q_{2\rho}^\lambda(z_0)} H\left(\frac{|u-u_{Q_{2\rho}^\lambda(z_0)}|}{2\rho}\right)\, dz.
    \end{aligned}
    $$
    We apply Lemma \ref{lem : Gagliardo-Nirenberg inequality} with $\Phi=G^{\frac{1}{2}}$, $q_1=\frac{\sigma_H}{2}$, $\gamma=2$, $q_2=2\theta$, $q_3=2$ and $\theta\in(\theta_0,1)\subset\left(\frac{n\sigma_H}{n\sigma_H+4},1\right)$ on the ball $B_{2\rho}(x_0)$ to $u(\cdot,t)-u_{Q_{2\rho}^\lambda(z_0)}$ for a.e. $t\in I_{2\rho}^\lambda(t_0)$ and then average over $t$. Since the spatial $L^2$ mean of $\frac{|u-u_{Q_{2\rho}^\lambda(z_0)}|}{2\rho}$ on each time slice is bounded by $S(u,Q_{2\rho}^\lambda(z_0))^{\frac12}$, this gives
    $$
    \begin{aligned}
    &\miint{Q_{2\rho}^\lambda(z_0)} G\left(\frac{|u-u_{Q_{2\rho}^\lambda(z_0)}|}{2\rho}\right)\, dz\\
    &\leq c\left(\miint{Q_{2\rho}^\lambda(z_0)} \left[G^\theta(|Du|)+G^\theta\left(\frac{|u-u_{Q_{2\rho}^\lambda(z_0)}|}{2\rho}\right)\right]\, dz\right)\left(G^{1-\theta}\left(S(u,Q_{2\rho}^\lambda(z_0))^\frac{1}{2}\right)\right).
    \end{aligned}
    $$
    Here, since $n\geq 2$ and $\sigma_H\geq 2$, we observe that
    $$
    \frac{1}{2}\geq \frac{\theta(n-2\theta)}{2\theta n}+\frac{\sigma_H(1-\theta)}{4}\quad \iff \quad \frac{n\sigma_H}{n\sigma_H+4}\leq \theta
    $$
    and
    $$
    2\theta >\frac{2n\sigma_H}{n\sigma_H+4}\geq 1.
    $$
    By Lemma \ref{lem : estimate S in p case with s}, we obtain
    $$
    \begin{aligned}
    &\miint{Q_{2\rho}^\lambda(z_0)} G\left(\frac{|u-u_{Q_{2\rho}^\lambda(z_0)}|}{2\rho}\right)\, dz\\
    &\qquad\qquad\leq cG^{1-\theta}(\lambda)\miint{Q_{2\rho}^\lambda(z_0)} \left[G^\theta(|Du|)+G^\theta\left(\frac{|u-u_{Q_{2\rho}^\lambda(z_0)}|}{2\rho}\right)\right]\, dz.
    \end{aligned}
    $$
    Next, Lemma \ref{lem : Gagliardo-Nirenberg inequality} with $\Phi=H^{\frac{1}{2}}$, $q_1=\frac{\sigma_H}{2}$, $\gamma=2$, $q_2=2\theta$, $q_3=2$ and $\theta\in\left(\frac{n\sigma_H}{n\sigma_H+4},1\right)$, applied in the same way, gives
    $$
    \begin{aligned}
    &\inf_{Q_{2\rho}^\lambda(z_0)}a(\cdot)\miint{Q_{2\rho}^\lambda(z_0)} H\left(\frac{|u-u_{Q_{2\rho}^\lambda(z_0)}|}{2\rho}\right)\, dz\\
    &\quad\leq c\inf_{Q_{2\rho}^\lambda(z_0)}a(\cdot)^\theta\miint{Q_{2\rho}^\lambda(z_0)} \left[ H^\theta(|Du|)+H^\theta\left(\frac{|u-u_{Q_{2\rho}^\lambda(z_0)}|}{2\rho}\right)\right]\, dz\\
    &\quad\qquad \times \inf_{Q_{2\rho}^\lambda(z_0)}a(\cdot)^{1-\theta}H^{1-\theta}\left(S(u,Q_{2\rho}^\lambda(z_0))^\frac{1}{2}\right).
    \end{aligned}
    $$
    By \eqref{cond : properties in p intrinsic case}$_{1}$ and Lemma \ref{lem : estimate S in p case with s}, we have
    $$
    \begin{aligned}
        &\inf_{Q_{2\rho}^\lambda(z_0)}a(\cdot)\miint{Q_{2\rho}^\lambda(z_0)} H\left(\frac{|u-u_{Q_{2\rho}^\lambda(z_0)}|}{2\rho}\right)\, dz\\
        &\leq c G^{1-\theta}(\lambda)\inf_{Q_{2\rho}^\lambda(z_0)}a(\cdot)^\theta\miint{Q_{2\rho}^\lambda(z_0)} \left[ H^\theta(|Du|)+H^\theta\left(\frac{|u-u_{Q_{2\rho}^\lambda(z_0)}|}{2\rho}\right)\right]\, dz.
    \end{aligned}
    $$
    Finally, by Lemma \ref{lem : Gap inequality with s}, we obtain
    $$
    \begin{aligned}
        &\omega(\rho) \miint{Q_{2\rho}^\lambda(z_0)} H\left(\frac{|u-u_{Q_{2\rho}^\lambda(z_0)}|}{2\rho}\right)\, dz\\
        &\quad \leq c \miint{Q_{2\rho}^\lambda(z_0)} G\left(\frac{|u-u_{Q_{2\rho}^\lambda(z_0)}|}{2\rho}\right)\, dz\\
        &\qquad + c \miint{Q_{2\rho}^\lambda(z_0)} \rho \left(\frac{|u-u_{Q_{2\rho}^\lambda(z_0)}|}{2\rho}\right)^\frac{s}{n+s}G\left(\frac{|u-u_{Q_{2\rho}^\lambda(z_0)}|}{2\rho}\right)\, dz.
    \end{aligned}
    $$
    The first term on the right-hand side has already been estimated above. To estimate the second term, we abbreviate $B\coloneq B_{2\rho}(x_0)$, $I\coloneq I_{2\rho}^\lambda(t_0)$ and
    $$
    v\coloneq \frac{|u-u_{Q_{2\rho}^\lambda(z_0)}|}{2\rho},
    $$
    and we argue slicewise in time. By Lemma \ref{lem : Slicewise L^s bounds}, applied with $(w,\rho_w,\lambda_w)$ replaced by $(z_0,\rho,\lambda)$, for a.e. $t\in I$ we have
    $$
    v\leq \frac{|u-u_{B}(t)|}{2\rho}+\frac{|u_{B}(t)-u_{Q_{2\rho}^\lambda(z_0)}|}{2\rho}\leq \frac{|u-u_{B}(t)|}{2\rho}+c\rho^{-\frac{n+s}{s}}
    $$
    and
    $$
    \dashint_{B}\left(\frac{|u-u_{B}(t)|}{2\rho}\right)^s dx\leq c\rho^{-(n+s)},
    $$
    and hence
    \begin{equation}\label{eq : slicewise L^s bound of v in G phase}
    \rho\left(\dashint_{B} v(x,t)^s\, dx\right)^{\frac{1}{n+s}}\leq c\rho\left(\rho^{-(n+s)}\right)^{\frac{1}{n+s}}= c
    \end{equation}
    for a.e. $t\in I$, where $c=c(n,s,\|u\|_{L^\infty(0,T;L^s(\Omega,\mr^N))})$. Then H\"{o}lder's inequality with exponents $n+s$ and $\frac{n+s}{n+s-1}$ in the spatial variable and \eqref{eq : slicewise L^s bound of v in G phase} give, for a.e. $t\in I$,
    $$
    \begin{aligned}
         \dashint_{B} \rho\, v^\frac{s}{n+s}G(v)\, dx
         &\leq \rho\left(\dashint_{B} v^s\, dx\right)^{\frac{1}{n+s}}\left(\dashint_{B}G^{\frac{n+s}{n+s-1}}(v)\, dx\right)^{\frac{n+s-1}{n+s}}\\
         &\leq c\left(\dashint_{B}G^{\frac{n+s}{n+s-1}}(v)\, dx\right)^{\frac{n+s-1}{n+s}}.
    \end{aligned}
    $$
    Note that
    $$
    \frac{n(\sigma_H(n+s)+2)}{(n+s)(n\sigma_H+4)}=\frac{n^2\sigma_H+n\sigma_H s+ 2n}{n^2\sigma_H+n\sigma_H s + 4n+ 4s}<1
    $$
    and
    $$
    \frac{2n(\sigma_H(n+s)+2)}{(n+s)(n\sigma_H+4)}=\frac{2n^2\sigma_H+2n\sigma_H s+ 4n}{n^2\sigma_H+n\sigma_H s + 4n+ 4s}>1.
    $$
    Thus, choosing $\Phi=G^{\frac{1}{2}}$, $\gamma=\frac{2(n+s)}{n+s-1}$, $q_1=\frac{\sigma_H}{2}$, $q_2=2\theta$, $q_3=2$ and $\theta \in \left(\frac{n(\sigma_H(n+s)+2)}{(n+s)(n\sigma_H+4)},1\right)$, we see that
    $$
    \frac{n+s-1}{2(n+s)}\geq \theta\frac{n-2\theta}{2\theta n}+(1-\theta)\frac{\sigma_H}{4}\quad \iff\quad \frac{n(\sigma_H(n+s)+2)}{(n+s)(n\sigma_H+4)}\leq \theta
    $$
    is satisfied. Therefore, applying Lemma \ref{lem : Gagliardo-Nirenberg inequality} on the ball $B$ to $f=u(\cdot,t)-u_{Q_{2\rho}^\lambda(z_0)}$ with $r=2\rho$ for a.e. $t\in I$, and then using $\big(\dashint_B v(x,t)^2\, dx\big)^{\frac{1}{2}}\leq S(u,Q_{2\rho}^\lambda(z_0))^{\frac{1}{2}}$ and Lemma \ref{lem : estimate S in p case with s}, we obtain
    $$
    \begin{aligned}
        \left(\dashint_{B}G^{\frac{n+s}{n+s-1}}(v)\, dx\right)^{\frac{n+s-1}{n+s}}
        &\leq c\left(\dashint_{B} \left[G^\theta(|Du|)+G^\theta(v)\right] dx\right) G^{1-\theta}\left(S(u,Q_{2\rho}^\lambda(z_0))^\frac{1}{2}\right)\\
        &\leq c G^{1-\theta}(\lambda)\dashint_{B} \left[G^\theta(|Du|)+G^\theta(v)\right] dx
    \end{aligned}
    $$
    for a.e. $t\in I$. Combining the last two displays and taking the average over $t\in I$, we conclude that
    $$
    \miint{Q_{2\rho}^\lambda(z_0)} \rho\, v^\frac{s}{n+s}G(v)\, dz\leq c G^{1-\theta}(\lambda)\miint{Q_{2\rho}^\lambda(z_0)} \left[G^\theta(|Du|)+G^\theta\left(\frac{|u-u_{Q_{2\rho}^\lambda(z_0)}|}{2\rho}\right)\right] dz.
    $$
    Combining the above results yields the desired conclusion.
\end{proof}
\begin{lemma}\label{lem : estimate of Psi in G phase with s}
    Let $u$ be a weak solution to \eqref{eq : main equation}. Then there exist constants $c=c(\operatorname{data})$ and $\theta_0=\theta_0(n,i_H,\sigma_H,s)$ such that if $Q_{2\kappa\rho}^\lambda(z_0)\subset \Omega_T$ satisfies \eqref{cond : properties in p intrinsic case},
    $$
    \miint{Q_{\rho}^\lambda(z_0)} \Psi(z,|Du|)\, dz \leq c\left(\miint{Q_{2\rho}^\lambda(z_0)} \Psi^\theta(z,|Du|)\, dz\right)^\frac{1}{\theta}
    $$
    for any $\theta\in (\theta_0,1)$.
\end{lemma}
\begin{proof}
    We use Lemma \ref{lem : Caccioppoli inequality} and get
    $$
    \begin{aligned}
        &\miint{Q_\rho^\lambda(z_0)} \Psi(z,|Du|)\, dz\\ &\qquad\leq c \miint{Q_{2\rho}^\lambda (z_0)}\left[G\left(\frac{|u-u_{Q_{2\rho}^\lambda(z_0)}|}{2\rho}\right)+a(z)H\left(\frac{|u-u_{Q_{2\rho}^\lambda(z_0)}|}{2\rho}\right)\right]\, dz\\
        &\qquad\qquad +c\frac{G(\lambda)}{\lambda^2}\miint{Q_{2\rho}^\lambda(z_0)} \frac{|u-u_{Q_{2\rho}^\lambda(z_0)}|^2}{(2\rho)^2}\, dz.
    \end{aligned}
    $$
    By Lemma \ref{lem : estimate of Psi(|u-(u)_Q|) in G intrinsic case with s}, we have
    $$
    \begin{aligned}
        &\miint{Q_{2\rho}^\lambda (z_0)}\left[G\left(\frac{|u-u_{Q_{2\rho}^\lambda(z_0)}|}{2\rho}\right)+a(z)H\left(\frac{|u-u_{Q_{2\rho}^\lambda(z_0)}|}{2\rho}\right)\right]\, dz\\
        &\qquad \leq c G^{1-\theta}(\lambda) \miint{Q_{2\rho}^\lambda(z_0)} \Psi^\theta(z,|Du|) \, dz\\
        &\qquad\qquad + c G^{1-\theta}(\lambda) \miint{Q_{2\rho}^\lambda(z_0)} G^\theta\left(\frac{|u-u_{Q_{2\rho}^\lambda(z_0)}|}{2\rho}\right)\, dz\\
        &\qquad\qquad + c G^{1-\theta}(\lambda) \inf_{Q_{2\rho}^\lambda(z_0)}a(\cdot)^\theta\miint{Q_{2\rho}^\lambda(z_0)} H^\theta\left(\frac{|u-u_{Q_{2\rho}^\lambda(z_0)}|}{2\rho}\right)\, dz
    \end{aligned}
    $$
    for some $c=c(\operatorname{data})$ and $\theta\in (\theta_0,1)$, where $\theta_0\coloneq\max\{\theta_0',\frac{1}{i_H},1-\frac{1}{\sigma_H}\}\in(0,1)$ and $\theta_0'$ denotes the constant $\theta_0$ of Lemma \ref{lem : estimate of Psi(|u-(u)_Q|) in G intrinsic case with s}. Since such $\theta$ satisfies \eqref{eq : choice of theta in G phase}, Lemmas \ref{lem : G phase Poincare} and \ref{lem : H phase Poincare} apply with $\sigma=2\rho$ and give
    $$
    \begin{aligned}
        &\miint{Q_{2\rho}^\lambda (z_0)}\left[G\left(\frac{|u-u_{Q_{2\rho}^\lambda(z_0)}|}{2\rho}\right)+a(z)H\left(\frac{|u-u_{Q_{2\rho}^\lambda(z_0)}|}{2\rho}\right)\right]\, dz\\
        &\qquad \leq c G^{1-\theta}(\lambda)\miint{Q_{2\rho}^\lambda(z_0)} \Psi^\theta(z,|Du|)\, dz.
    \end{aligned}
    $$
    Next, we observe from Lemma \ref{lem : estimate S in p case with s} that 
    $$
    \begin{aligned}
    \miint{Q_{2\rho}^\lambda(z_0)} \frac{|u-u_{Q_{2\rho}^\lambda(z_0)}|^2}{(2\rho)^2}\, dz &\leq \left(S(u, Q_{2\rho}^\lambda(z_0))\right)^{\frac{1}{2}}\left(\miint{Q_{2\rho}^\lambda(z_0)} \frac{|u-u_{Q_{2\rho}^\lambda(z_0)}|^2}{(2\rho)^2}\, dz\right)^\frac{1}{2}\\
    &\leq c\lambda\left(\miint{Q_{2\rho}^\lambda(z_0)} \frac{|u-u_{Q_{2\rho}^\lambda(z_0)}|^2}{(2\rho)^2}\, dz\right)^\frac{1}{2}.
    \end{aligned}
    $$
    Then it follows from \eqref{eq : comparison with Phi and Phi'}, \eqref{Youngs-ineq} with $\Phi=G^*$ and \eqref{eq : Phi*(Phi')} that
    $$
    \begin{aligned}
        \frac{G(\lambda)}{\lambda^2} \miint{Q_{2\rho}^\lambda(z_0)} \frac{|u-u_{Q_{2\rho}^\lambda(z_0)}|^2}{(2\rho)^2}\, dz&\leq c\frac{G(\lambda)}{\lambda}\left(\miint{Q_{2\rho}^\lambda(z_0)} \frac{|u-u_{Q_{2\rho}^\lambda(z_0)}|^2}{(2\rho)^2}\, dz\right)^\frac{1}{2}\\
        &\leq \frac{1}{4}G(\lambda)+c\tilde{G}\left(\miint{Q_{2\rho}^\lambda(z_0)} \frac{|u-u_{Q_{2\rho}^\lambda(z_0)}|^2}{(2\rho)^2}\, dz\right).
    \end{aligned}
    $$
    By Jensen's inequality \cite[Lemma 4.3.1]{Harjulehto2019} and the above result, we have
    $$
    \begin{aligned}
        \tilde{G}\left(\miint{Q_{2\rho}^\lambda(z_0)} \frac{|u-u_{Q_{2\rho}^\lambda(z_0)}|^2}{(2\rho)^2}\, dz\right)&\leq c\miint{Q_{2\rho}^\lambda(z_0)} G\left(\frac{|u-u_{Q_{2\rho}^\lambda(z_0)}|}{2\rho}\right)\, dz\\
        &\leq c G^{1-\theta}(\lambda)\miint{Q_{2\rho}^\lambda(z_0)} \Psi^\theta(z,|Du|)\, dz.
    \end{aligned}
    $$
    Combining the above inequalities and using Young's inequality, we have
    $$
    \begin{aligned}
        \miint{Q_\rho^\lambda(z_0)} \Psi(z,|Du|)\, dz&\leq \frac{1}{4}G(\lambda) + c G^{1-\theta}(\lambda)\miint{Q_{2\rho}^\lambda(z_0)} \Psi^\theta(z,|Du|)\, dz\\
        &\leq \frac{1}{2}G(\lambda)+c\left(\miint{Q_{2\rho}^\lambda(z_0)} \Psi^\theta(z,|Du|)\, dz\right)^\frac{1}{\theta}.
    \end{aligned}
    $$
    By \eqref{cond : properties in p intrinsic case}$_2$, we obtain the desired conclusion.
\end{proof}
\begin{lemma}\label{lem : estimate of psi in level set with s}
    Let $u$ be a weak solution to \eqref{eq : main equation}. If $Q_{2\kappa \rho}^\lambda(z_0)\subset\Omega_T$ satisfies \eqref{cond : properties in p intrinsic case}, then there exist constants $c=c(\operatorname{data})$ and $\theta_0=\theta_0(n,i_H,\sigma_H,s)\in (0,1)$ such that for any $\theta\in (\theta_0,1)$,
    $$
    \iints{Q_{2\kappa\rho}^\lambda(z_0)}\Psi(z,|Du|)\, dz \leq c\Lambda^{1-\theta}\iints{Q_{2\rho}^\lambda(z_0)\cap \Upsilon(c^{-1}\Lambda)} \Psi^\theta(z,|Du|)\, dz.
    $$
\end{lemma}
\begin{proof}
    By \eqref{cond : properties in p intrinsic case}$_3$, we have
    $$
    \left(\miint{Q_{2\rho}^\lambda(z_0)}\Psi^\theta(z,|Du|)\,dz\right)^\frac{1}{\theta}\leq G^{1-\theta}(\lambda)\miint{Q_{2\rho}^\lambda(z_0)}\Psi^\theta(z,|Du|)\, dz.
    $$
    Since $\kappa>1$,
    $$
    \begin{aligned}
        &G^{1-\theta}(\lambda)\miint{Q_{2\rho}^\lambda(z_0)}\Psi^\theta(z,|Du|)\, dz\\
        &\qquad\leq \frac{G^{1-\theta}(\lambda)}{|Q_{2\rho}^\lambda|}\iints{Q_{2\rho}^\lambda(z_0)\setminus \Upsilon((2\kappa^{n+2}c)^{-\frac{1}{\theta}}G(\lambda))}\Psi^\theta(z,|Du|)\, dz\\
        &\qquad\qquad+\frac{G^{1-\theta}(\lambda)}{|Q_{2\rho}^\lambda|}\iints{Q_{2\rho}^\lambda(z_0)\cap \Upsilon((2\kappa^{n+2}c)^{-\frac{1}{\theta}}G(\lambda))}\Psi^\theta(z,|Du|)\, dz\\
        &\qquad\leq \frac{1}{2c} G(\lambda)+\frac{G^{1-\theta}(\lambda)}{|Q_{2\rho}^\lambda|}\iints{Q_{2\rho}^\lambda(z_0)\cap \Upsilon((2\kappa^{n+2}c)^{-\frac{1}{\theta}}G(\lambda))}\Psi^\theta(z,|Du|)\, dz
    \end{aligned}
    $$
    for any $c>0$. By Lemma \ref{lem : estimate of Psi in G phase with s}, we choose $\tilde{c}=\tilde{c}(\operatorname{data})$ and $\theta_0=\theta_0(n,i_H,\sigma_H,s)\in (0,1)$ such that for $\theta\in (\theta_0,1)$,
    $$
    \begin{aligned}
    \miint{Q_\rho^\lambda(z_0)} \Psi(z,|Du|)\, dz&\leq \tilde{c}\left(\miint{Q_{2\rho}^\lambda(z_0)}\Psi^\theta(z,|Du|)\, dz\right)^\frac{1}{\theta}\\
    &\leq \frac{1}{2} G(\lambda)\\
    &\quad+\frac{\tilde{c}G^{1-\theta}(\lambda)}{|Q_{2\rho}^\lambda|}\iints{Q_{2\rho}^\lambda(z_0)\cap \Upsilon((2\kappa^{n+2}\tilde{c})^{-\frac{1}{\theta}}G(\lambda))}\Psi^\theta(z,|Du|)\, dz.
    \end{aligned}
    $$
    By \eqref{cond : properties in p intrinsic case}$_2$ and \eqref{cond : properties in p intrinsic case}$_3$, for $\theta\in(\theta_0,1)$, we have
    $$
    \begin{aligned}
    \miint{Q_{2\kappa\rho}^\lambda(z_0)} \Psi(z,|Du|)\, dz&\leq 2\tilde{c}\frac{G^{1-\theta}(\lambda)}{|Q_{2\rho}^\lambda|}\iints{Q_{2\rho}^\lambda(z_0)\cap \Upsilon((2\kappa^{n+2}\tilde{c})^{-\frac{1}{\theta}}G(\lambda))}\Psi^\theta(z,|Du|)\, dz
    \end{aligned}
    $$
    and hence
    $$
    \iints{Q_{2\kappa\rho}^\lambda(z_0)}\Psi(z,|Du|)\, dz\leq 2\kappa^{n+2}\tilde{c} G^{1-\theta}(\lambda)\iints{Q_{2\rho}^\lambda(z_0)\cap \Upsilon((2\kappa^{n+2}\tilde{c})^{-\frac{1}{\theta}}G(\lambda))}\Psi^\theta(z,|Du|)\, dz.
    $$
    Since $\Lambda=\Psi_{z_0}(\lambda)\leq (1+K)G(\lambda)\leq 2KG(\lambda)$ by \eqref{cond : properties in p intrinsic case}$_1$, for $\theta\in (\theta_0,1)$,
    $$
    \frac{G(\lambda)}{(2\kappa^{n+2}\tilde{c})^\frac{1}{\theta}}\geq\frac{G(\lambda)}{(2\kappa^{n+2}\tilde{c})^\frac{1}{\theta_0}}\geq \frac{\Lambda}{2K(2\kappa^{n+2}\tilde{c})^\frac{1}{\theta_0}}, 
    $$
    and so,
    $$
    \Upsilon((2\kappa^{n+2}\tilde{c})^{-\frac{1}{\theta}}G(\lambda))\subset \Upsilon((2K(2\kappa^{n+2}\tilde{c})^\frac{1}{\theta_0})^{-1}\Lambda).
    $$
    Note that $2\kappa^{n+2}\tilde{c}\leq 2K(2\kappa^{n+2}\tilde{c})^\frac{1}{\theta_0}$ and $G(\lambda)\leq\Lambda$. Writing $c\coloneq 2K(2\kappa^{n+2}\tilde{c})^\frac{1}{\theta_0}$, we have
    $$
    \iints{Q_{2\kappa\rho}^\lambda(z_0)}\Psi(z,|Du|)\, dz\leq c \Lambda^{1-\theta}\iints{Q_{2\rho}^\lambda(z_0)\cap \Upsilon(c^{-1}\Lambda)}\Psi^\theta(z,|Du|)\, dz.
    $$
\end{proof}
\subsubsection{\texorpdfstring{$G$}{G}-phase with \eqref{cond : main assumption with infty}}
Now, we assume that \eqref{cond : main assumption with infty} holds.
\begin{lemma}\label{lem : estimate S in p case with infty}
    Let $u$ be a weak solution to \eqref{eq : main equation}. If $Q_{2\kappa\rho}^\lambda (z_0)\subset \Omega_T$ satisfies \eqref{cond : properties in p intrinsic case}, there exists a constant $c=c(\operatorname{data})$ such that
    $$
    S(u,Q_{2\rho}^\lambda (z_0))=\sup_{I_{2\rho}^\lambda (t_0)} \dashint_{B_{2\rho}(x_0)} \frac{|u-u_{Q_{2\rho}^\lambda (z_0)}|^2}{(2\rho)^2}\,dx\leq c\lambda^2.
    $$
\end{lemma}
\begin{proof}
    As in the proof of Lemma \ref{lem : estimate S in p case with s}, we have
    \begin{equation}\label{eq : estimate of G in G phase with infty}
    \miint{Q_{4\rho}^\lambda (z_0)} G\left(\frac{|u-u_{Q_{4\rho}^\lambda (z_0)}|}{4\rho}\right)\, dz\leq cG(\lambda)
    \end{equation}
    and
    $$
    \frac{G(\lambda)}{\lambda^2}S(u,Q_{2\rho}^\lambda (z_0))\leq cG(\lambda)+c\,\omega(4\rho)\miint{Q_{4\rho}^\lambda (z_0)}H\left(\frac{|u-u_{Q_{4\rho}^\lambda (z_0)}|}{4\rho}\right) \,dz,
    $$
    where $c=c(\operatorname{data}).$
    Since
    \begin{equation}\label{eq : bound of u-(u)}
    |u-u_{Q_{4\rho}^\lambda (z_0)}|\leq 2\|u\|_{L^\infty(\Omega_T)},
    \end{equation}
    by \eqref{cond : main assumption with infty}, \eqref{eq : estimate of G in G phase with infty} and Lemma \ref{lem : Gap inequality with infty}, we get
    $$
    \begin{aligned}
    &\omega(4\rho)\miint{Q_{4\rho}^\lambda (z_0)}H\left(\frac{|u-u_{Q_{4\rho}^\lambda (z_0)}|}{4\rho}\right) \,dz\\
    &\qquad\leq c\miint{Q_{4\rho}^\lambda (z_0)}G\left(\frac{|u-u_{Q_{4\rho}^\lambda (z_0)}|}{4\rho}\right) \,dz\\
    &\qquad\qquad+ c\rho\miint{Q_{4\rho}^\lambda (z_0)}\left(\frac{|u-u_{Q_{4\rho}^\lambda (z_0)}|}{4\rho}\right)G\left(\frac{|u-u_{Q_{4\rho}^\lambda (z_0)}|}{4\rho}\right) \,dz\\
    &\qquad\leq cG(\lambda).
    \end{aligned}
    $$
    Combining the above inequalities, we obtain $S(u,Q_{2\rho}^\lambda (z_0))\leq c\lambda^2$.
\end{proof}
\begin{lemma}\label{lem : estimate of Psi(|u-(u)_Q|) in G intrinsic case with infty}
    Let $u$ be a weak solution to \eqref{eq : main equation}. If $Q_{2\kappa\rho}^\lambda(z_0)\subset \Omega_T$ satisfies \eqref{cond : properties in p intrinsic case}, then there exists a constant $c=c(\operatorname{data})$ such that for any $\theta\in(\theta_0,1)$, where $\theta_0\coloneq\frac{n\sigma_H}{n\sigma_H+4}$,
    $$
    \begin{aligned}
        &\miint{Q_{2\rho}^\lambda(z_0)} \left[G\left(\frac{|u-u_{Q_{2\rho}^\lambda(z_0)}|}{2\rho}\right)+a(z)H\left(\frac{|u-u_{Q_{2\rho}^\lambda(z_0)}|}{2\rho}\right)\right]\, dz\\
        &\quad\leq cG^{1-\theta}(\lambda) \miint{Q_{2\rho}^\lambda(z_0)}\left[G^\theta\left(\frac{|u-u_{Q_{2\rho}^\lambda(z_0)}|}{2\rho}\right)+G^\theta(|Du|)\right]\, dz\\
        &\qquad +cG^{1-\theta}(\lambda)\inf_{Q_{2\rho}^\lambda(z_0)}a(\cdot)^\theta\miint{Q_{2\rho}^\lambda(z_0)} \left[ H^\theta\left(\frac{|u-u_{Q_{2\rho}^\lambda(z_0)}|}{2\rho}\right)+H^\theta(|Du|)\right]\, dz.
    \end{aligned}
    $$
\end{lemma}
\begin{proof}
    We note that
    $$
    \begin{aligned}
        &\miint{Q_{2\rho}^\lambda(z_0)} \left[G\left(\frac{|u-u_{Q_{2\rho}^\lambda(z_0)}|}{2\rho}\right)+a(z)H\left(\frac{|u-u_{Q_{2\rho}^\lambda(z_0)}|}{2\rho}\right)\right]\, dz\\
        &\quad\leq \miint{Q_{2\rho}^\lambda(z_0)} G\left(\frac{|u-u_{Q_{2\rho}^\lambda(z_0)}|}{2\rho}\right)\, dz\\
        &\qquad + \inf_{Q_{2\rho}^\lambda(z_0)}a(\cdot)\miint{Q_{2\rho}^\lambda(z_0)} H\left(\frac{|u-u_{Q_{2\rho}^\lambda(z_0)}|}{2\rho}\right)\, dz\\
        &\qquad +c\omega(\rho) \miint{Q_{2\rho}^\lambda(z_0)} H\left(\frac{|u-u_{Q_{2\rho}^\lambda(z_0)}|}{2\rho}\right)\, dz.
    \end{aligned}
    $$
    By applying the same argument as in the proof of Lemma \ref{lem : estimate of Psi(|u-(u)_Q|) in G intrinsic case with s}, with Lemma \ref{lem : estimate S in p case with infty} in place of Lemma \ref{lem : estimate S in p case with s}, we obtain
    \begin{equation}\label{eq : (4.23)}
    \begin{aligned}
    &\miint{Q_{2\rho}^\lambda(z_0)} G\left(\frac{|u-u_{Q_{2\rho}^\lambda(z_0)}|}{2\rho}\right)\, dz\\
    &\qquad\qquad\qquad \leq c G^{1-\theta}(\lambda)\miint{Q_{2\rho}^\lambda(z_0)} \left[G^\theta(|Du|)+G^\theta\left(\frac{|u-u_{Q_{2\rho}^\lambda(z_0)}|}{2\rho}\right)\right]\, dz
    \end{aligned}
    \end{equation}
    and
    $$
    \begin{aligned}
    &\inf_{Q_{2\rho}^\lambda(z_0)}a(\cdot)\miint{Q_{2\rho}^\lambda(z_0)} H\left(\frac{|u-u_{Q_{2\rho}^\lambda(z_0)}|}{2\rho}\right)\, dz\\
    &\qquad \leq c\inf_{Q_{2\rho}^\lambda(z_0)}a(\cdot)^\theta\miint{Q_{2\rho}^\lambda(z_0)} \left[ H^\theta(|Du|)+H^\theta\left(\frac{|u-u_{Q_{2\rho}^\lambda(z_0)}|}{2\rho}\right)\right]\, dz\\
    &\quad\qquad \times \inf_{Q_{2\rho}^\lambda(z_0)}a(\cdot)^{1-\theta}H^{1-\theta}\left(S(u,Q_{2\rho}^\lambda(z_0))^\frac{1}{2}\right)
    \end{aligned}
    $$
    for any $\theta\in \left(\frac{n\sigma_H}{n\sigma_H+4},1\right)$.
    Finally, by Lemma \ref{lem : Gap inequality with infty} and \eqref{eq : bound of u-(u)}, we obtain
    $$
    \begin{aligned}
        &\omega(\rho) \miint{Q_{2\rho}^\lambda(z_0)} H\left(\frac{|u-u_{Q_{2\rho}^\lambda(z_0)}|}{2\rho}\right)\, dz\\
        &\qquad \leq c \miint{Q_{2\rho}^\lambda(z_0)} G\left(\frac{|u-u_{Q_{2\rho}^\lambda(z_0)}|}{2\rho}\right)\, dz\\
        &\qquad\qquad + c \miint{Q_{2\rho}^\lambda(z_0)} \rho \left(\frac{|u-u_{Q_{2\rho}^\lambda(z_0)}|}{2\rho}\right)G\left(\frac{|u-u_{Q_{2\rho}^\lambda(z_0)}|}{2\rho}\right)\, dz\\
        &\qquad\leq c\miint{Q_{2\rho}^\lambda(z_0)} G\left(\frac{|u-u_{Q_{2\rho}^\lambda(z_0)}|}{2\rho}\right)\, dz.
    \end{aligned}
    $$
    Moreover, Lemma \ref{lem : estimate S in p case with infty}, the $\Delta_2$-condition of $H$ and \eqref{cond : properties in p intrinsic case}$_1$ give
    $$
    \inf_{Q_{2\rho}^\lambda(z_0)}a(\cdot)^{1-\theta}H^{1-\theta}\left(S(u,Q_{2\rho}^\lambda(z_0))^\frac{1}{2}\right)\leq c\left(\sup_{Q_{10\rho}(z_0)}a(\cdot)H(\lambda)\right)^{1-\theta}\leq cG^{1-\theta}(\lambda).
    $$
    Using \eqref{eq : (4.23)} and combining the above inequalities, we obtain the conclusion.
\end{proof}
Proceeding similarly to the proofs of Lemmas \ref{lem : estimate of Psi in G phase with s} and \ref{lem : estimate of psi in level set with s}, with $\theta_0\coloneq\max\{\theta_0'',\frac{1}{i_H},1-\frac{1}{\sigma_H}\}$, where $\theta_0''$ denotes the constant $\theta_0$ of Lemma \ref{lem : estimate of Psi(|u-(u)_Q|) in G intrinsic case with infty}, yields the following two lemmas.
\begin{lemma}\label{lem : estimate of Psi in G phase with infty}
    Let $u$ be a weak solution to \eqref{eq : main equation}. Then there exist constants $c=c(\operatorname{data})$ and $\theta_0=\theta_0(n,i_H,\sigma_H)$ such that if $Q_{2\kappa\rho}^\lambda(z_0)\subset \Omega_T$ satisfies \eqref{cond : properties in p intrinsic case},
    $$
    \miint{Q_{\rho}^\lambda(z_0)} \Psi(z,|Du|)\, dz \leq c\left(\miint{Q_{2\rho}^\lambda(z_0)} \Psi^\theta(z,|Du|)\, dz\right)^\frac{1}{\theta}
    $$
    for any $\theta\in (\theta_0,1)$.
\end{lemma}
\begin{lemma}\label{lem : estimate of psi in level set with infty}
    Let $u$ be a weak solution to \eqref{eq : main equation}. If $Q_{2\kappa \rho}^\lambda(z_0)\subset\Omega_T$ satisfies \eqref{cond : properties in p intrinsic case}, then there exist constants $c=c(\operatorname{data})$ and $\theta_0=\theta_0(n,i_H,\sigma_H)\in (0,1)$ such that for any $\theta\in (\theta_0,1)$,
    $$
    \iints{Q_{2\kappa\rho}^\lambda(z_0)}\Psi(z,|Du|)\, dz \leq c\Lambda^{1-\theta}\iints{Q_{2\rho}^\lambda(z_0)\cap \Upsilon(c^{-1}\Lambda)} \Psi^\theta(z,|Du|)\, dz.
    $$
\end{lemma}

\subsection{The \texorpdfstring{$(G,H)$}{(G,H)}-phase case}\label{subsec : GH phase}
We now assume \ref{case : p,q-intrinsic}, so that $z_0$, $\rho$ and $\lambda$ satisfy
\eqref{cond : properties in p,q intrinsic case} and $a(z_0)>0$ by the discussion in Section \ref{section 3}. The relevant
cylinders are the $\Psi_{z_0}$-intrinsic ones introduced in
\eqref{eq : definition of p,q-intrinsic cylinder}. The point of this case is that the
modulating coefficient is comparable to the constant $a(z_0)$, so that the two phases may be
treated as a single Young function $\Psi_{z_0}$. The estimates then follow from exactly the
same mechanism as in the $G$-phase, with $\Psi_{z_0}$ in place of $G$ and with only one
$\Xi$-function instead of two.

\begin{lemma}\label{lem : Psi z_0 comparability}
    Assume \eqref{cond : properties in p,q intrinsic case}. Then
    \begin{equation}\label{eq : Psi z_0 comparability}
        \tfrac12\,\Psi_{z_0}(\varkappa)\leq \Psi(z,\varkappa)\leq 2\,\Psi_{z_0}(\varkappa)
        \qquad\text{for all }z\in \widetilde{Q}_{4\rho}^\lambda(z_0)\text{ and }\varkappa\geq0 .
    \end{equation}
    Moreover, $\Psi_{z_0}$ is an admissible Young function in the sense of
    \eqref{cond : admissible Phi}, and it satisfies $\Psi_{z_0}\in\Delta_2\cap\nabla_2$ with
    $\Delta_2(\Psi_{z_0})\leq\max\{\Delta_2(G),\Delta_2(H)\}$ and $\nabla_2(\Psi_{z_0})=2$.
    Finally, for $\theta\in(\frac{1}{i_H},1]$, the function
    $\varkappa\mapsto \Psi_{z_0}^\theta(\varkappa)/\varkappa$ is almost increasing on
    $(0,\infty)$ with constant $2\tilde{L}$.
\end{lemma}
\begin{proof}
    By \eqref{cond : properties in p,q intrinsic case}$_1$ we have
    $\frac{a(z_0)}{2}\leq a(z)\leq 2a(z_0)$ for $z\in \widetilde{Q}_{4\rho}^\lambda(z_0)$, whence
    $$
    \tfrac12\Psi_{z_0}(\varkappa)\leq G(\varkappa)+\tfrac{a(z_0)}{2}H(\varkappa)
    \leq \Psi(z,\varkappa)\leq G(\varkappa)+2a(z_0)H(\varkappa)\leq 2\Psi_{z_0}(\varkappa),
    $$
    which is \eqref{eq : Psi z_0 comparability}.

    Since $G$ and $H$ are admissible and $a(z_0)>0$, both
    $\frac{\Psi_{z_0}(\varkappa)}{\varkappa^{2}}
    =\frac{G(\varkappa)}{\varkappa^{2}}+a(z_0)\frac{H(\varkappa)}{\varkappa^{2}}$ and
    $\frac{\Psi_{z_0}(\varkappa)}{\varkappa^{\sigma_H}}
    =\frac{G(\varkappa)}{\varkappa^{\sigma_H}}+a(z_0)\frac{H(\varkappa)}{\varkappa^{\sigma_H}}$
    inherit the monotonicity of their summands, so $\Psi_{z_0}$ is admissible. The bound
    $$\Psi_{z_0}(2\varkappa)\leq \Delta_2(G)G(\varkappa)+\Delta_2(H)a(z_0)H(\varkappa)
    \leq\max\{\Delta_2(G),\Delta_2(H)\}\Psi_{z_0}(\varkappa)$$ gives the $\Delta_2$-condition,
    while the first condition in \eqref{cond : admissible Phi} yields
    $\Psi_{z_0}(2\varkappa)\geq 4\Psi_{z_0}(\varkappa)=2\cdot 2\,\Psi_{z_0}(\varkappa)$, that is,
    $\Psi_{z_0}\in\nabla_2$ with $\nabla_2(\Psi_{z_0})=2$.

    Finally, let $0<\varkappa_1<\varkappa_2$. Using $(x+y)^\theta\leq x^\theta+y^\theta$ for
    $\theta\in(0,1]$, then \eqref{cond : ainc_1 of G and H}, and finally
    $G^\theta\leq\Psi_{z_0}^\theta$ and $a(z_0)^\theta H^\theta=(a(z_0)H)^\theta\leq
    \Psi_{z_0}^\theta$, we get
    \[
    \frac{\Psi_{z_0}^\theta(\varkappa_1)}{\varkappa_1}
    \leq \frac{G^\theta(\varkappa_1)}{\varkappa_1}
    +a(z_0)^\theta\frac{H^\theta(\varkappa_1)}{\varkappa_1}
    \leq \tilde{L}\left[\frac{G^\theta(\varkappa_2)}{\varkappa_2}
    +a(z_0)^\theta\frac{H^\theta(\varkappa_2)}{\varkappa_2}\right]
    \leq 2\tilde{L}\,\frac{\Psi_{z_0}^\theta(\varkappa_2)}{\varkappa_2}. \qedhere
    \]
\end{proof}

\begin{lemma}\label{lem : Poincare for Psi z_0}
    Lemma \ref{lem : semi-Parabolic Poincare inequality} remains valid with $\Phi$ replaced by
    $\Psi_{z_0}$, with a constant
    $c=c(n,N,L,\Delta_2(G),\nabla_2(G),\Delta_2(H),\nabla_2(H))>1$.
\end{lemma}
\begin{proof}
    The proof of Lemma \ref{lem : semi-Parabolic Poincare inequality} uses only the following
    properties of $\Phi$: it is a Young function belonging to $\Delta_2\cap\nabla_2$, and
    $\Phi^\theta(\varkappa)/\varkappa$ is almost increasing, the latter being needed for the
    Sobolev--Poincar\'{e} and Jensen inequalities of \cite{Harjulehto2019}. The remaining
    ingredient, Lemma \ref{lem : gluing lemma}, does not involve $\Phi$.
    All of these properties are supplied by Lemma \ref{lem : Psi z_0 comparability}, with constants
    depending only on the quantities listed, so the proof applies verbatim.
\end{proof}

\begin{lemma}\label{lem : GH phase Poincare}
    Let $u$ be a weak solution to \eqref{eq : main equation} and assume that
    $\widetilde{Q}_{4\rho}^\lambda(z_0)\subset\Omega_T$ satisfies \eqref{cond : properties in p,q intrinsic case}. Then,
    for $\sigma\in[2\rho,4\rho]$ and $\theta$ as in
    \eqref{eq : choice of theta in G phase}, there exists a constant
    $c=c(\operatorname{data})>1$ such that
    $$
    \miint{\widetilde{Q}_\sigma^\lambda(z_0)}\Psi_{z_0}^{\theta}
    \left(\frac{\bigl|u-u_{\widetilde{Q}_\sigma^\lambda(z_0)}\bigr|}{\sigma}\right) dz
    \leq c\miint{\widetilde{Q}_\sigma^\lambda(z_0)} \Psi_{z_0}^\theta(|Du|)\, dz .
    $$
\end{lemma}
\begin{proof}
    We abbreviate $\mathcal{Q}\coloneq \widetilde{Q}_\sigma^\lambda(z_0)$ and put
    $$
    Y\coloneq \miint{\mathcal{Q}} \Psi_{z_0}^\theta(|Du|)\, dz ,
    $$
    which we may assume to be positive, since otherwise $Du=0$ almost everywhere in $\mathcal{Q}$ and
    both sides vanish.

    We first derive the intrinsic bound. Since $\kappa\geq10$ we have
    $2\kappa\rho\geq 20\rho>\sigma$, so that
    \eqref{cond : properties in p,q intrinsic case}$_3$ applies and gives
    $$\miint{\mathcal{Q}}\Psi(z,|Du|)\, dz\leq \Psi_{z_0}(\lambda).$$ Combining this with
    \eqref{eq : Psi z_0 comparability}, which is available because
    $\mathcal{Q}\subset \widetilde{Q}_{4\rho}^\lambda(z_0)$, and with Jensen's inequality, we obtain
    \begin{equation}\label{eq : Y leq Psi lambda}
        Y\leq\left(\miint{\mathcal{Q}}\Psi_{z_0}(|Du|)\, dz\right)^\theta
        \leq \left(2\miint{\mathcal{Q}}\Psi(z,|Du|)\, dz\right)^\theta
        \leq 2\,\Psi_{z_0}^\theta(\lambda).
    \end{equation}
    We set $\epsilon\coloneq Y/\bigl(2\Psi_{z_0}^\theta(\lambda)\bigr)\in(0,1]$.

    Next, we estimate the gluing term. By Lemma \ref{lem : Phi and Phi'},
    \eqref{eq : Psi z_0 comparability} and \eqref{eq : Xi pointwise} applied to the admissible
    function $\Psi_{z_0}$,
    $$
    G'(|Du|)+a(z)H'(|Du|)\leq c\,\frac{\Psi(z,|Du|)}{|Du|}
    \leq c\,\frac{\Psi_{z_0}(|Du|)}{|Du|}
    =c\,\widehat{\Xi}_{\Psi_{z_0}}\!\bigl(\Psi_{z_0}^\theta(|Du|)\bigr)
    $$
    with $c=c(\Delta_2(G),\Delta_2(H))$. By Remark
    \ref{rem : Xi quasiconcave consequences} and Lemma \ref{lem : quasiconcave Jensen},
    followed by \eqref{eq : Xi scaling}, \eqref{eq : Y leq Psi lambda}, the second assertion of
    Remark \ref{rem : Xi quasiconcave consequences} with $\mu=2$, and finally
    \eqref{eq : Xi pointwise},
    $$
    \begin{aligned}
    \miint{\mathcal{Q}}\left[G'(|Du|)+a(z)H'(|Du|)\right] dz
    &\leq c\,\widehat{\Xi}_{\Psi_{z_0}}(Y)
    =c\,\widehat{\Xi}_{\Psi_{z_0}}\!\bigl(\epsilon\cdot 2\Psi_{z_0}^\theta(\lambda)\bigr)\\
    &\leq c\,\epsilon^{\frac{1}{2\theta}}\,
    \widehat{\Xi}_{\Psi_{z_0}}\!\bigl(2\Psi_{z_0}^\theta(\lambda)\bigr)
    \leq c\,\epsilon^{\frac{1}{2\theta}}\,
    \widehat{\Xi}_{\Psi_{z_0}}\!\bigl(\Psi_{z_0}^\theta(\lambda)\bigr)\\
    &=c\,\epsilon^{\frac{1}{2\theta}}\,\frac{\Psi_{z_0}(\lambda)}{\lambda}.
    \end{aligned}
    $$
    Consequently,
    \begin{equation}\label{eq : GH phase gluing term}
        \frac{\lambda^2}{\Psi_{z_0}(\lambda)}
        \miint{\mathcal{Q}}\left[G'(|Du|)+a(z)H'(|Du|)\right] dz
        \leq c\,\epsilon^{\frac{1}{2\theta}}\,\lambda .
    \end{equation}
    We stress that the only parts of \eqref{cond : properties in p,q intrinsic case} used so
    far are the comparability \eqref{cond : properties in p,q intrinsic case}$_1$ and the
    intrinsic condition \eqref{cond : properties in p,q intrinsic case}$_3$.

    Finally, we apply Lemma \ref{lem : Poincare for Psi z_0} with
    $R=\sigma$ and $\ell=\frac{\lambda^2}{\Psi_{z_0}(\lambda)}\sigma^2$, so that
    $Q_{R,\ell}(z_0)=\mathcal{Q}$ and $\frac{\ell}{R^2}=\frac{\lambda^2}{\Psi_{z_0}(\lambda)}$.
    This gives
    $$
    \begin{aligned}
    \miint{\mathcal{Q}}\Psi_{z_0}^{\theta}
    \left(\frac{\bigl|u-u_{\mathcal{Q}}\bigr|}{\sigma}\right) dz
    &\leq c\miint{\mathcal{Q}} \Psi_{z_0}^\theta(|Du|)\, dz\\
    &\qquad+c\,\Psi_{z_0}^\theta\!\left(\frac{\lambda^2}{\Psi_{z_0}(\lambda)}
    \miint{\mathcal{Q}}\left[G'(|Du|)+a(z)H'(|Du|)\right] dz\right).
    \end{aligned}
    $$
    The first term equals $cY$. For the second one we insert
    \eqref{eq : GH phase gluing term}, absorb the constant $c\geq1$ by
    \eqref{eq : Phi(Lambda t) vs Phi(t)} applied to $\Psi_{z_0}\in\Delta_2$, and then use the
    first condition in \eqref{cond : admissible Phi} in the form
    $\Psi_{z_0}(c\varkappa)\leq c^2\Psi_{z_0}(\varkappa)$ for $c\leq1$:
    $$
    \Psi_{z_0}^\theta\!\left(c\,\epsilon^{\frac{1}{2\theta}}\lambda\right)
    \leq c\,\Psi_{z_0}^\theta\!\left(\epsilon^{\frac{1}{2\theta}}\lambda\right)
    \leq c\left(\epsilon^{\frac{1}{2\theta}}\right)^{2\theta}\Psi_{z_0}^\theta(\lambda)
    =c\,\epsilon\,\Psi_{z_0}^\theta(\lambda)=\tfrac{c}{2}\,Y .
    $$
    Combining the last three displays completes the proof.
\end{proof}

\begin{lemma}\label{lem : estimate of S in (G,H) phase}
    Let $u$ be a weak solution to \eqref{eq : main equation}. If $\widetilde{Q}_{2\kappa \rho}^\lambda (z_0)\subset \Omega_T$ satisfies \eqref{cond : properties in p,q intrinsic case}, there exists a constant $c=c(\operatorname{data})>1$ such that
    $$
    S(u,\widetilde{Q}_{2\rho}^\lambda(z_0))=\sup_{J_{2\rho}^\lambda(t_0)}\dashint_{B_{2\rho}(x_0)}\frac{|u-u_{\widetilde{Q}_{2\rho}^\lambda(z_0)}|^2}{(2\rho)^2}\, dx\leq c \lambda^2.
    $$
\end{lemma}
\begin{proof}
    We apply Lemma \ref{lem : Caccioppoli inequality} with $R=4\rho$, $r=2\rho$, $\ell=\frac{\lambda^2}{\Psi_{z_0}(\lambda)}(4\rho)^2$ and $\tau=\frac{\lambda^2}{\Psi_{z_0}(\lambda)}(2\rho)^2$, so that $Q_{R,\ell}(z_0)=\widetilde{Q}_{4\rho}^\lambda(z_0)$ and $Q_{r,\tau}(z_0)=\widetilde{Q}_{2\rho}^\lambda(z_0)$. Arguing as in the proof of Lemma \ref{lem : estimate S in p case with s} and using \eqref{eq : Psi z_0 comparability}, we obtain
    $$
    \begin{aligned}
    \frac{\Psi_{z_0}(\lambda)}{\lambda^2}S(u,\widetilde{Q}_{2\rho}^\lambda (z_0))&\leq c\miint{\widetilde{Q}_{4\rho}^\lambda (z_0)} \Psi_{z_0}\left(\frac{|u-u_{\widetilde{Q}_{4\rho}^\lambda (z_0)}|}{4\rho}\right)\, dz\\ &\qquad+ \frac{c\Psi_{z_0}(\lambda)}{\lambda^2}\miint{\widetilde{Q}_{4\rho}^\lambda (z_0)} \frac{|u-u_{\widetilde{Q}_{4\rho}^\lambda (z_0)}|^2}{(4\rho)^2}\, dz
    \end{aligned}
    $$
    for some $c=c(n,\nu,L,\Delta_2(G),\nabla_2(G),\Delta_2(H),\nabla_2(H))>1$.
    By Lemma \ref{lem : GH phase Poincare}, \eqref{eq : Psi z_0 comparability} and \eqref{cond : properties in p,q intrinsic case}$_3$, we have
    $$
    \begin{aligned}
        \miint{\widetilde{Q}_{4\rho}^\lambda (z_0)} \Psi_{z_0}\left(\frac{|u-u_{\widetilde{Q}_{4\rho}^\lambda (z_0)}|}{4\rho}\right)\, dz&\leq c\miint{\widetilde{Q}_{4\rho}^\lambda (z_0)} \Psi_{z_0}\left(|Du|\right)\, dz\leq c\Psi_{z_0}(\lambda).
    \end{aligned}
    $$
    By Lemma \ref{lem : Young function of Psi_z_0}, we can use Jensen's inequality \cite[Lemma 4.3.1]{Harjulehto2019} and get
    $$
    \begin{aligned}
    \Psi_{z_0}\left(\left(\miint{\widetilde{Q}_{4\rho}^\lambda (z_0)} \frac{|u-u_{\widetilde{Q}_{4\rho}^\lambda (z_0)}|^2}{(4\rho)^2}\, dz\right)^\frac{1}{2}\right)&\leq c\miint{\widetilde{Q}_{4\rho}^\lambda (z_0)} \Psi_{z_0}\left(\frac{|u-u_{\widetilde{Q}_{4\rho}^\lambda (z_0)}|}{4\rho}\right)\, dz\\
    &\leq c\Psi_{z_0}(\lambda).
    \end{aligned}
    $$
    Hence, we have
    $$
    \miint{\widetilde{Q}_{4\rho}^\lambda (z_0)} \frac{|u-u_{\widetilde{Q}_{4\rho}^\lambda (z_0)}|^2}{(4\rho)^2}\, dz\leq c\lambda^2.
    $$
    Combining the above inequalities, we obtain $S(u,\widetilde{Q}_{2\rho}^\lambda (z_0))\leq c\lambda^2$.
\end{proof}
\begin{lemma}\label{lem : estimate of Psi u in (G,H) phase}
    Let $u$ be a weak solution to \eqref{eq : main equation}. If $\widetilde{Q}_{2\kappa \rho}^\lambda(z_0)\subset\Omega_T$ satisfies \eqref{cond : properties in p,q intrinsic case}, then there exists a constant $c=c(\operatorname{data})>1$ such that for any $\theta\in \left(\frac{n\sigma_H}{n\sigma_H+4},1\right)$,
    $$
    \begin{aligned}
        &\miint{\widetilde{Q}_{2\rho}^\lambda (z_0)}\Psi\left(z,\frac{|u-u_{\widetilde{Q}_{2\rho}^\lambda(z_0)}|}{2\rho}\right)\, dz\\
        &\qquad \leq c\Psi_{z_0}^{1-\theta}(\lambda)\miint{\widetilde{Q}_{2\rho}^\lambda (z_0)}\left[\Psi_{z_0}^\theta\left(\frac{|u-u_{\widetilde{Q}_{2\rho}^\lambda(z_0)}|}{2\rho}\right)+\Psi_{z_0}^\theta\left(|Du|\right)\right]\, dz.
    \end{aligned}
    $$
\end{lemma}
\begin{proof}
    By \eqref{eq : Psi z_0 comparability}, we obtain
    $$
    \begin{aligned}
        \miint{\widetilde{Q}_{2\rho}^\lambda (z_0)}\Psi\left(z,\frac{|u-u_{\widetilde{Q}_{2\rho}^\lambda(z_0)}|}{2\rho}\right)\, dz\leq 2\miint{\widetilde{Q}_{2\rho}^\lambda (z_0)}\Psi_{z_0}\left(\frac{|u-u_{\widetilde{Q}_{2\rho}^\lambda(z_0)}|}{2\rho}\right)\, dz.
    \end{aligned}
    $$
    Using Lemma \ref{lem : Gagliardo-Nirenberg inequality} with $\Phi=\Psi_{z_0}^\frac{1}{2}$, $q_1=\frac{\sigma_H}{2}$, $\gamma=2$, $q_2=2\theta$, $q_3=2$ and $\theta\in \left(\frac{n\sigma_H}{n\sigma_H+4},1\right)$ gives
    $$
    \begin{aligned}
        &\miint{\widetilde{Q}_{2\rho}^\lambda (z_0)}\Psi_{z_0}\left(\frac{|u-u_{\widetilde{Q}_{2\rho}^\lambda(z_0)}|}{2\rho}\right)\, dz\\
        &\quad\leq c\left(\miint{\widetilde{Q}_{2\rho}^\lambda (z_0)}\left[\Psi_{z_0}^\theta\left(\frac{|u-u_{\widetilde{Q}_{2\rho}^\lambda(z_0)}|}{2\rho}\right)+\Psi_{z_0}^\theta(|Du|)\right]\, dz\right)\Psi_{z_0}^{1-\theta}(S(u,\widetilde{Q}_{2\rho}^\lambda(z_0))^\frac{1}{2}).
    \end{aligned}
    $$
    By Lemma \ref{lem : estimate of S in (G,H) phase}, we have
    $$
    \begin{aligned}
        &\miint{\widetilde{Q}_{2\rho}^\lambda (z_0)}\Psi_{z_0}\left(\frac{|u-u_{\widetilde{Q}_{2\rho}^\lambda(z_0)}|}{2\rho}\right)\, dz\\
        &\quad\leq c\Psi_{z_0}^{1-\theta}(\lambda)\miint{\widetilde{Q}_{2\rho}^\lambda (z_0)}\left[\Psi_{z_0}^\theta\left(\frac{|u-u_{\widetilde{Q}_{2\rho}^\lambda(z_0)}|}{2\rho}\right)+\Psi_{z_0}^\theta(|Du|)\right]\, dz.
    \end{aligned}
    $$
\end{proof}
\begin{lemma}\label{lem : estimate in level set in (G,H) phase}
    Let $u$ be a weak solution to \eqref{eq : main equation}. If $\widetilde{Q}_{2\kappa \rho}^\lambda(z_0)\subset \Omega_T$ satisfies \eqref{cond : properties in p,q intrinsic case}, then there exist constants $c=c(\operatorname{data})$ and $\theta_0=\theta_0(n,i_H,\sigma_H)\in (0,1)$ such that for any $\theta\in (\theta_0,1)$,
    $$
    \begin{aligned}
        \miint{\widetilde{Q}_{\rho}^\lambda(z_0)} \Psi_{z_0}(|Du|)\, dz \leq c\left(\miint{\widetilde{Q}_{2\rho}^\lambda(z_0)} \Psi_{z_0}^\theta (|Du|)\, dz\right)^\frac{1}{\theta}
    \end{aligned}
    $$
    and
    $$
    \iints{\widetilde{Q}_{2\kappa \rho}^\lambda(z_0)} \Psi(z,|Du|)\, dz\leq c\Lambda^{1-\theta}\iints{\widetilde{Q}_{2\rho}^\lambda(z_0)\cap \Upsilon(c^{-1}\Lambda)} \Psi^\theta(z,|Du|)\, dz.
    $$
\end{lemma}
\begin{proof}
    Once the first inequality is established, the second one follows by the same argument as in the proof of Lemma \ref{lem : estimate of psi in level set with s}, with \eqref{eq : Psi z_0 comparability} and \eqref{cond : properties in p,q intrinsic case} in place of \eqref{cond : properties in p intrinsic case}.

    Now, we prove the first inequality in the statement. Throughout the proof we take
    $$
    \theta_0\coloneq \max\left\{\frac{n\sigma_H}{n\sigma_H+4},\,\frac{1}{i_H},\,1-\frac{1}{\sigma_H}\right\}\in(0,1),
    $$
    so that every $\theta\in(\theta_0,1)$ is admissible in Lemma \ref{lem : estimate of Psi u in (G,H) phase} and satisfies \eqref{eq : choice of theta in G phase}. By Lemma \ref{lem : Caccioppoli inequality} and \eqref{eq : Psi z_0 comparability},
    $$
    \begin{aligned}
        \miint{\widetilde{Q}_{\rho}^\lambda(z_0)} \Psi(z,|Du|)\, dz &\leq c\miint{\widetilde{Q}_{2\rho}^\lambda(z_0)} \Psi_{z_0}\left(\frac{|u-u_{\widetilde{Q}_{2\rho}^\lambda(z_0)}|}{2\rho}\right)\, dz\\
        &\qquad+c\frac{\Psi_{z_0}(\lambda)}{\lambda^2}\miint{\widetilde{Q}_{2\rho}^\lambda(z_0)} \frac{|u-u_{\widetilde{Q}_{2\rho}^\lambda(z_0)}|^2}{(2\rho)^2}\, dz.
    \end{aligned}
    $$
    Applying Lemmas \ref{lem : estimate of Psi u in (G,H) phase} and \ref{lem : GH phase Poincare} together with \eqref{eq : Psi z_0 comparability}, we have
    $$
    \begin{aligned}
        \miint{\widetilde{Q}_{2\rho}^\lambda(z_0)} \Psi_{z_0}\left(\frac{|u-u_{\widetilde{Q}_{2\rho}^\lambda(z_0)}|}{2\rho}\right)\, dz&\leq c\Psi_{z_0}^{1-\theta}(\lambda)\miint{\widetilde{Q}_{2\rho}^\lambda(z_0)} \Psi_{z_0}^\theta\left(|Du|\right)\, dz.
    \end{aligned}
    $$
    By Lemma \ref{lem : estimate of S in (G,H) phase}, we have 
    $$
    \begin{aligned}
    \miint{\widetilde{Q}_{2\rho}^\lambda(z_0)} \frac{|u-u_{\widetilde{Q}_{2\rho}^\lambda(z_0)}|^2}{(2\rho)^2}\, dz &\leq \left(S(u, \widetilde{Q}_{2\rho}^\lambda(z_0))\right)^{\frac{1}{2}}\left(\miint{\widetilde{Q}_{2\rho}^\lambda(z_0)} \frac{|u-u_{\widetilde{Q}_{2\rho}^\lambda(z_0)}|^2}{(2\rho)^2}\, dz\right)^\frac{1}{2}\\
    &\leq c\lambda\left(\miint{\widetilde{Q}_{2\rho}^\lambda(z_0)} \frac{|u-u_{\widetilde{Q}_{2\rho}^\lambda(z_0)}|^2}{(2\rho)^2}\, dz\right)^\frac{1}{2}.
    \end{aligned}
    $$
    Then it follows from \eqref{eq : comparison with Phi and Phi'}, \eqref{Youngs-ineq} with $\Phi=\Psi_{z_0}^*$ and \eqref{eq : Phi*(Phi')} that
    $$
    \begin{aligned}
        &\frac{\Psi_{z_0}(\lambda)}{\lambda^2} \miint{\widetilde{Q}_{2\rho}^\lambda(z_0)} \frac{|u-u_{\widetilde{Q}_{2\rho}^\lambda(z_0)}|^2}{(2\rho)^2}\, dz\leq c\frac{\Psi_{z_0}(\lambda)}{\lambda}\left(\miint{\widetilde{Q}_{2\rho}^\lambda(z_0)} \frac{|u-u_{\widetilde{Q}_{2\rho}^\lambda(z_0)}|^2}{(2\rho)^2}\, dz\right)^\frac{1}{2}\\
        &\qquad\qquad\qquad\qquad\leq \frac{1}{4}\Psi_{z_0}(\lambda)+c\Psi_{z_0}\left(\left(\miint{\widetilde{Q}_{2\rho}^\lambda(z_0)} \frac{|u-u_{\widetilde{Q}_{2\rho}^\lambda(z_0)}|^2}{(2\rho)^2}\, dz\right)^\frac{1}{2}\right).
    \end{aligned}
    $$
    By Lemma \ref{lem : Young function of Psi_z_0}, Jensen's inequality \cite[Lemma 4.3.1]{Harjulehto2019} and the above result, we have
    $$
    \begin{aligned}
        \Psi_{z_0}\left(\left(\miint{\widetilde{Q}_{2\rho}^\lambda(z_0)} \frac{|u-u_{\widetilde{Q}_{2\rho}^\lambda(z_0)}|^2}{(2\rho)^2}\, dz\right)^\frac{1}{2}\right)&\leq c\miint{\widetilde{Q}_{2\rho}^\lambda(z_0)} \Psi_{z_0}\left(\frac{|u-u_{\widetilde{Q}_{2\rho}^\lambda(z_0)}|}{2\rho}\right)\, dz\\
        &\leq c \Psi_{z_0}^{1-\theta}(\lambda)\miint{\widetilde{Q}_{2\rho}^\lambda(z_0)} \Psi_{z_0}^\theta(|Du|)\, dz.
    \end{aligned}
    $$
    Thus, combining the above inequalities and using Young's inequality, we have
    $$
    \begin{aligned}
    \miint{\widetilde{Q}_{\rho}^\lambda(z_0)} \Psi(z,|Du|)\, dz &\leq \frac{1}{4}\Psi_{z_0}(\lambda) + c \Psi_{z_0}^{1-\theta}(\lambda)\miint{\widetilde{Q}_{2\rho}^\lambda(z_0)} \Psi_{z_0}^\theta(|Du|)\, dz\\
    &\leq \frac{1}{2}\Psi_{z_0}(\lambda) + c\left(\miint{\widetilde{Q}_{2\rho}^\lambda(z_0)} \Psi_{z_0}^\theta(|Du|)\, dz\right)^\frac{1}{\theta}.
    \end{aligned}
    $$
    Since the left-hand side equals $\Psi_{z_0}(\lambda)$ by \eqref{cond : properties in p,q intrinsic case}$_2$, we obtain
    $$
    \Psi_{z_0}(\lambda)\leq 2c\left(\miint{\widetilde{Q}_{2\rho}^\lambda(z_0)} \Psi_{z_0}^\theta(|Du|)\, dz\right)^\frac{1}{\theta},
    $$
    and hence, by \eqref{eq : Psi z_0 comparability} and \eqref{cond : properties in p,q intrinsic case}$_2$ again,
    $$
    \begin{aligned}
    \miint{\widetilde{Q}_{\rho}^\lambda(z_0)} \Psi_{z_0}(|Du|)\, dz&\leq 2\miint{\widetilde{Q}_{\rho}^\lambda(z_0)} \Psi(z,|Du|)\, dz=2\Psi_{z_0}(\lambda)\\
    &\leq c\left(\miint{\widetilde{Q}_{2\rho}^\lambda(z_0)} \Psi_{z_0}^\theta(|Du|)\, dz\right)^\frac{1}{\theta},
    \end{aligned}
    $$
    which is the first inequality in the statement.
\end{proof}
\section{Proof of the main theorems}\label{section 5}
\subsection{Vitali covering argument}\label{subsec : Vitali Covering lemma}
For each $w\in \Upsilon(\Lambda,r_1)$ that is a Lebesgue point as in Section \ref{section 3}, we consider
$$
\mathcal{Q}(w)= \begin{cases}
    Q_{2\rho_w}^{\lambda_w}(w)\quad &\text{if \ref{case : p-intrinsic} holds},\\
    \widetilde{Q}_{2\varsigma_w}^{\lambda_w}(w)\quad &\text{otherwise},
\end{cases} 
$$
where $\lambda_w$, $\rho_w$ and $\varsigma_w$ are defined in Section \ref{section 3}. According to \cite[Subsection 5.2]{2023_Gradient_Higher_Integrability_for_Degenerate_Parabolic_Double-Phase_Systems}, one can find a countable family $\mathcal{G}$ of pairwise disjoint cylinders within $\mathcal{F} \coloneq \{\mathcal{Q}(w) : w \in \Upsilon(\Lambda, r_1)\}$ such that:
\begin{itemize}
    \item For each $\mathcal{Q}(z_1) \in \mathcal{F}$, there is a cylinder $\mathcal{Q}(z_2) \in \mathcal{G}$ for which
    \begin{equation}\label{eq : intersect of Q(z_1) and Q(z_2)}
    \mathcal{Q}(z_1) \cap \mathcal{Q}(z_2) \neq \emptyset.
    \end{equation}
    \item The points $z_1$ and $z_2$ fulfill
    \begin{equation}\label{eq : relation of ells}
    \ell_{z_1} \le 2\ell_{z_2},
    \end{equation}
    where $\ell_{z_i}\coloneq 2\rho_{z_i}$ if \ref{case : p-intrinsic} holds at $z_i$ and $\ell_{z_i}\coloneq 2\varsigma_{z_i}$ otherwise.
\end{itemize}
Then, we assert that for the above points $z_1$ and $z_2$, the inclusion
$$
\mathcal{Q}(z_1)\subset \kappa \mathcal{Q}(z_2)
$$
holds, where for $\mathcal{Q}=Q^\lambda_\ell(z)$ or $\mathcal{Q}=\widetilde{Q}^\lambda_\ell(z)$ we write $\kappa\mathcal{Q}\coloneq Q^\lambda_{\kappa\ell}(z)$ or $\kappa\mathcal{Q}\coloneq\widetilde{Q}^\lambda_{\kappa\ell}(z)$, respectively. Thus $\kappa\mathcal{Q}$ has the same center and the same scaling parameter as $\mathcal{Q}$, while its radius is multiplied by $\kappa$ and its time length by $\kappa^2$. We establish this assertion by examining the following four cases:
\begin{enumerate}[label = (\roman*)]
    \item\label{case : G phases} $\mathcal{Q}(z_2)= Q_{\ell_{z_2}}^{\lambda_{z_2}}(z_2)$ and $\mathcal{Q}(z_1)= Q_{\ell_{z_1}}^{\lambda_{z_1}}(z_1)$,
    \item\label{case : (G,H) phases} $\mathcal{Q}(z_2)= \widetilde{Q}_{\ell_{z_2}}^{\lambda_{z_2}}(z_2)$ and $\mathcal{Q}(z_1)= \widetilde{Q}_{\ell_{z_1}}^{\lambda_{z_1}}(z_1)$,
    \item\label{case : G phase and (G,H) phase} $\mathcal{Q}(z_2)= \widetilde{Q}_{\ell_{z_2}}^{\lambda_{z_2}}(z_2)$ and $\mathcal{Q}(z_1)= Q_{\ell_{z_1}}^{\lambda_{z_1}}(z_1)$,
    \item\label{case : (G,H) phase and G phase} $\mathcal{Q}(z_2)= Q_{\ell_{z_2}}^{\lambda_{z_2}}(z_2)$ and $\mathcal{Q}(z_1)= \widetilde{Q}_{\ell_{z_1}}^{\lambda_{z_1}}(z_1)$.
\end{enumerate}
By \cite[Subsection 5.2]{2023_Gradient_Higher_Integrability_for_Degenerate_Parabolic_Double-Phase_Systems}, the spatial part of the assertion follows from \eqref{eq : intersect of Q(z_1) and Q(z_2)} and \eqref{eq : relation of ells} by scaling the radius around $z_2$ by a factor of $5$. To show the assertion for the time variable, a comparison between $\lambda_{z_1}$ and $\lambda_{z_2}$ is required. We first observe that
\begin{equation}\label{eq : z_1 in Q_10 rho of z_2}
z_1\in Q_{10\rho_{z_2}}(z_2).
\end{equation}
Indeed, write $z_i=(x_i,t_i)$ for $i=1,2$. Since $\lambda_{z_i}\geq\lambda_0\geq 1$, Lemma \ref{lem : quadratic ainc} gives $\frac{\lambda_{z_i}^2}{\Psi_{z_i}(\lambda_{z_i})}\leq\frac{\lambda_{z_i}^2}{G(\lambda_{z_i})}\leq 1$, so that the half-length of the time interval of $\mathcal{Q}(z_i)$ is at most $\ell_{z_i}^2$. Hence \eqref{eq : intersect of Q(z_1) and Q(z_2)}, \eqref{eq : relation of ells} and $\ell_{z_2}\leq 2\rho_{z_2}$ imply
$$
|x_1-x_2|<\ell_{z_1}+\ell_{z_2}\leq 3\ell_{z_2}\leq 6\rho_{z_2}
$$
and
$$
|t_1-t_2|<\ell_{z_1}^2+\ell_{z_2}^2\leq 5\ell_{z_2}^2\leq 20\rho_{z_2}^2<(10\rho_{z_2})^2,
$$
which is \eqref{eq : z_1 in Q_10 rho of z_2}. Consequently, if \ref{case : p-intrinsic} holds at $z_2$, then
\begin{equation}\label{eq : a(z_1) H bound from z_2}
a(z_1)H(\lambda_{z_2})\leq \sup_{Q_{10\rho_{z_2}}(z_2)}a(\cdot)H(\lambda_{z_2})\leq KG(\lambda_{z_2}),
\end{equation}
while if \ref{case : p,q-intrinsic} holds at $z_2$, then \eqref{eq : comparison of a in p,q phase} gives
\begin{equation}\label{eq : a(z_1) comparable to a(z_2)}
\frac{a(z_2)}{2}\leq a(z_1)\leq 2a(z_2).
\end{equation}
Moreover, if \ref{case : p-intrinsic} holds at $z_i$, then $a(z_i)H(\lambda_{z_i})\leq KG(\lambda_{z_i})$, and therefore \eqref{eq : definition of lambda_w} yields
\begin{equation}\label{eq : Lambda leq (1+K) G in G phase}
\Lambda=\Psi_{z_i}(\lambda_{z_i})\leq (1+K)G(\lambda_{z_i}).
\end{equation}

We claim that
\begin{equation}\label{eq : lambda_z_2 lesssim lambda_z_1}
    (2\tilde{L}(1+K))^{-\frac{1}{i_H}}\lambda_{z_1}\leq\lambda_{z_2}\leq (2\tilde{L}(1+K))^\frac{1}{i_H}\lambda_{z_1}.
\end{equation}
For the sake of contradiction, suppose that the second inequality in \eqref{eq : lambda_z_2 lesssim lambda_z_1} fails, that is, $\lambda_{z_2}>(2\tilde{L}(1+K))^\frac{1}{i_H}\lambda_{z_1}$. Since $\varkappa\mapsto \frac{G(\varkappa)}{\varkappa^{i_H}}$ is almost increasing with constant $\tilde{L}$, we get
$$
2(1+K)G(\lambda_{z_1})\leq G\left((2\tilde{L}(1+K))^\frac{1}{i_H}\lambda_{z_1}\right)\leq G(\lambda_{z_2}),
$$
that is,
$$
G(\lambda_{z_1})\leq \frac{1}{2(1+K)}G(\lambda_{z_2}).
$$
Also, since $\lambda_{z_1}<\lambda_{z_2}$, Lemma \ref{lem208} gives
$$
H(\lambda_{z_1})=\frac{H(\lambda_{z_1})}{G(\lambda_{z_1})}G(\lambda_{z_1})\leq\frac{H(\lambda_{z_2})}{G(\lambda_{z_2})}\cdot\frac{1}{2(1+K)}G(\lambda_{z_2})=\frac{1}{2(1+K)}H(\lambda_{z_2}).
$$
We now distinguish three cases. If \ref{case : p-intrinsic} holds at $z_1$, then \eqref{eq : Lambda leq (1+K) G in G phase} and \eqref{eq : definition of lambda_w} give
$$
\Lambda\leq (1+K)G(\lambda_{z_1})\leq \frac{1}{2}G(\lambda_{z_2})\leq \frac{1}{2}\Psi_{z_2}(\lambda_{z_2})=\frac{1}{2}\Lambda.
$$
If \ref{case : p,q-intrinsic} holds at both $z_1$ and $z_2$, then \eqref{eq : a(z_1) comparable to a(z_2)} and \eqref{eq : definition of lambda_w} give
$$
\Lambda=G(\lambda_{z_1})+a(z_1)H(\lambda_{z_1})\leq \frac{G(\lambda_{z_2})+2a(z_2)H(\lambda_{z_2})}{2(1+K)}\leq \frac{\Psi_{z_2}(\lambda_{z_2})}{1+K}\leq\frac{1}{2}\Lambda.
$$
If \ref{case : p,q-intrinsic} holds at $z_1$ and \ref{case : p-intrinsic} holds at $z_2$, then \eqref{eq : a(z_1) H bound from z_2} and \eqref{eq : definition of lambda_w} give
$$
\Lambda=G(\lambda_{z_1})+a(z_1)H(\lambda_{z_1})\leq \frac{G(\lambda_{z_2})+a(z_1)H(\lambda_{z_2})}{2(1+K)}\leq \frac{(1+K)G(\lambda_{z_2})}{2(1+K)}\leq \frac{1}{2}\Lambda.
$$
In each case we arrive at $\Lambda\leq\frac{1}{2}\Lambda$, which contradicts $\Lambda>0$. Hence the second inequality in \eqref{eq : lambda_z_2 lesssim lambda_z_1} holds.

Next, suppose that the first inequality in \eqref{eq : lambda_z_2 lesssim lambda_z_1} fails, that is, $\lambda_{z_1}>(2\tilde{L}(1+K))^\frac{1}{i_H}\lambda_{z_2}$. Exchanging the roles of $z_1$ and $z_2$ in the above argument, we obtain
$$
G(\lambda_{z_2})\leq \frac{1}{2(1+K)}G(\lambda_{z_1})\quad\text{and}\quad H(\lambda_{z_2})\leq\frac{1}{2(1+K)}H(\lambda_{z_1}).
$$
If \ref{case : p-intrinsic} holds at $z_2$, then \eqref{eq : Lambda leq (1+K) G in G phase} and \eqref{eq : definition of lambda_w} give $\Lambda\leq (1+K)G(\lambda_{z_2})\leq\frac{1}{2}G(\lambda_{z_1})\leq\frac{1}{2}\Lambda$. If \ref{case : p,q-intrinsic} holds at both $z_1$ and $z_2$, then \eqref{eq : a(z_1) comparable to a(z_2)} and \eqref{eq : definition of lambda_w} give
$$
\Lambda=G(\lambda_{z_2})+a(z_2)H(\lambda_{z_2})\leq \frac{G(\lambda_{z_1})+2a(z_1)H(\lambda_{z_1})}{2(1+K)}\leq \frac{\Psi_{z_1}(\lambda_{z_1})}{1+K}\leq\frac{1}{2}\Lambda.
$$
If \ref{case : p,q-intrinsic} holds at $z_2$ and \ref{case : p-intrinsic} holds at $z_1$, then \eqref{eq : a(z_1) comparable to a(z_2)}, $a(z_1)H(\lambda_{z_1})\leq KG(\lambda_{z_1})$ and \eqref{eq : definition of lambda_w} give
$$
\begin{aligned}
\Lambda=G(\lambda_{z_2})+a(z_2)H(\lambda_{z_2})&\leq \frac{G(\lambda_{z_1})+2a(z_1)H(\lambda_{z_1})}{2(1+K)}\\
&\leq \frac{(1+2K)G(\lambda_{z_1})}{2(1+K)}\leq\frac{1+2K}{2+2K}\Lambda<\Lambda.
\end{aligned}
$$
In each case we obtain a contradiction, and \eqref{eq : lambda_z_2 lesssim lambda_z_1} is proved. Since $2(1+K)\leq 4K$, \eqref{eq : lambda_z_2 lesssim lambda_z_1} implies
\begin{equation}\label{eq : comparison of lambda_z_1 and lambda_z_2}
    (4\tilde{L}K)^{-\frac{1}{i_H}}\lambda_{z_1}\leq \lambda_{z_2}\leq (4\tilde{L}K)^\frac{1}{i_H}\lambda_{z_1},
\end{equation}
and then \eqref{cond : monotonicity of Phi over powers} together with \eqref{cond : relation of sigma G and sigma H} implies
\begin{equation}\label{eq : comparison of G(lambda_z_1) and G(lambda_z_2)}
    (4\tilde{L}K)^{-\frac{\sigma_H}{i_H}}G(\lambda_{z_1})\leq G(\lambda_{z_2})\leq (4\tilde{L}K)^\frac{\sigma_H}{i_H}G(\lambda_{z_1}).
\end{equation}

Now, we prove the assertion in the time variable. We note that, since $\frac{\sigma_H+2}{i_H}\geq 1$, $\tilde{L}\geq 1$ and $K\geq 1$,
\begin{equation}\label{eq : constants below kappa square}
\begin{aligned}
(32\tilde{L}K)^{\frac{\sigma_H+2}{i_H}}+1&\leq (32\tilde{L}K^2)^{\frac{\sigma_H+2}{i_H}}+1\leq (32\tilde{L}K^2+68\tilde{L}K^2)^{\frac{\sigma_H+2}{i_H}}\\
&\leq (100\tilde{L}^2K^2)^{\frac{\sigma_H+2}{i_H}}=\kappa^2,
\end{aligned}
\end{equation}
where we used $(a+b)^{\gamma}\geq a^{\gamma}+b^{\gamma}$ for $a,\,b\geq0$ and $\gamma\geq1$. Also, $8(4\tilde{L}K)^{\frac{2}{i_H}}+1\leq (32\tilde{L}K)^{\frac{\sigma_H+2}{i_H}}+1$, because $\frac{\sigma_H+2}{i_H}\geq \max\{1,\frac{2}{i_H}\}$. Together with \eqref{eq : constants below kappa square}, this justifies the last inequality in each of the following four cases.

\ref{case : G phases} For any $\tau\in I_{\ell_{z_1}}^{\lambda_{z_1}}(t_1)$, by \eqref{eq : intersect of Q(z_1) and Q(z_2)}, \eqref{eq : relation of ells}, \eqref{eq : comparison of lambda_z_1 and lambda_z_2} and \eqref{eq : comparison of G(lambda_z_1) and G(lambda_z_2)}, we have
$$
\begin{aligned}
    |\tau-t_2|&\leq |\tau - t_1| + |t_1 - t_2|\leq 2 \frac{\lambda_{z_1}^{2}}{G(\lambda_{z_1})}\ell_{z_1}^2+\frac{\lambda_{z_2}^{2}}{G(\lambda_{z_2})}\ell_{z_2}^2\\
    &\leq \left((32\tilde{L}K)^{\frac{\sigma_H+2}{i_H}}+1\right)\frac{\lambda_{z_2}^2}{G(\lambda_{z_2})}\ell_{z_2}^2\leq \frac{\lambda_{z_2}^2}{G(\lambda_{z_2})} (\kappa \ell_{z_2})^2. 
\end{aligned}
$$

\ref{case : (G,H) phases} For any $\tau\in J_{\ell_{z_1}}^{\lambda_{z_1}}(t_1)$, we have
$$
|\tau-t_2|\leq |\tau-t_1|+|t_1-t_2|\leq 2\frac{\lambda_{z_1}^2}{\Psi_{z_1}(\lambda_{z_1})}\ell_{z_1}^2+\frac{\lambda_{z_2}^2}{\Psi_{z_2}(\lambda_{z_2})}\ell_{z_2}^2.
$$
By \eqref{eq : definition of lambda_w} and \eqref{eq : comparison of lambda_z_1 and lambda_z_2}, we have
$$
\frac{\lambda_{z_1}^2}{\Psi_{z_1}(\lambda_{z_1})}=\frac{\lambda_{z_1}^2}{\Lambda}\leq (4\tilde{L}K)^\frac{2}{i_H}\frac{\lambda_{z_2}^2}{\Lambda}\leq (4\tilde{L}K)^{\frac{\sigma_H+2}{i_H}}\frac{\lambda_{z_2}^2}{\Psi_{z_2}(\lambda_{z_2})}.
$$
By \eqref{eq : relation of ells}, we obtain
$$
|\tau-t_2|\leq\left((32\tilde{L}K)^{\frac{\sigma_H+2}{i_H}}+1\right)\frac{\lambda_{z_2}^2}{\Psi_{z_2}(\lambda_{z_2})}\ell_{z_2}^2\leq \frac{\lambda_{z_2}^2}{\Psi_{z_2}(\lambda_{z_2})}(\kappa\ell_{z_2})^2.
$$

\ref{case : G phase and (G,H) phase} For any $\tau\in I_{\ell_{z_1}}^{\lambda_{z_1}}(t_1)$, 
$$
|\tau-t_2|\leq |\tau-t_1|+|t_1-t_2|\leq 2\frac{\lambda_{z_1}^2}{G(\lambda_{z_1})}\ell_{z_1}^2+\frac{\lambda_{z_2}^2}{\Psi_{z_2}(\lambda_{z_2})}\ell_{z_2}^2.
$$
Since $KG(\lambda_{z_1})\geq a(z_1)H(\lambda_{z_1})$, we get
$$
\frac{\lambda_{z_1}^2}{G(\lambda_{z_1})}\leq \frac{2K\lambda_{z_1}^2}{\Psi_{z_1}(\lambda_{z_1})}=\frac{2K\lambda_{z_1}^2}{\Lambda}\leq 2K(4\tilde{L}K)^{\frac{2}{i_H}}\frac{\lambda_{z_2}^2}{\Lambda}\leq (4\tilde{L}K)^{\frac{\sigma_H+2}{i_H}}\frac{\lambda_{z_2}^2}{\Lambda}.
$$
Hence we have
$$
|\tau-t_2|\leq \left((32\tilde{L}K^2)^{\frac{\sigma_H+2}{i_H}}+1\right)\frac{\lambda_{z_2}^2}{\Psi_{z_2}(\lambda_{z_2})}\ell_{z_2}^2\leq \frac{\lambda_{z_2}^2}{\Psi_{z_2}(\lambda_{z_2})}(\kappa\ell_{z_2})^2.
$$

\ref{case : (G,H) phase and G phase} For any $\tau \in J_{\ell_{z_1}}^{\lambda_{z_1}}(t_1)$, we get from \eqref{eq : definition of lambda_w} and \eqref{eq : comparison of lambda_z_1 and lambda_z_2} that
$$
\begin{aligned}
    |\tau-t_2|&\leq |\tau-t_1|+|t_1-t_2|\leq 2\frac{\lambda_{z_1}^2}{\Lambda}\ell_{z_1}^2 + \frac{\lambda_{z_2}^2}{G(\lambda_{z_2})}\ell_{z_2}^2\\
    &\leq \left(8(4\tilde{L}K)^\frac{2}{i_H}+1\right)\frac{\lambda_{z_2}^2}{G(\lambda_{z_2})}\ell_{z_2}^2\leq \frac{\lambda_{z_2}^2}{G(\lambda_{z_2})} (\kappa\ell_{z_2})^2.
\end{aligned}
$$
Thus, we obtain
$$
\mathcal{Q}(z_1)\subset \kappa \mathcal{Q}(z_2).
$$
\subsection{Proof of Theorems \ref{thm : main theorem with infty} and \ref{thm : main theorem with s}}
This proof is analogous to that in \cite[Subsection 5.3]{2023_Gradient_Higher_Integrability_for_Degenerate_Parabolic_Double-Phase_Systems}. First, let $\mathcal{G}$ be the countable pairwise disjoint collection of intrinsic cylinders constructed in Subsection \ref{subsec : Vitali Covering lemma}, whose members are denoted by $\mathcal{Q}_j = \mathcal{Q}(z_j)$ with $z_j\in \Upsilon(\Lambda,r_1)$ and $j\in\mathbb{N}$. Let $\theta_0\in(0,1)$ be the maximum of the constant $\theta_0$ in Lemma \ref{lem : estimate in level set in (G,H) phase} and the constant $\theta_0$ in Lemma \ref{lem : estimate of psi in level set with s} or Lemma \ref{lem : estimate of psi in level set with infty}, according as \eqref{cond : main assumption with s} or \eqref{cond : main assumption with infty} is assumed, so that $\theta_0=\theta_0(n,i_H,\sigma_H,s)$ in the former case and $\theta_0=\theta_0(n,i_H,\sigma_H)$ in the latter, and set $\theta\coloneq\frac{\theta_0+1}{2}$. Since $\kappa\mathcal{Q}_j\subset Q_{r_2}(z_0)\subset\Omega_T$ by the choice of $\rho_{z_j}$ in Section \ref{section 3}, these lemmas provide a constant $c=c(\operatorname{data})$ such that
$$
\iints{\kappa \mathcal{Q}_j} \Psi(z,|Du|)\, dz \leq c \Lambda^{1-\theta}\iints{\mathcal{Q}_j\cap \Upsilon(c^{-1}\Lambda)} \Psi^\theta(z,|Du|)\, dz
$$
for every $j\in\mathbb{N}$. Since almost every point of $\Upsilon(\Lambda,r_1)$ is such a Lebesgue point $w$ and $w\in\mathcal{Q}(w)\subset\kappa\mathcal{Q}_j$ for some $j$, summing this inequality over $j$ and using that the cylinders $\mathcal{Q}_j\subset Q_{r_2}(z_0)$ are pairwise disjoint, we obtain 
$$
\iints{\Upsilon(\Lambda,r_1)} \Psi(z,|Du|)\, dz\leq c \Lambda^{1-\theta}\iints{\Upsilon(c^{-1}\Lambda,r_2)} \Psi^\theta(z,|Du|)\, dz.
$$
Then, in view of 
$$
\iints{\Upsilon(c^{-1}\Lambda,r_1)\setminus \Upsilon(\Lambda,r_1)} \Psi(z,|Du|)\, dz\leq \Lambda^{1-\theta}\iints{\Upsilon(c^{-1}\Lambda,r_2)} \Psi^\theta(z,|Du|)\, dz,
$$
we obtain
$$
\iints{\Upsilon(c^{-1}\Lambda,r_1)} \Psi(z,|Du|)\, dz\leq c \Lambda^{1-\theta}\iints{\Upsilon(c^{-1}\Lambda,r_2)} \Psi^\theta(z,|Du|)\, dz.
$$
For $k\in\mathbb{N}$, let
$$
\Psi_k(z,|Du|) = \min\{\Psi(z,|Du|),k\}
$$
and
$$
\Upsilon_k(\Lambda,\rho)=\{z\in Q_\rho(z_0): \Psi_k(z,|Du(z)|)>\Lambda\}.
$$
We observe that if $\Lambda\geq k$, then $\Upsilon_k(\Lambda,\rho)=\emptyset$ and if $\Lambda< k$, then $\Upsilon_k(\Lambda,\rho)=\Upsilon(\Lambda,\rho)$. This yields 
$$
\iints{\Upsilon_k(c^{-1}\Lambda,r_1)} \Psi^{1-\theta}_k(z,|Du|)\Psi^\theta(z,|Du|)\, dz\leq c \Lambda^{1-\theta}\iints{\Upsilon_k(c^{-1}\Lambda,r_2)} \Psi^\theta(z,|Du|)\, dz.
$$
In view of the lower bound on $\Lambda$ imposed in Section \ref{section 3}, we put
$$
\Lambda_1= c^{-1}\left(\frac{4\kappa r}{r_2-r_1}\right)^{\frac{\sigma_H(n+2)}{2}}\Lambda_0.
$$
Then for any $\Lambda >\Lambda_1$, we have
$$
\iints{\Upsilon_k(\Lambda,r_1)} \Psi^{1-\theta}_k(z,|Du|)\Psi^\theta(z,|Du|)\, dz\leq c \Lambda^{1-\theta}\iints{\Upsilon_k(\Lambda,r_2)} \Psi^\theta(z,|Du|)\, dz.
$$
Let $\varepsilon\in(0,1)$ be a parameter to be fixed later. Multiplying the above inequality by $\Lambda^{\varepsilon - 1}$ and integrating over $(\Lambda_1, \infty)$, we obtain 
$$
\begin{aligned}
    \mathrm{I}&=\int_{\Lambda_1}^\infty \Lambda^{\varepsilon-1}\iints{\Upsilon_k(\Lambda,r_1)} \Psi^{1-\theta}_k(z,|Du|)\Psi^\theta(z,|Du|)\, dz\, d\Lambda\\
    &\leq c \int_{\Lambda_1}^\infty \Lambda^{\varepsilon-\theta}\iints{\Upsilon_k(\Lambda,r_2)} \Psi^\theta(z,|Du|)\, dz\,d\Lambda=\mathrm{II}.
\end{aligned}
$$
Then, Fubini's theorem gives
$$
\begin{aligned}
\mathrm{I}&=\frac{1}{\varepsilon}\iints{\Upsilon_k(\Lambda_1,r_1)} \Psi^{1-\theta+\varepsilon}_k
(z,|Du|)\Psi^\theta(z,|Du|)\, dz\\
&\quad- \frac{1}{\varepsilon}\Lambda_1^\varepsilon\iints{\Upsilon_k(\Lambda_1,r_1)} \Psi^{1-\theta}_k
(z,|Du|)\Psi^\theta(z,|Du|)\, dz
\end{aligned}
$$
and
$$
\mathrm{II}\leq \frac{c}{1-\theta+\varepsilon}\iints{Q_{r_2}(z_0)} \Psi^{1-\theta+\varepsilon}_k
(z,|Du|)\Psi^\theta(z,|Du|)\, dz.
$$
Moreover, since
$$
\begin{aligned}
&\iints{Q_{r_1}(z_0)\setminus \Upsilon_k(\Lambda_1,r_1)} \Psi^{1-\theta+\varepsilon}_k
(z,|Du|)\Psi^\theta(z,|Du|)\, dz\\
&\qquad \leq \Lambda_1^\varepsilon \iints{Q_{2r}(z_0)} \Psi^{1-\theta}_k
(z,|Du|)\Psi^\theta(z,|Du|)\, dz,
\end{aligned}
$$
we get
$$
\begin{aligned}
\mathrm{I}&\geq\frac{1}{\varepsilon}\iints{Q_{r_1}(z_0)} \Psi^{1-\theta+\varepsilon}_k
(z,|Du|)\Psi^\theta(z,|Du|)\, dz\\
&\quad- \frac{2}{\varepsilon}\Lambda_1^\varepsilon\iints{Q_{2r}(z_0)} \Psi^{1-\theta}_k
(z,|Du|)\Psi^\theta(z,|Du|)\, dz.
\end{aligned}
$$
Combining these estimates, we have
$$
\begin{aligned}
    &\iints{Q_{r_1}(z_0)} \Psi^{1-\theta+\varepsilon}_k
    (z,|Du|)\Psi^\theta(z,|Du|)\, dz\\
    &\qquad\leq \frac{c\varepsilon}{1-\theta+\varepsilon}\iints{Q_{r_2}(z_0)} \Psi^{1-\theta+\varepsilon}_k
    (z,|Du|)\Psi^\theta(z,|Du|)\, dz\\
    &\qquad\qquad +c\Lambda_1^\varepsilon\iints{Q_{2r}(z_0)} \Psi^{1-\theta}_k
(z,|Du|)\Psi^\theta(z,|Du|)\, dz.
\end{aligned}
$$
We choose $\varepsilon_0=\varepsilon_0(\operatorname{data})\in(0,1)$ so that for any $\varepsilon\in (0,\varepsilon_0)$,
$$
\frac{c\varepsilon}{1-\theta+\varepsilon}\leq \frac{1}{2}.
$$
Set $\gamma\coloneq\frac{\sigma_H(n+2)}{2}$. Since $\frac{4\kappa r}{r_2-r_1}\geq 4\kappa>1$ and $\varepsilon<1$, we have $\Lambda_1^\varepsilon\leq c\left(\frac{4\kappa r}{r_2-r_1}\right)^{\gamma}\Lambda_0^\varepsilon$ with $c$ independent of $\varepsilon$. We apply Lemma \ref{lem : iteration lemma} on $[r,2r]$ with $\vartheta=\frac12$, $B=0$, the bounded function
$$
g(\varrho)\coloneq\iints{Q_\varrho(z_0)}\Psi^{1-\theta+\varepsilon}_k(z,|Du|)\Psi^\theta(z,|Du|)\, dz
$$
and the constant
$$
A\coloneq c(4\kappa r)^{\gamma}\Lambda_0^\varepsilon\iints{Q_{2r}(z_0)} \Psi^{1-\theta}_k(z,|Du|)\Psi^\theta(z,|Du|)\, dz,
$$
and we get
$$
\begin{aligned}
    \iints{Q_r(z_0)} \Psi^{1-\theta+\varepsilon}_k
    (z,|Du|)\Psi^\theta(z,|Du|)\, dz\leq c\Lambda_0^\varepsilon\iints{Q_{2r}(z_0)} \Psi^{1-\theta}_k
(z,|Du|)\Psi^\theta(z,|Du|)\, dz.
\end{aligned}
$$
Letting $k\rightarrow\infty$, we have
$$
    \iints{Q_r(z_0)} \Psi^{1+\varepsilon}
    (z,|Du|)\, dz\leq c\Lambda_0^\varepsilon\iints{Q_{2r}(z_0)} \Psi(z,|Du|)\, dz.
$$
Recalling \eqref{def : lambda_0 and Lambda_0}, since $G\prec H$ and $\lambda_0\geq 1$, Lemma \ref{lem205} gives $G(\lambda_0)\leq H(\lambda_0)$ and hence $\Lambda_0^\varepsilon\leq cH^\varepsilon(\lambda_0)$ for some $c=c(\operatorname{data},\|a\|_{L^\infty(\Omega_T)})$. Moreover, by \eqref{cond : monotonicity of Phi over powers}, 
$$
\Lambda_0^\varepsilon\leq c\lambda_0^{\varepsilon\sigma_H}.
$$
Since \eqref{def : lambda_0 and Lambda_0} gives
$$
\lambda_0^{\varepsilon\sigma_H}\leq c\left[\left(\miint{Q_{2r}(z_0)}\Psi(z,|Du|)\, dz\right)^{\frac{\varepsilon\sigma_H}{2}}+1\right],
$$
dividing by $|Q_r|$ and using $x\leq x^{1+\frac{\varepsilon\sigma_H}{2}}+1$ for $x\geq 0$, we obtain the desired conclusion when $2r\leq 1$.

Finally, we remove the restriction $2r\leq 1$. Let $Q_{2r}(z_0)\subset\Omega_T$ with $2r>1$. Since $B_{2r}(x_0)\subset\Omega$, we have $r\leq \frac{1}{4}\operatorname{diam}(\Omega)$. We choose points $\hat{z}_1,\dots,\hat{z}_M\in Q_r(z_0)$ with $M\leq c(n)\operatorname{diam}(\Omega)^{n+2}$ such that $Q_r(z_0)\subset\bigcup_{i=1}^M Q_{1/4}(\hat{z}_i)$. Since $r\geq\frac{1}{2}$, we have $Q_{1/2}(\hat{z}_i)\subset Q_{r+\frac{1}{2}}(z_0)\subset Q_{2r}(z_0)$ for every $i$. Applying the estimate already proved to the cylinders $Q_{1/2}(\hat{z}_i)$, we obtain
$$
\begin{aligned}
&\iints{Q_r(z_0)}\Psi(z,|Du|)^{1+\varepsilon}\, dz\\
&\quad\leq \sum_{i=1}^M\iints{Q_{1/4}(\hat{z}_i)}\Psi(z,|Du|)^{1+\varepsilon}\, dz\\
&\quad\leq c|Q_{1/4}|\sum_{i=1}^M\left[\left(\miint{Q_{1/2}(\hat{z}_i)}\Psi(z,|Du|)\, dz\right)^{1+\frac{\varepsilon\sigma_H}{2}}+1\right]\\
&\quad\leq cM|Q_{1/4}|\left[\left(\frac{|Q_{2r}|}{|Q_{1/2}|}\miint{Q_{2r}(z_0)}\Psi(z,|Du|)\, dz\right)^{1+\frac{\varepsilon\sigma_H}{2}}+1\right].
\end{aligned}
$$
Since $\frac{|Q_{2r}|}{|Q_{1/2}|}=(4r)^{n+2}\leq \operatorname{diam}(\Omega)^{n+2}$, $\varepsilon\leq 1$ and $M|Q_{1/4}|\leq c(n)\operatorname{diam}(\Omega)^{n+2}|Q_r|$, dividing by $|Q_r|$ gives the desired conclusion with a constant depending on $\operatorname{data}$, $\operatorname{diam}(\Omega)$ and $\|a\|_{L^\infty(\Omega_T)}$. \qed


\begin{thebibliography}{10}

\bibitem{BaasandorjByun2023ASNS}
S.~Baasandorj and S.-S. Byun.
\newblock Gradient estimates for {O}rlicz double-phase problems.
\newblock {\em Ann. Sc. Norm. Super. Pisa Cl. Sci. (5)}, 24:2215--2268, 2023.

\bibitem{Baasandorj2020}
S.~Baasandorj, S.-S. Byun, and J.~Oh.
\newblock Calder\'{o}n-{Z}ygmund estimates for generalized double phase
  problems.
\newblock {\em J. Funct. Anal.}, 279(7):108670, 57, 2020.

\bibitem{BaroniColomboMingione2016SPMJ}
P.~Baroni, M.~Colombo, and G.~Mingione.
\newblock Non-autonomous functionals, borderline cases and related function
  classes.
\newblock {\em St. Petersburg Math. J.}, 27(3):347--379, 2016.

\bibitem{Baroni2018}
P.~Baroni, M.~Colombo, and G.~Mingione.
\newblock Regularity for general functionals with double phase.
\newblock {\em Calc. Var. Partial Differential Equations}, 57(2):Paper No. 62,
  48, 2018.

\bibitem{BogeleinDuzaarMarcellini2013ARMA}
V.~B\"{o}gelein, F.~Duzaar, and P.~Marcellini.
\newblock Parabolic systems with {$p,q$}-growth: a variational approach.
\newblock {\em Arch. Ration. Mech. Anal.}, 210(1):219--267, 2013.

\bibitem{Byun2020}
S.-S. Byun and J.~Oh.
\newblock Regularity results for generalized double phase functionals.
\newblock {\em Anal. PDE}, 13(5):1269--1300, 2020.

\bibitem{Colombo2015a}
M.~Colombo and G.~Mingione.
\newblock Bounded minimisers of double phase variational integrals.
\newblock {\em Arch. Ration. Mech. Anal.}, 218(1):219--273, 2015.

\bibitem{Colombo2015}
M.~Colombo and G.~Mingione.
\newblock Regularity for double phase variational problems.
\newblock {\em Arch. Ration. Mech. Anal.}, 215(2):443--496, 2015.

\bibitem{DeFilippisMingione2021AMP}
C.~De~Filippis and G.~Mingione.
\newblock Interpolative gap bounds for nonautonomous integrals.
\newblock {\em Anal. Math. Phys.}, 11(3):Paper No. 117, 2021.

\bibitem{1993_Degenerate_parabolic_equations_DiBenedetto}
E.~DiBenedetto.
\newblock {\em Degenerate parabolic equations}.
\newblock Universitext. Springer-Verlag, New York, 1993.

\bibitem{2011Diening}
L.~Diening, P.~Harjulehto, P.~H\"{a}st\"{o}, and M.~R\r{u}\v{z}i\v{c}ka.
\newblock {\em Lebesgue and {S}obolev spaces with variable exponents}, volume
  2017 of {\em Lecture Notes in Mathematics}.
\newblock Springer, Heidelberg, 2011.

\bibitem{Esposito2004}
L.~Esposito, F.~Leonetti, and G.~Mingione.
\newblock Sharp regularity for functionals with {$(p,q)$} growth.
\newblock {\em J. Differential Equations}, 204(1):5--55, 2004.

\bibitem{Gehring1973}
F.~W. Gehring.
\newblock The {$L^{p}$}-integrability of the partial derivatives of a
  quasiconformal mapping.
\newblock {\em Acta Math.}, 130:265--277, 1973.

\bibitem{GiaquintaModica1979}
M.~Giaquinta and G.~Modica.
\newblock Regularity results for some classes of higher order nonlinear
  elliptic systems.
\newblock {\em J. Reine Angew. Math.}, 311/312:145--169, 1979.

\bibitem{2003_Giusti_Direct_methods_in_the_calculus_of_variations}
E.~Giusti.
\newblock {\em Direct methods in the calculus of variations}.
\newblock World Scientific Publishing Co., Inc., River Edge, NJ, 2003.

\bibitem{Harjulehto2019}
P.~Harjulehto and P.~H\"ast\"o.
\newblock {\em Orlicz spaces and generalized {O}rlicz spaces}, volume 2236 of
  {\em Lecture Notes in Mathematics}.
\newblock Springer, Cham, 2019.

\bibitem{Hasto_2021}
P.~H\"{a}st\"{o} and J.~Ok.
\newblock Higher integrability for parabolic systems with {O}rlicz growth.
\newblock {\em J. Differential Equations}, 300:925--948, 2021.

\bibitem{Haestoe2022a}
P.~H\"{a}st\"{o} and J.~Ok.
\newblock Maximal regularity for local minimizers of non-autonomous
  functionals.
\newblock {\em J. Eur. Math. Soc. (JEMS)}, 24(4):1285--1334, 2022.

\bibitem{HastoOk2026}
P.~H{\"a}st{\"o} and J.~Ok.
\newblock Higher integrability for parabolic {PDE}s with generalized {O}rlicz
  growth.
\newblock {\em arXiv:2511.19758, to appear in JMPA}, 2026.

\bibitem{KimOh2024JDE}
B.~Kim and J.~Oh.
\newblock Higher integrability for weak solutions to parabolic multi-phase
  equations.
\newblock {\em J. Differential Equations}, 409:223--298, 2024.

\bibitem{kim2025bounded}
B.~Kim and J.~Oh.
\newblock Bounded solutions and interpolative gap bounds for degenerate
  parabolic double phase problems.
\newblock {\em arXiv preprint arXiv:2511.13454}, 2025.

\bibitem{kim2026interpolative}
B.~Kim and J.~Oh.
\newblock Interpolative refinement of gap bound conditions for singular
  parabolic double phase problems.
\newblock {\em arXiv preprint arXiv:2601.01571}, 2026.

\bibitem{2023_Gradient_Higher_Integrability_for_Degenerate_Parabolic_Double-Phase_Systems}
W.~Kim, J.~Kinnunen, and K.~Moring.
\newblock Gradient higher integrability for degenerate parabolic double-phase
  systems.
\newblock {\em Arch. Ration. Mech. Anal.}, 247(5):Paper No. 79, 46, 2023.

\bibitem{KIM2025110738}
W.~Kim, J.~Kinnunen, and L.~Särkiö.
\newblock Lipschitz truncation method for parabolic double-phase systems and
  applications.
\newblock {\em J. Funct. Anal.}, 288(3):Paper No. 110738, 2025.

\bibitem{KIM2025113231}
W.~Kim, K.~Moring, and L.~Särkiö.
\newblock Hölder regularity for degenerate parabolic double-phase equations.
\newblock {\em J. Differential Equations}, 434:Paper No. 113231, 2025.

\bibitem{KS24}
W.~Kim and L.~S\"arki\"o.
\newblock Gradient higher integrability for singular parabolic double-phase
  systems.
\newblock {\em NoDEA Nonlinear Differential Equations Appl.}, 31(3):Paper No.
  40, 38, 2024.

\bibitem{Kinnunen2000}
J.~Kinnunen and J.~L. Lewis.
\newblock Higher integrability for parabolic systems of {$p$}-{L}aplacian type.
\newblock {\em Duke Math. J.}, 102(2):253--271, 2000.

\bibitem{1961Rutickii}
M.~A. Krasnosel'ski\v{i} and J.~B. Ruticki\v{i}.
\newblock {\em Convex functions and {O}rlicz spaces}.
\newblock P. Noordhoff Ltd., Groningen, 1961.

\bibitem{Marcellini1989}
P.~Marcellini.
\newblock Regularity of minimizers of integrals of the calculus of variations
  with non standard growth conditions.
\newblock {\em Arch. Ration. Mech. Anal.}, 105:267--284, 1989.

\bibitem{OhSen2026CVPDE}
J.~Oh and A.~Sen.
\newblock Gradient higher integrability for degenerate parabolic double phase
  systems with two modulating coefficients.
\newblock {\em Calc. Var. Partial Differential Equations}, 65(6):Paper No. 187,
  2026.

\bibitem{Ok2020}
J.~Ok.
\newblock Regularity for double phase problems under additional integrability
  assumptions.
\newblock {\em Nonlinear Anal.}, 194:111408, 13, 2020.

\bibitem{OSS24}
J.~Ok, G.~Scilla, and B.~Stroffolini.
\newblock Regularity theory for parabolic systems with {U}hlenbeck structure.
\newblock {\em J. Math. Pures Appl. (9)}, 182:116--163, 2024.

\bibitem{S25}
A.~Sen.
\newblock Gradient higher integrability for degenerate/singular parabolic
  multi-phase problems.
\newblock {\em J. Geom. Anal.}, 35(6):Paper No. 170, 95, 2025.

\bibitem{Zhikov1986}
V.~V. Zhikov.
\newblock Averaging of functionals of the calculus of variations and elasticity
  theory.
\newblock {\em Izv. Akad. Nauk SSSR Ser. Mat.}, 50(4):675--710, 877, 1986.

\bibitem{Zhikov1995}
V.~V. Zhikov.
\newblock On {L}avrentiev's phenomenon.
\newblock {\em Russian J. Math. Phys.}, 3(2):249--269, 1995.

\end{thebibliography}
\end{document}